\documentclass[12pt]{report}
\usepackage[letterpaper, margin = 1.1in]{geometry}
\usepackage{amsmath}
\usepackage{amsthm}
\usepackage{amssymb}
\usepackage{array}
\usepackage{makecell}
\usepackage{mathtools}
\usepackage{paracol}

\usepackage{enumerate}
\usepackage{enumitem}
\usepackage{tikz-cd}
\usepackage{parskip}
\usepackage{graphicx}
\usepackage{boxedminipage}
\graphicspath{images/}

\usepackage[
backend=biber,
style=alphabetic,
sorting=ynt
]{biblatex}
\usepackage{hyperref}
\usepackage[nameinlink,capitalize]{cleveref}

\newcommand{\sbb}[2]{\left( \frac{#1}{#2} \right)}
\newcommand{\sbx}[4]{\left(\frac{#1}{#2},\frac{#3}{#4}\right)}
\renewcommand{\hom}{\text{Hom}}
\newcommand{\free}{\text{free}}
\newcommand{\sett}{\textbf{Set}}

\theoremstyle{definition}
\newtheorem*{defn}{Definition}
\newtheorem*{defns}{Definitions}

\newtheorem{uass}{Assumption}[chapter]

\newtheorem{prop}{Proposition}[chapter]

\newtheorem{clm}[prop]{Claim}
\newtheorem{rk}[prop]{Remark}
\newtheorem{lem}[prop]{Lemma}
\newcommand{\pf}{{\it Proof. }}
\newtheorem{cor}[prop]{Corollary}
\newtheorem{obs}[prop]{Observation}

\newcommand{\qo}{$\qed_{\text{F1}}$}
\newcommand{\qtw}{$\qed_{\text{F2}}$}
\newcommand{\qth}{$\qed_{\text{F3}}$}
\newcommand{\qedone}{$\qed_{\text{Case 1}}$}
\newcommand{\qedo}{$\qed_{\text{Case 0}}$}
\newcommand{\qcl}{$\qed_{\text{claim}}$}

\newcommand{\E}{\text{E}}
\newcommand{\RHS}{\text{RHS}}
\newcommand{\LHS}{\text{LHS}}

\newcommand{\fone}{\text{F1}}
\newcommand{\ftwo}{\text{F2}}
\newcommand{\fth}{\text{F3}}
\newcommand{\mfone}{\text{mF1}}
\newcommand{\mftwo}{\text{mF2}}
\newcommand{\mfth}{\text{mF3}}
\newcommand{\im}{\text{Im }}
\newcommand{\prehom}{\text{Prehom} }
\newcommand{\FP}{\text{FP}}
\newcommand{\PD}{\text{P}}
\newcommand{\pr}{\ensuremath{^{\prime}}}
\newcommand{\sub}{\text{Sub}}
\newcommand{\syn}{\text{Syn}}

\newcommand{\tkref}{\text{TKREF}}

\newcommand{\ins}{\text{INS}\ensuremath{\bigwedge} }
\newcommand{\ex}{\text{EX}\ensuremath{\bigwedge} }
\newcommand{\reord}{\text{REORD} }
\newcommand{\frlS}{\ensuremath{\forall}\text{S} }
\newcommand{\frlA}{\ensuremath{\forall}\text{A} }
\newcommand{\eS}{\ensuremath{\exists}\text{S} }
\newcommand{\eA}{\ensuremath{\exists}\text{A} }
\newcommand{\eD}{\ensuremath{\exists}\text{D} }
\newcommand{\id}{\text{ID}}
\newcommand{\ctr}{\text{CTR}}
\newcommand{\cp}{\text{CP}}
\newcommand{\expl}{\text{EXP}}

\newcommand{\pc}{\text{PC} }
\newcommand{\frelim}{\ensuremath{\forall}\text{ELIM }}
\newcommand{\eelim}{\ensuremath{\exists}\text{ELIM }}
\newcommand{\meta}{\text{META} }
\newcommand{\sdrop}{\text{SDROP} }
\newcommand{\andd}{\text{D}\ensuremath{\bigwedge} }
\newcommand{\oorl}{\text{L}\ensuremath{\bigvee}}
\newcommand{\oorr}{\text{R}\ensuremath{\bigvee}}
\newcommand{\andrr}{\text{R}\ensuremath{\bigwedge} }

\newcommand{\lo}{\text{LOAD} }
\newcommand{\un}{\text{UNLOAD} }
\newcommand{\ch}{\text{CHAIN} }
\newcommand{\assm}{\text{ASSM} }
\newcommand{\subs}{\text{SUBS} }
\newcommand{\ant}{\text{ANT} }
\newcommand{\sep}{\text{SEP}\ensuremath{\bigvee}}
\newcommand{\unsep}{\text{UNSEP}\ensuremath{\bigvee}}

\newcommand{\DeclareAutoPairedDelimiter}[3]{%
  \expandafter\DeclarePairedDelimiter\csname Auto\string#1\endcsname{#2}{#3}%
  \begingroup\edef\x{\endgroup
    \noexpand\DeclareRobustCommand{\noexpand#1}{%
      \expandafter\noexpand\csname Auto\string#1\endcsname*}}%
  \x}
\DeclareAutoPairedDelimiter{\paren}{(}{)}
\DeclareAutoPairedDelimiter{\abs}{\vert}{\vert}
\DeclareAutoPairedDelimiter{\rak}{[}{]}
\DeclareAutoPairedDelimiter{\crul}{\{}{\}}

\title{A Syntactic and Categorical Derivation of G\"{o}del's Completeness Theorem}
\author{Hugo A. Jenkins}
\date{September 2021}

\usepackage{subfiles} 

\begin{document}

\begin{titlepage}
    \begin{center}
        \vspace*{1cm}
        
        \Huge
        \textbf{A Syntactic and Categorical Derivation of G\"{o}del's Completeness Theorem}  
        
        \vspace{0.5cm}
        \LARGE

        \vspace{1.5cm}
        
        \textbf{Hugo Jenkins}
        
        \vfill
        
        \vspace{0.8cm}
        
        \Large
        \copyright \hspace{2mm} Hugo Jenkins September 2021
        
    \end{center}
\end{titlepage}

\begin{abstract}
\thispagestyle{plain}
\newcommand{\calm}{\mathcal{M}}
\newcommand{\set}{\textbf{Set}}
We give a ``fully syntactic'' construction of the category $\syn(T)$\footnote{Referred to in our primary reference \cite{lu} as the {\it weak syntactic category} $\text{Syn}_{0}(T)$.} associated to a consistent theory $T$ over a pure predicate language in first-order logic with equality; show it is a consistent coherent category; and show that a morphism of coherent categories from $\syn(T)\to \set$ gives rise to a model $\calm$ of $T$ in the usual sense. Thus, we show a main case	 of the classical G\"{o}del Completeness Theorem via Deligne's Completeness Theorem for coherent categories, as outlined in \cite{lu}, Lecture 6.

See \cite{di} (in French) for similar work to ours but with a broader context.

\end{abstract}

\chapter*{Dedication}
This work is dedicated to Siri E. Huntoon, my mother. Without her it would not have come to be.

\chapter*{Acknowledgements}
My utmost thanks to Professor Thomas W. Scanlon, for his excellent mentorship and patient guidance. Conversation with and encouragement from him has allowed this work to emerge. Thanks to my mother, and to the texnical support dream team: Nicholas and Owen Jenkins. Thanks also to the librarian at the Universit\'{e} de Montr\'{e}al who was kind enough to send me a photocopy of almost 80 pages of \cite{di}.

\tableofcontents

\chapter{Introduction and Deductive Calculus}\label{ch-intro}

\section{Language}
For the entirety of this paper, we will be considering a fixed first-order language $L$ with equality (see e.g. \cite{eft}) {\it over a signature which contains only predicate (relation) symbols}. Thus, we have no nontrivial terms. The variables in our language will be considered to be indexed by the {\it integers $\mathbb{Z}$}, and will be written as $x_{k}$ for $k\in \mathbb{Z}$. We will use capital Roman letters for metavariables, i.e. placeholders for some partially or fully general $x_{i}$, as the context may be.\\

$E$ will refer to the ``emptiness sentence'' $\neg \exists x_{0}\rak{x_{0} = x_{0}}$.

\section{Sequents}
We will formalize provability via the {\it one-sided sequent calculus}, a system due to Gentzen 
in which the basic objects are of the form\[
\Gamma \Rightarrow \varphi,
\]
where $\Gamma$ is a finite set of w.f.f.s and $\varphi$ is a w.f.f. These objects are called {\it sequents}. They are defined by their left-hand-sides and right-hand-sides, and there are no restrictions on what these can be (i.e. no further restrictions on the finite set $\Gamma$ or the formula $\varphi$.)

We call a sequent {\it derivable} 
if there exists a finite derivation sequence for it in a certain formal calculus, the rules of which are given below. A derivation sequence is essentially a sequence of these rules, each term augmented with the data of the premises and conclusion, and perhaps some other data witnessing a legal application.

When we speak of obtaining/showing/deriving a sequent, we mean giving a mathematical argument in the metalanguage for why it is derivable.\\

Often we will explicitly indicate certain elements or subsets of the LHS by writing them separated by commas, e.g.\[
x_{1} = x_{2}, \exists x_{1} \rak{x_{1} = x_{1}}, Px_{1}x_{2} \Rightarrow \varphi
\]or\[
\exists x_{1} \rak{x_{1} = x_{1}}, x_{1} = x_{2}, Px_{1}x_{2} \Rightarrow \varphi
\]or
\[\crul{x_{1} = x_{2}, \exists x_{1} \rak{x_{1} = x_{1}}}, Px_{1}x_{2}, Px_{1}x_{2} \Rightarrow \varphi.
\]
The LHS is always considered as a collection of formulas in the natural way and so the above all refer to the same sequent.

For a more general use of sequents in categorical logic, see \cite{mr}. For some context on the goal of the present paper, see especially the remarks on pages 241-242.

\section{Substitution and alpha-equivalence}\label{seg-subalph}
The notion of changing every free occurrence of $x_{j}$ in a formula $\varphi$ to an occurrence of $x_{i}$, for arbitrary $i$, $j$, and $\varphi$, will be essential to us. However we want a substitution operation, like the one defined in \cite{eft}, which makes ``variable capture'' impossible. Thus we wish to define $\varphi\sbb{x_{i}}{x_{j}}$ in such a way that, for example, if $\varphi$ is not a sentence then $\varphi\sbb{x_{i}}{x_{j}}$ also cannot be. For this we could indeed take $\sbb{\cdot}{\cdot}$ to be defined as in \cite{eft} (written there without the parentheses). We note that this definition, as far as we can see, is not technically justified using the set-theoretic Recursion Theorem stated in \cite{end}.\\

However, we in fact want something slightly more than just safe substitution. In the remainder of this paper, we will be (semi-secretly) considering each w.f.f. $\varphi$ only up to {\it alpha-equivalence}. This is an equivalence relation on w.f.f.s which identifies, for example, the {\it alphabetic variants} $\exists x_{0}\rak{x_{1} = x_{0}}$ and $\exists x_{7}\rak{x_{1} = x_{7}}$; or $\exists x_{1} \exists x_{2} \rak{x_{1} = x_{2}}$ and $\exists x_{2} \exists x_{1} \rak{x_{2} = x_{1}}$; but it does not identify these latter with $\exists x_{2} \exists x_{1}\rak{ x_{1} = x_{2}}$. The bound variables of an alpha-equivalence class are not well-defined. A precise and satisfactory construction of the alpha-equivalence relation using the Recursion Theorem (allowing, for instance, a verification that it respects formula constructors and that the sequent calculus is well-defined mod it) is, in the author's opinion, quite nontrivial and in fact occupied approximately the first third of his time spent on this project. It will not be given here as it is of no interest.\\

Instead, we will just assume the necessary properties of alpha-equivalence, and assume we have a substitution $\varphi\sbb{x_{i}}{x_{j}}$ which is defined by first taking an alphabetic variant $\varphi'$ in which neither of $x_{i}$, $x_{j}$ occur bound (which certainly exists; cf \cite{end} Theorem 24I for a proof of almost this statement), and then performing ``raw substitution'' on the level of strings, changing {\it every} $x_{j}$ in $\varphi'$ to $x_{i}$. Two key properties of this setup we will need are that \begin{enumerate}
\item If $T$ is distinct from $X$ and $Y$, then
\[
(\exists T \varphi)\sbb{Y}{X}\quad \text{and} \quad \exists T \rak{\varphi\sbb{Y}{X}} \quad \text{are alpha-equivalent.}
\]
\item If $N$ does not occur at all (free or bound) in $\varphi$, then
\[
\exists X \varphi \quad \text{and} \quad \exists N \rak{\varphi\sbb{N}{X}} \quad \text{are alpha-equivalent.}
\]
\end{enumerate}

From now on we speak as little about alpha-equivalence as possible (almost not at all outside this chapter). See \cite{mr}, p. 70 for a brief discussion which differs from ours but also essentially states property 2 above.\\

We shall need the following properties of $\sbb{\cdot}{\cdot}$, all of which hold essentially because they hold for raw substitution on strings:\begin{enumerate}
\item $\varphi\sbb{X}{X} = \varphi$
\item If $N$ does not occur in $\varphi$, then $\varphi\sbb{N}{X}\sbb{Y}{N} = \varphi\sbb{Y}{X}$
\item If $X$, $Y$, $Z$ are pairwise distinct and $X$, $Y$, $T$ are pairwise distinct, then $\varphi\sbb{Z}{X}\sbb{T}{Y} = \varphi\sbb{T}{Y}\sbb{Z}{X}$.
\end{enumerate}

\subsection{Simultaneous substitution}\label{seg-simsub}
Using the above, we can define the {\it simultaneous substitution} $\varphi \sbb{Y_{i}}{X_{i}}$, where $(Y_{i})$ and $(X_{i})$ are variable {\it sequences} of equal length $n\geq 1$, and the $X_{i}$ are distinct. Namely, we define it as \[
\varphi\sbb{N_{1}}{X_{1}}\dots \sbb{N_{n}}{X_{n}}\sbb{Y_{1}}{N_{1}}\dots \sbb{Y_{n}}{N_{n}},
\]
where the $N_{i}$ are any distinct variables not occurring in $\varphi$ and not equal to any $X$'s or $Y$'s. We must check that this is well-defined.

But indeed, by substitution property 3 we can commute the first $n$ substitutions amongst themselves, and also the last $n$. We can then see that for each $i$ the iterated substitution is not dependent on the choice of $N_{i}$, using substitution property 2. For instance, to show independence from $N_{1}$, we rewrite as\[
\varphi \dots \sbb{N_{n}}{X_{n}}\sbb{N_{1}}{X_{1}}\sbb{Y_{1}}{N_{1}}\dots \sbb{Y_{n}}{N_{n}}
\]\[
 = \varphi \dots \sbb{N_{n}}{X_{n}} \sbb{Y_{1}}{X_{1}}\dots \sbb{Y_{n}}{N_{n}}, \tag{1}
\]
by property 2, justified as $N_{1}$ does not occur in the formula $\varphi\dots \sbb{N_{n}}{X_{n}}$ which is being substituted into. Now $N_{1}$ no longer appears in the expression (1).

It follows that $\varphi\sbb{Y_{i}}{X_{i}}$ is well-defined.\qed\\

We shall need the following not-entirely-obvious property of simultaneous substitution:\begin{enumerate}
\item If $(X_{i})$ and $(T_{i})$ are sequences of distinct variables such that all free variables of $\varphi$ occur among the $X$'s, and $(Y_{i})$ is any sequence, then\[
\varphi\sbb{T_{i}}{X_{i}}\sbb{Y_{i}}{T_{i}} = \varphi \sbb{Y_{i}}{X_{i}}.
\]\end{enumerate}
\pf Let $(N_{i})$ be a sequence of distinct variables not in $\varphi$ and disjoint from all sequences above. We can write the left-hand side as\[
\varphi\sbb{N_{1}}{X_{1}}\dots \sbb{N_{n}}{X_{n}}\sbb{T_{1}}{N_{1}}\dots \sbb{T_{n}}{N_{n}}\sbb{N_{1}}{T_{1}}\dots \sbb{N_{n}}{T_{n}}\sbb{Y_{1}}{N_{1}}\dots \sbb{Y_{1}}{N_{n}} \tag{2}
\]
since the $N$'s do not occur in $\varphi \sbb{T_{i}}{X_{i}} = \varphi\sbb{N_{1}}{X_{1}}\dots \sbb{N_{n}}{X_{n}}\sbb{T_{1}}{N_{1}}\dots \sbb{T_{n}}{N_{n}}$. Commuting $\sbb{N_{n}}{T_{n}}$ leftwards $n-1$ times by substitution property 3, we rewrite (2) so that following the second ellipsis we have the pair \[
\sbb{T_{n}}{N_{n}}\sbb{N_{n}}{T_{n}}.\]
We see that $T_{n}$ does not occur in the formula to the left of this pair: the formula $\varphi\sbb{N_{1}}{X_{1}}\dots \sbb{N_{n}}{X_{n}}$ does not contain it, having only free variables among the $N$'s; and it is not introduced by any of the intervening substitutions. Therefore by substitution property 2 the pair acts as just \[
\sbb{N_{n}}{N_{n}} \qquad = \text{id, by subs. property 1.}
\]
Therefore we have that (2) is equal to\[
\varphi\sbb{N_{1}}{X_{1}}\dots \sbb{N_{n}}{X_{n}}\sbb{T_{1}}{N_{1}}\dots \sbb{T_{n-1}}{N_{n-1}}\sbb{N_{1}}{T_{1}}\dots \sbb{N_{n-1}}{T_{n-1}}\sbb{Y_{1}}{N_{1}}\dots \sbb{Y_{1}}{N_{n}}.
\]
Iterating this we get rid of the middle $2n$ substitutions in (2), and we get something equal to the right-hand side by definition. \qed\\

If $(Y_{i})$, $(Y\pr_{j})$, $(X_{i})$, $(X\pr_{j})$ are sequences such that the $X$'s and $X\pr$'s are all mutually distinct, we shall sometimes form a simultaneous substitution using the concatenation of the $Y$'s on top, and that of the $X$'s on bottom. We will write this as $\varphi\sbx{Y_{i}}{X_{i}}{Y\pr_{j}}{X\pr_{j}}$.

\section{Sequent calculus}
\subsection{Rules}
\subsubsection{Premise-free rules}
\begin{center}
\begin{boxedminipage}{0.35\textwidth}
\[\overline{\Gamma \Rightarrow \varphi} \tag{\assm}\]
if $\varphi \in \Gamma$.\end{boxedminipage}\qquad \qquad \qquad \begin{boxedminipage}{0.35\textwidth}
\[\overline{\Rightarrow X=X} \tag{\id}\]
\end{boxedminipage}\end{center}

\subsubsection{Propositional rules}\begin{center}
\begin{boxedminipage}{0.35\textwidth}
\begin{align*}
\Gamma &\Rightarrow \varphi\\
\cline{1-2}
\Gamma' &\Rightarrow \varphi \tag{\ant}
\end{align*}
if $\Gamma\subset \Gamma'$.
\end{boxedminipage}\qquad \qquad
\begin{boxedminipage}{0.35\textwidth}
\begin{align*}
\Gamma, \psi &\Rightarrow \varphi\\
\Gamma, \neg \psi &\Rightarrow \varphi\\
\cline{1-2}
\Gamma &\Rightarrow \varphi \tag{\pc}
\end{align*}
\end{boxedminipage}

\begin{boxedminipage}{0.35\textwidth}
\begin{align*}
\Gamma, \neg \varphi &\Rightarrow \psi\\
\Gamma, \neg \varphi &\Rightarrow \neg \psi\\
\cline{1-2}
\Gamma &\Rightarrow \varphi \tag{\ctr}
\end{align*}
\end{boxedminipage}\qquad \qquad
\begin{boxedminipage}{0.35\textwidth}
\begin{align*}
\Gamma, \varphi &\Rightarrow \chi\\
\Gamma, \psi &\Rightarrow \chi \\
\cline{1-2}
\Gamma, \varphi \vee \psi &\Rightarrow \chi \tag{\oorl}
\end{align*}
\end{boxedminipage}

\begin{boxedminipage}{0.35\textwidth}
\begin{align*}
\Gamma &\Rightarrow \varphi\\
\cline{1-2}
\Gamma &\Rightarrow \varphi \vee \psi \tag{\oorr a}
\end{align*}
\end{boxedminipage}\qquad \qquad
\begin{boxedminipage}{0.35\textwidth}
\begin{align*}
\Gamma &\Rightarrow \varphi\\
\cline{1-2}
\Gamma &\Rightarrow \psi \vee \varphi \tag{\oorr b}
\end{align*}
\end{boxedminipage}
\end{center}

\subsubsection{Quantifier and substitution rules}
\begin{center}
\begin{boxedminipage}{0.4\textwidth}
\begin{align*} \Gamma, \varphi \sbb{Y}{X}	 &\Rightarrow \psi\\
\cline{1-2}
\Gamma, \exists X \varphi &\Rightarrow \psi \tag{\eA}
\end{align*} if $Y$ does not occur free in the conclusion sequent.\end{boxedminipage} \qquad \qquad
\begin{boxedminipage}{0.35\textwidth}
\begin{align*} \Gamma &\Rightarrow \varphi\sbb{Y}{X}\\ 
\cline{1-2}
\Gamma, \neg E  &\Rightarrow \exists X \varphi \tag{\eS} 
\end{align*}
\end{boxedminipage}

\begin{boxedminipage}{0.43\textwidth}
\begin{align*} \Gamma, \psi\sbb{Y}{X} &\Rightarrow \varphi\sbb{Y\pr}{X\pr}\\ 
\cline{1-2}
\Gamma, \exists X \psi  &\Rightarrow \exists X\pr \varphi \tag{\eD} 
\end{align*} if $Y$ does not occur free in the conclusion sequent.
\end{boxedminipage}\qquad \qquad
\begin{boxedminipage}{0.4\textwidth}
\begin{align*} \Gamma &\Rightarrow \varphi\sbb{Y}{X}\\ 
\cline{1-2}
\Gamma, Y = Y' &\Rightarrow \varphi \sbb{Y'}{X} \tag{\subs}
\end{align*}
\end{boxedminipage}
\end{center}

\bigskip

These rules are a weakening of those of \cite{eft} which is designed to allow for soundness with respect to the empty structure. (Although we have no need to say what precisely this means.) Indeed, our version of \eS is weaker than the one given in \cite{eft}, since it just adds an antecedent. And our rule \eD simply consists of the {\it usual} \eS, followed by \eA.

Therefore we have consistency.

\subsection{Derivable rules (propositional)}
In practice we shall use more sophisticated rules than the above, which are nevertheless fully derivable in terms of these. Each derivable rule corresponds to a well-defined sequence (its derivation) of the basic rules, and this sequence could, if desired, be inserted in an algorithmically simple way into any derivation sequence containing the derivable rule, to produce a ``basic'' derivation sequence containing only the rules above. Thus, to show derivability of a particular sequent in the calculus, it suffices to establish the existence of a derivation sequence for it using the expanded ruleset (cf \cite{eft}).

We shall not give full justifications for all derivable rules here, but instead indicate how the development goes and fully justify the most nontrivial ones.\\

One first obtains several simple propositional rules:

\newcommand{\dn}{\text{DN}}

\begin{center}
\begin{boxedminipage}{0.35\textwidth}
\begin{align*}
\Gamma &\Rightarrow \varphi\\
\Gamma', \varphi &\Rightarrow \psi\\
\cline{1-2}
\Gamma \cup \Gamma' &\Rightarrow \psi \tag{\ch}
\end{align*}\end{boxedminipage}\qquad \qquad
\begin{boxedminipage}{0.35\textwidth}
\begin{align*}
\Gamma, \varphi &\Rightarrow \psi\\
\cline{1-2}
\Gamma, \neg \psi &\Rightarrow \neg \varphi \tag{\cp}
\end{align*}\end{boxedminipage}

\begin{boxedminipage}{0.35\textwidth}
\begin{align*}
\Gamma &\Rightarrow \varphi\\
\cline{1-2}
\Gamma &\Rightarrow \neg \neg \varphi \tag{\dn a}
\end{align*}\end{boxedminipage}\qquad \qquad
\begin{boxedminipage}{0.35\textwidth}
\begin{align*}
\Gamma &\Rightarrow \neg \neg \varphi\\
\cline{1-2}
\Gamma &\Rightarrow \varphi \tag{\dn b}
\end{align*}\end{boxedminipage}

\begin{boxedminipage}{0.35\textwidth}
\begin{align*}
\Gamma, \varphi&\Rightarrow \psi\\
\cline{1-2}
\Gamma, \neg \neg \varphi &\Rightarrow \psi \tag{\dn c}
\end{align*}\end{boxedminipage}\qquad \qquad
\begin{boxedminipage}{0.35\textwidth}
\begin{align*}
\Gamma, \neg \neg \varphi &\Rightarrow \psi\\
\cline{1-2}
\Gamma, \varphi &\Rightarrow \psi \tag{\dn d}
\end{align*}
\end{boxedminipage}\end{center}

\bigskip

One then defines the notations $\varphi \wedge \psi \equiv \neg (\neg \varphi \vee \neg \psi)$ and $\varphi \rightarrow \psi \equiv \neg \varphi \vee \psi$, and obtains the dual introduction rules governing $\wedge$. Importantly one also obtains\\

\begin{center}
\begin{boxedminipage}{0.35\textwidth}
\begin{align*}
\Gamma, \varphi \wedge \psi  &\Rightarrow \chi\\ 
\cline{1-2}
\Gamma, \varphi, \psi &\Rightarrow \chi \tag{\lo}
\end{align*}\end{boxedminipage}\qquad \qquad
\begin{boxedminipage}{0.35\textwidth}
\begin{align*}
\Gamma, \varphi, \psi &\Rightarrow \chi\\
\cline{1-2}
\Gamma, \varphi \wedge \psi &\Rightarrow \chi \tag{\un}
\end{align*}\end{boxedminipage}

\begin{boxedminipage}{0.35\textwidth}
\begin{align*}
\Gamma, \varphi &\Rightarrow \psi\\
\cline{1-2}
\Gamma &\Rightarrow \varphi \rightarrow \psi \tag{\meta a}
\end{align*}\end{boxedminipage}\qquad \qquad
\begin{boxedminipage}{0.35\textwidth}
\begin{align*}
\Gamma &\Rightarrow \varphi \rightarrow \psi \\
\cline{1-2}
\Gamma, \varphi &\Rightarrow \psi \tag{\meta b}
\end{align*}
\end{boxedminipage}\end{center}

\bigskip

We shall refer to rule pairs like the above two, as well as (\dn a) and (\dn b)/(\dn c) and (\dn d), as {\it reversible rules}. They let us freely pass between two sequents in a proof. The \meta pair and \lo/\un justify thinking of derivable sequents as 
implicational tautologies, with a finite (possibly empty) conjunction of antecedents on the left and a succedent on the right.\\

We also have reversible rules corresponding to commutativity and associativity for $\vee$ and $\wedge$. Beyond this chapter, we will gloss over almost all reversible rule applications.

\subsection{Predicate derivable rules\textemdash sequences}\label{seg-predderiv}
\newcommand{\meA}{\ensuremath{\exists \text{A}\pr} }
\newcommand{\meS}{\ensuremath{\exists \text{S}\pr} }
\newcommand{\msubs}{\text{SUBS}\pr }
\newcommand{\meD}{\ensuremath{\exists \text{D}\pr} }
We need versions of \eA, \eS, \eD and \subs which apply to finite sequences of variables. These will be straightforward to obtain, given the remarks of \cref{seg-subalph}.

Indeed, we claim (using the notation for simultaneous substitutions, and abbreviating the nesting $\exists X_{1} \dots \exists X_{n}$ as $\exists X_{i}$) derivability of the rules\\

\begin{center}
\begin{boxedminipage}{0.35\textwidth}
\begin{align*} \Gamma &\Rightarrow \varphi\sbb{Y_{i}}{X_{i}}\\ 
\cline{1-2}
\Gamma, \neg E  &\Rightarrow \exists X_{i} \varphi \tag{\meS} 
\end{align*} if the $X_{i}$ are distinct.\end{boxedminipage}\qquad \qquad
\begin{boxedminipage}{0.45\textwidth}
\begin{align*} \Gamma &\Rightarrow \varphi\sbb{Y_{i}}{X_{i}}\\ 
 \cline{1-2}
 \Gamma, \bigwedge_{i}Y_{i} = Y_{i}'  &\Rightarrow \varphi\sbb{Y_{i}'}{X_{i}}\tag{\msubs} 
\end{align*} if the $X_{i}$ are distinct.\end{boxedminipage}


\begin{boxedminipage}{0.35\textwidth}
\begin{align*} \Gamma, \varphi &\Rightarrow \psi\\
\cline{1-2}
\Gamma, \exists X_{i} \varphi &\Rightarrow \psi \tag{\meA}
\end{align*} if no $X_{i}$ occurs free in $\Gamma$ or $\psi$.\end{boxedminipage}\qquad \qquad
\begin{boxedminipage}{0.35\textwidth}
\begin{align*} \Gamma, \varphi &\Rightarrow \psi\\
\cline{1-2}
\Gamma, \exists X_{i} \varphi &\Rightarrow \exists X_{i}\psi \tag{\meD}
\end{align*} if no $X_{i}$ occurs free in $\Gamma$.\end{boxedminipage}
\end{center}

\bigskip

\pf For \meS, we do the argument for a sequence $X_{1}, X_{2}$ of length two; the general case is identical. By definition of simultaneous substitution, the premise is $\Gamma \Rightarrow \varphi\sbb{N_{1}}{X_{1}}\sbb{N_{2}}{X_{2}}\sbb{Y_{1}}{N_{1}}\sbb{Y_{2}}{N_{2}}$. By $\eS$, this gives $\Gamma, \neg E \Rightarrow \exists N_{2}\rak{\varphi\sbb{N_{1}}{X_{1}}\sbb{N_{2}}{X_{2}}\sbb{Y_{1}}{N_{1}}}$. The succedent is identified with $\exists N_{2}\rak{\varphi\sbb{N_{1}}{X_{1}}\sbb{N_{2}}{X_{2}}}\sbb{Y_{1}}{N_{1}}$ by setup property 1 in \cref{seg-subalph}, because neither of $Y_{1}$, $N_{1}$ are equal to $N_{2}$. Therefore by \eS again we obtain\[
\Gamma, \neg E \Rightarrow \exists N_{1}\rak{\exists N_{2}\rak{\varphi\sbb{N_{1}}{X_{1}}\sbb{N_{2}}{X_{2}}}}.
\]
Now $N_{2}$ is not in $\varphi\sbb{N_{1}}{X_{1}}$, so the inner quantification in the succedent is identified with the formula $\exists X_{2}\rak{\varphi\sbb{N_{1}}{X_{1}}}$, by setup property 2. Therefore the succedent is identified with $\exists N_{1}\rak{\exists X_{2}\rak{\varphi\sbb{N_{1}}{X_{1}}}}$. Neither of $N_{1}$, $X_{1}$ are equal to $X_{2}$, so this is identified with $\exists N_{1}\rak{\exists X_{2}\rak{\varphi}\sbb{N_{1}}{X_{1}}}$ by property 1. $N_{1}$ is not in $\exists X_{2} \varphi$, so this is identified with $\exists X_{1} \exists X_{2} \varphi$ by property 2. \qed\\

For \subs\pr, we also just do the length-two case. The premise is $\Gamma \Rightarrow \varphi\sbb{N_{1}}{X_{1}}\sbb{N_{2}}{X_{2}}\sbb{Y_{1}}{N_{1}}\sbb{Y_{2}}{N_{2}}$, and by \subs on this we have \[
\Gamma, Y_{2} = Y_{2}\pr \Rightarrow  \varphi\sbb{N_{1}}{X_{1}}\sbb{N_{2}}{X_{2}}\sbb{Y_{1}}{N_{1}}\sbb{Y_{2}\pr}{N_{2}}.
\] As $N_{1}$ and $N_{2}$ are distinct from all $Y$'s/$Y\pr$'s and each other, we can commute the last two substitutions by substitution rule 3, and this is the same sequent as\[
\Gamma, Y_{2} = Y_{2}\pr \Rightarrow  \varphi\sbb{N_{1}}{X_{1}}\sbb{N_{2}}{X_{2}}\sbb{Y_{2}\pr}{N_{2}}\sbb{Y_{1}}{N_{1}}.
\]
Applying \subs again and then switching the order back, we have\[
\Gamma, Y_{1} = Y_{1} \pr, Y_{2} = Y_{2}\pr \Rightarrow  \varphi\sbb{N_{1}}{X_{1}}\sbb{N_{2}}{X_{2}}\sbb{Y_{1}\pr}{N_{1}}\sbb{Y_{2}\pr}{N_{2}}
\]\[
\qquad \qquad  \qquad  \qquad \qquad  \overline{\Gamma, \bigwedge_{i = 1}^{2}Y_{i} = Y_{i}\pr \Rightarrow  \varphi\sbb{N_{1}}{X_{1}}\sbb{N_{2}}{X_{2}}\sbb{Y_{1}\pr}{N_{1}}\sbb{Y_{2}\pr}{N_{2}}} \qquad (\un)
\] which is the claimed conclusion. \qed\\

For \meA, we have stated a slightly weaker version than possible, because using it in practice we won't need to make nontrivial substitutions of witnesses. For the derivation we can simply observe that if $X\not\in \free(\Gamma)\cup \free(\psi)$, \eA with $Y =X$ becomes \begin{align*}
\Gamma, \varphi \Rightarrow \psi \\
\cline{1-2}
\Gamma, \exists X \varphi \Rightarrow \psi
\end{align*}
and we can iterate.\qed\\

Similarly for \meD we have not stated the strongest possible result. However to derive it we can just take $X_{i} = X_{i}\pr = Y_{i} = Y_{i}\pr$ in $\eD$ for each $i$ and iterate.\qed\\

Henceforth we shall drop the primes and assume \eS, \subs, \eA and \eD refer to these multiple-variable derivable versions. Of course the sequence of variables could have length one. It shall usually be assumed to be nonempty, however.

\section{Provable equivalence}
We will make use of the following concepts extensively throughout the paper.\\

\begin{defn} Let $S$ be a set of sentences, $\varphi$ a formula. If for some finite $\Gamma\subset S$ we have a derivable sequent\[
\Gamma \Rightarrow \varphi
\]
then we say that $S$ {\it proves} $\varphi$.
\end{defn}

\bigskip

\begin{defn} Let $S$ be a set of sentences, $\varphi$ and $\varphi'$ formulas. If for some finite $\Gamma\subset S$ we have a derivable sequent\[
\Gamma, \varphi \Rightarrow \varphi'
\]
then we say that $\varphi$ {\it $S$-provably implies $\varphi'$.}
\end{defn}
If $\varphi$ and $\varphi'$ $S$-provably imply each other then we say they are {\it $S$-provably equivalent}. It is obvious that this is an equivalence relation on formulas, for any $S$.

\bigskip

\begin{rk}\label{dec-metaisgood} By \meta (and taking unions of finite sets), we have $S$-provable implication (resp. equivalence) iff 
\[\text{$S$ proves $\varphi \rightarrow \varphi'$}\qquad \text{(resp. $S$ proves $\varphi \leftrightarrow \varphi'$).}\]
\end{rk}

\bigskip

\begin{lem}[Implication]\label{dec-implemma} Let $\varphi$ $S$-provably imply $\varphi'$. Then\begin{enumerate}[label=(\alph*)]
\item $\varphi \vee \psi$ $S$-provably implies $\varphi' \vee \psi$, and $\psi \vee \varphi$ $S$-provably implies $\psi	 \vee \varphi'$, for any $\psi$
\item $\exists X \varphi$ $S$-provably implies $\exists X \varphi'$, for any $X$
\item $\neg \varphi'$ $S$-provably implies $\neg \varphi$.
\end{enumerate}
\end{lem}

\pf (a): \[
\Gamma, \varphi \Rightarrow \varphi'
\]\[
\qquad \qquad \quad \overline{\Gamma, \varphi \Rightarrow \varphi' \vee \psi} \qquad (\oorr a) \tag{1}
\]

\[\overline{\psi \Rightarrow \psi} \qquad (\assm)
\]\[
\overline{\psi \Rightarrow \varphi' \vee \psi} \qquad (\oorr b)
\]\[
\overline{\Gamma, \psi \Rightarrow \varphi' \vee \psi} \qquad (\ant) \tag{2}
\]

\[(1), (2)
\]\[
\qquad \qquad \overline{\Gamma, \varphi \vee \psi \Rightarrow \varphi' \vee \psi} \qquad (\oorl).\qed
\]
The proof of the other claim is identical. \qed\\

(b): \[
\Gamma, \varphi \Rightarrow \varphi'\qquad
\]\[
\overline{\Gamma, \exists X \varphi \Rightarrow \exists X \varphi'} \qquad (\eD),
\]
justified as $\Gamma$ is a set of sentences, so $X$ does not occur free in it.\qed\\

(c): \[\Gamma, \varphi \Rightarrow \varphi'\qquad \qquad \qquad 
\]\[
\overline{\Gamma, \neg \varphi' \Rightarrow \neg \varphi} \qquad (\cp).
\]\\

\begin{lem}[Equivalence]\label{dec-eqlemma}
Let $\varphi$ and $\varphi'$ be $S$-provably equivalent, and let $C_{1}, \dots C_{n}$ be a sequence of unary formula building operations in the sense of \cite{end}. (So, each $C$ is equal either to $\mathcal{E}_{\vee\psi}/\mathcal{E}_{\psi \vee}$ for some $\psi$; to $\mathcal{E}_{\neg}$; or to $Q_{i}$ for some $i$.) Then $(C_{1}\circ \dots \circ C_{n})(\varphi)$ is $S$-provably equivalent to $(C_{1}\circ \dots \circ C_{n})(\varphi')$.
\end{lem}

\pf The case $n = 1$ follows directly from the Implication Lemma. The general result follows by induction, using the $n = 1$ case in the inductive step. \qed\\

\section{Universals}\label{seg-univ}
We define the notation $\forall X \equiv \neg \exists X \neg$. We shall derive the introduction rules\\

\begin{center}
\begin{boxedminipage}{0.35\textwidth}
\begin{align*} \Gamma &\Rightarrow \varphi\\ 
\cline{1-2}
\Gamma &\Rightarrow \forall X \varphi \tag{\frlS} 
\end{align*} if $X$ does not occur free in $\Gamma$.\end{boxedminipage}\qquad \qquad
\begin{boxedminipage}{0.35\textwidth}
\begin{align*} \Gamma, \varphi\sbb{Y}{X} &\Rightarrow \psi\\
\cline{1-2}
\Gamma, \neg E, \forall X \varphi &\Rightarrow \psi \tag{\frlA}
\end{align*}
\end{boxedminipage}
\end{center}

\bigskip

However first we obtain a key control on emptiness:\\

\begin{clm}\label{dec-existstononempty} Let $\varphi$ be a formula, $X$ a variable. Then we have a derivable sequent\[
\exists X \varphi \Rightarrow \neg E.
\] \end{clm}

\pf \[\overline{ \Rightarrow X = X} \qquad (\id)
\]\[
\overline{\varphi\Rightarrow X = X} \qquad (\ant)
\]\[
\overline{\exists X \varphi \Rightarrow \exists X \rak{X = X}} \qquad (\eD)
\]
where the last inference is justified as $\Gamma = \emptyset$. The (double-negation of the) succedent is identified with $\neg E$ by alpha-equivalence, clearly. \qed\\

Derivation of \frlS: Let $T$ be a variable not free in $\Gamma$ and unequal to $X$.

\[
\Gamma \Rightarrow \varphi \qquad \qquad
\]\[
\overline{\Gamma, T=T \Rightarrow \varphi} \qquad (\ant)
\]\[
\overline{\Gamma, \neg \varphi \Rightarrow \neg T = T} \qquad (\cp)
\]\[
\overline{\Gamma, \exists X\neg \varphi \Rightarrow \neg T = T} \qquad (\eA)
\]\[
\overline{\Gamma, \neg \neg T=T \Rightarrow \neg \exists X\neg \varphi} \qquad (\cp) \tag{1}
\]

\[
\overline{\Rightarrow T = T}\qquad (\id)
\]\[
\overline{\Rightarrow \neg \neg T = T} \qquad (\dn a) \tag{2}
\]

\[(2), (1)
\]\[
\qquad \qquad \qquad \qquad \overline{\Gamma \Rightarrow \neg \exists X\neg \varphi} \qquad (\ch). \qed
\]

Derivation of \frlA:

\[
\Gamma, \varphi\sbb{Y}{X} \Rightarrow \psi\qquad 
\]\[
\overline{\Gamma, \neg \psi \Rightarrow \paren{\neg \varphi}\sbb{Y}{X}} \qquad (\cp) \tag{1},
\]
 justified as we have $\neg \paren{\varphi\sbb{Y}{X}} = \paren{\neg\varphi}\sbb{Y}{X}$.\[
 (1)\]\[
\overline{\Gamma, \neg \psi, \neg E \Rightarrow \exists X \neg \varphi} \qquad (\eS)
\]\[
\overline{\Gamma, \neg E, \neg \exists X \neg \varphi \Rightarrow \neg \neg \psi} \qquad (\cp)
\]\[
\qquad \qquad \overline{\Gamma, \neg E, \neg \exists X \neg \varphi \Rightarrow \psi} \qquad (\dn b). \qed
\]\\

The generalization of \frlS to a finite sequence $(X_{i})$ of variables not free in $\Gamma$, is just by iteration, exactly as \eA was generalized above. For the multivariable version of \frlA, we must note that we again have alpha-equivalence properties\begin{enumerate}
\item $(\forall T \varphi)\sbb{Y}{X} \sim \forall T \rak{\varphi\sbb{Y}{X}} $ if $T$ is unequal to $X$ and $Y$.
\item $\forall X\varphi \sim \forall N \rak{\varphi\sbb{N}{X}}$ if $N$ is a variable not occurring free or bound in $\varphi$.
\end{enumerate} 

Indeed, for 2, the right-hand side is equal to \[
\neg \exists N \rak{\neg \paren{\varphi\sbb{N}{X}}} = \neg \exists N \rak{(\neg \varphi)\sbb{N}{X}}.
\]
By the original property 2 the subformula $\exists N \rak{(\neg \varphi)\sbb{N}{X}}$ is identified with $\exists X \rak{\neg \varphi}$. This means its negation is identified with $\neg \exists X \rak{\neg \varphi} = \forall X \varphi$, as desired.

For 1, we just use the original property 1 and the fact that $\sbb{\cdot}{\cdot}$ obviously commutes with negation.\qed

\frlA is now generalized exactly as \eS was, the argument given in \cref{seg-predderiv} having only relied on these and the properties of substitution.\qed

We shall have no need of a rule ``$\forall\text{D}$''. \frlS and \frlA shall refer to the multivariable versions, which we state for the record\\ \begin{center}
\begin{boxedminipage}{0.35\textwidth}
\begin{align*} \Gamma &\Rightarrow \varphi\\
\cline{1-2}
\Gamma &\Rightarrow \forall X_{i} \varphi \tag{\frlS final} 
\end{align*} if no $X_{i}$ occurs free in $\Gamma$.\end{boxedminipage}\qquad \qquad
\begin{boxedminipage}{0.4\textwidth}
\begin{align*} \Gamma, \varphi\sbb{Y_{i}}{X_{i}} &\Rightarrow \psi\\
\cline{1-2}
\Gamma, \neg E, \forall X_{i} \varphi &\Rightarrow \psi \tag{\frlA final}
\end{align*} if the $X_{i}$ are distinct.\end{boxedminipage}
\end{center}

\bigskip\bigskip

\begin{rk}\label{dec-emptytofrl} \cref{dec-existstononempty} (when applied to $\neg \varphi$, and technically using the Equivalence Lemma for the $\emptyset$-provably equivalent formulas $\varphi$ and $\neg \neg \varphi$) yields that for any $\varphi$ and $X$ we have\[
E \Rightarrow \forall X \varphi.
\]This will have the effect of making $\neg E$ a rather mild assumption going forward.
\end{rk}

\section{Further quantifier rules}
We derive more rules about quantifiers, using the language just defined. An $\emptyset$-provable equivalence of two different types of formulas always gives a reversible rule.\\

$(X_{i})$ will always be a sequence of length $n \geq 1$.

\subsection{Elimination}
\begin{prop}\label{dec-frelim}
If the $X_{i}$ are distinct and $(Y_{i})$ is any sequence, then $\forall X_{i}\varphi$ $\crul{\neg E}$-provably implies $\varphi\sbb{Y_{i}}{X_{i}}$.
\end{prop}

\pf \[
\overline{\varphi\sbb{Y_{i}}{X_{i}}\Rightarrow \varphi\sbb{Y_{i}}{X_{i}}}\qquad (\assm)
\]\[
\overline{\neg E, \forall X_{i} \varphi \Rightarrow \varphi\sbb{Y_{i}}{X_{i}}}\qquad (\frlA)\qed
\]

\begin{prop}\label{dec-eelim}
If the $X_{i}$ are distinct and $(Y_{i})$ is any sequence, then $\varphi\sbb{Y_{i}}{X_{i}}$ $\crul{\neg E}$-provably implies $\exists X_{i}\varphi$.
\end{prop}
\pf \[
\overline{\varphi\sbb{Y_{i}}{X_{i}}\Rightarrow \varphi\sbb{Y_{i}}{X_{i}}}\qquad (\assm)
\]\[
\overline{\neg E, \varphi\sbb{Y_{i}}{X_{i}}\Rightarrow \exists X_{i} \varphi}\qquad (\eS)\qed
\]

By combining the outputs of Propositions \ref*{dec-frelim} and \ref*{dec-eelim} with \ch, we get the following ``elimination rules'':\\

\begin{center}\begin{boxedminipage}{0.45\textwidth}
\begin{align*}
\Gamma &\Rightarrow \forall X_{i} \varphi\\
\cline{1-2}
\Gamma, \neg E&\Rightarrow \varphi\sbb{Y_{i}}{X_{i}}\tag{\frelim}
\end{align*} if the $X_{i}$ are distinct.\end{boxedminipage}\qquad
\begin{boxedminipage}{0.45\textwidth}
\begin{align*}
\Gamma, \exists X_{i} \varphi &\Rightarrow \psi\\
\cline{1-2}
\Gamma, \neg E, \varphi\sbb{Y_{i}}{X_{i}}&\Rightarrow \psi \tag{\eelim}
\end{align*} if the $X_{i}$ are distinct.\end{boxedminipage}

\end{center}

\bigskip

We can think of these as not-quite-reverses of \frlS and \eA respectively.\\

\subsection{Rearrangement}
\newcommand{\andr}{\text{\ensuremath{\bigwedge}R}}
\begin{prop}\label{dec-insex} If $X_{i}\not\in \free(\varphi)$, then the formulas $\exists X_{i}\rak{\varphi\wedge \psi}$ and $\varphi\wedge \exists X_{i}\psi$ are $\emptyset$-provably equivalent.
\end{prop}

\pf \[\overline{\varphi, \psi \Rightarrow \varphi \wedge \psi} \qquad (\assm)
\]\[\overline{\varphi, \exists X_{i}\psi \Rightarrow \exists X_{i}\rak{\varphi \wedge \psi}} \qquad (\eD),
\]
justified as $X_{i}\not \in \free(\varphi)$. Applying \un we have the direction 
\[\varphi \wedge \exists X_{i} \psi \Rightarrow \exists X_{i} \rak{\varphi\wedge \psi}.
\]

Conversely\[
\overline{\varphi, \psi \Rightarrow \varphi} \qquad (\assm)
\]\[
\qquad \overline{\varphi\wedge \psi \Rightarrow \varphi} \qquad (\un)
\]\[
\overline{\exists X_{i}\rak{\varphi\wedge \psi} \Rightarrow \varphi} \qquad (\eA) \tag{1},
\]
justified as $X_{i}\not\in \free(\varphi)$.
Similarly except using \eD instead of \eA in the last step, we obtain \[
\exists X_{i}\rak{\varphi\wedge \psi} \Rightarrow \exists X_{i} \psi \tag{2}.
\]

\[
(1), (2)
\]\[\qquad \qquad \overline{\exists X_{i}\rak{\varphi\wedge \psi} \Rightarrow \varphi \wedge \exists X_{i} \psi}\qquad (\andr)\qed
\]

We call the associated reversible rule to \cref*{dec-insex} \ex or \ins, depending on which direction we are using.\\

We note that even without $X_{i} \not\in \free(\varphi)$, we could have used \eD the first time in the converse direction above (instead of \eA) to obtain\[
\exists X_{i} \rak{\varphi \wedge \psi} \Rightarrow \exists X_{i} \varphi\tag{1\pr}
\]
and combined this with (2). Therefore we have\\

\begin{cor}\label{dec-purex}
For any $\varphi$ and $\psi$, $\exists X_{i}\rak{\varphi\wedge \psi}$ $\emptyset$-provably implies $\exists X_{i} \varphi\wedge \exists X_{i}\psi$.
\end{cor}

\bigskip

\begin{prop}\label{dec-binreord} For any $\varphi$, $i$ and $j$, $\exists x_{i}\exists x_{j} \varphi$ and $\exists x_{j} \exists x_{i} \varphi$ are $\emptyset$-provably equivalent.\end{prop}

\pf \[
\qquad \overline{\varphi \Rightarrow \varphi\sbb{x_{i}}{x_{i}}}\qquad (\assm)
\]\[
\overline{\neg E, \varphi \Rightarrow \exists x_{i}\varphi}\qquad (\eS)
\]\[
\overline{\neg E, \exists x_{j}\varphi \Rightarrow \exists x_{j} \exists x_{i} \varphi} \qquad (\eD) \tag{1},
\]
justified as $x_{j}$ is not free in (1). Now by \cref{dec-existstononempty} we have
\[\exists x_{j} \varphi \Rightarrow \neg E \tag{2}.
\]
\[(2), (1)
\]\[
\qquad \qquad \quad \overline{\exists x_{j} \varphi \Rightarrow \exists x_{j}\exists x_{i}\varphi} \qquad (\ch)
\]\[
\qquad \qquad \overline{\exists x_{i} \exists x_{j}\varphi \Rightarrow \exists x_{j} \exists x_{i} \varphi} \qquad (\eA).
\]

The converse is obtained symmetrically. \qed\\

\begin{cor}\label{dec-reord} For any $\varphi$ and permutation $\pi$, $\exists X_{i} \varphi$ and $\exists X_{\pi(i)}\varphi$ are $\emptyset$-provably equivalent.
\end{cor}

\pf Follows from the last Prop via the Equivalence Lemma, since the symmetric group is generated by adjacent transpositions. 

In more detail, $\emptyset$-p.e. is an equivalence relation. The Proposition allows us to swap the leading two quantifiers in any nesting, and by the Equivalence Lemma, we can actually do this for any two adjacent quantifiers. Thus we can pass between arbitrary elements of the symmetric group and stay within the same $\emptyset$-p.e. equivalence class.\qed\\

We will call the associated reversible rule \reord.

\section{Further misc propositional rules}\label{seg-randomprop}
We shall also reference the derived rules\\ \begin{center}
\begin{boxedminipage}{0.35\textwidth}
\begin{align*}
\Gamma, \varphi &\Rightarrow \psi\\
\cline{1-2}
\Gamma, \varphi \wedge \alpha &\Rightarrow \psi \wedge \alpha \tag{\andd}\\
\end{align*}\end{boxedminipage} \qquad \qquad \qquad
\begin{boxedminipage}{0.35\textwidth}
\begin{align*}
\Gamma &\Rightarrow\varphi \wedge \psi\\
\cline{1-2}
\Gamma &\Rightarrow  \varphi\tag{\sdrop}\\
\end{align*}\end{boxedminipage}\end{center}

\begin{center}
\begin{boxedminipage}{0.35\textwidth}
\begin{align*}
\Gamma &\Rightarrow \psi\\
\Gamma &\Rightarrow \neg \psi\\
\cline{1-2}
\Gamma &\Rightarrow \varphi \tag{\expl}\\
\end{align*}\end{boxedminipage}\end{center}

\chapter{Construction}\label{ch-construction}

\section{Setup}

In what follows, we shall sometimes write $|\varphi|$ for the number of distinct free variables occurring in $\varphi$. Thus $|s| = 0$ iff $s$ is a sentence. If $\varphi$ is a non-sentence and $(T_{i})$ is any variable sequence of length $n = |\varphi|$, we shall write $\varphi(T_{i})$ for the simultaneous substitution\[
\varphi\sbb{T_{i}}{X_{i}},
\] where $(X_{i})$ is the sequence of variables free in $\varphi$, in increasing index order.\\

If we write $\varphi(X_{i})$ without first specifying the sequence $(X_{i})$, we mean to let $X_{1}, \dots X_{n}$ {\it be} this increasing-index enumeration of $\free(\varphi)$.\\

For ease of expression, we fix the ``index-limit pairs'' $(i, n)$, $(j, m)$, and $(k, l)$. That is, unless otherwise specified, a sequence indexed by $i$ will always be taken to be of some length $n > 0$; one indexed by $j$ to be of length $m > 0$; and one indexed by $k$ to be of length $l > 0$.\\

Fix a consistent theory $T$ over $L$. That is: let $T$ be a set of sentences over $L$ such that for no formula $\psi$ do we have $T$-provability of both $\psi$ and $\neg \psi$. By virtue of the rule \expl (\cref{seg-randomprop}) this is equivalent to supposing there exists a formula which is not $T$-provable. In this chapter we begin the construction of the category $\syn(T)$.\\

$\text{Ob}(\syn(T))$ is the set of (alpha-equivalence classes of) formulas of the language of $T$.

\section{Premorphisms}
We generally denote sentences by $s, s', $ ect., and nonsentences by small Greek letters.\\

\begin{defns} Let $\mathfrak{A}$, $\mathfrak{B}$ be formulas of the language $L$ of $T$. A {\it premorphism} $\frak{A} \rightarrow \frak{B}$ is defined according to the following table.\end{defns}

\begin{center}
\begin{tabular}{|   c|   c|  c|}
\hline
Form of $\mathfrak{A} \rightarrow \mathfrak{B}$ & Premorphism type & Definition\\
\hline\hline
$\alpha(X_{i})\rightarrow \beta(Y_{j})$ & 1 & \makecell{Formula $\theta$ with free variables among $(x_{-i}) \cup (x_{j})$\\ such that $T$ proves \fone, \ftwo, \fth}\\
\hline
$\alpha(X_{i}) \rightarrow s$ &2 &Formal object $*_{\alpha s}$\\
\hline
$s \rightarrow \beta(Y_{j})$ & 3& \makecell{Formula $\theta$ with free variables among $(x_{j})$\\ such that $T$ proves \mfone, \mftwo, \mfth}\\
\hline
$s \rightarrow s'$ &4 & Formal object $*_{ss'}$\\
\hline
\end{tabular}
\end{center}

In the cases of types 1 and 3 respectively we define the sentences

\begin{itemize}
\item \fone $= \forall x_{-i}\forall x_{j}\rak{ \theta \rightarrow \alpha(x_{-i})\wedge \beta(x_{j})}$
\item \ftwo $= \forall x_{-i}\rak{ \alpha(x_{-i}) \rightarrow \exists x_{j} \theta}$
\item \fth $= \forall x_{-i}\forall x_{j}\forall x_{m + j}\rak{ \theta \wedge \theta\sbb{x_{m + j}}{x_{j}}\rightarrow \bigwedge_{j}x_{j} = x_{m + j}}$
\end{itemize}
and
\begin{itemize}
\item \mfone $=  \forall x_{j}\rak{\theta \rightarrow s \wedge \beta(x_{j})}$
\item \mftwo $= s \rightarrow \exists x_{j}\theta$
\item \mfth $= \forall x_{j}\forall x_{m + j} \rak{\theta \wedge \theta\sbb{x_{m + j}}{x_{j}}\rightarrow \bigwedge_{j}x_{j} = x_{m + j}}$.
\end{itemize}
\bigskip

In the cases of types 2 and 4 we stipulate that there is at most one such formal object. There is one \begin{itemize}
\item in case 2 if and only if $T$ proves the sentence $\forall x_{-i}\rak{ \alpha(x_{-i})\rightarrow s}$
\item in case 4 if and only if $T$ proves the sentence $s \rightarrow s'$.
\end{itemize}
\bigskip

In all cases we let $\prehom(\mathfrak{A}, \mathfrak{B})$ be the set of premorphisms.\footnote{Cf \cite{mr}, page 242.}\\

\begin{rk}\label{dec-funcforms} If a formula $\theta$ is a type (1) premorphism $\mathfrak{A} \to \mathfrak{B}$ for some $\mathfrak{A}$ and $\mathfrak{B}$, it will be called {\it functional}. We think of the sentences F1, F2, F3 as asserting that $\theta$ defines the graph of a function ``from assignments satisfying $\mathfrak{A}$ to those satisfying $\mathfrak{ B}$''.

If $\theta$ is a functional formula and $(X_{i})$ and $(Y_{j})$ are any variable sequences of the appropriate lengths, as natural in this context the notation $\theta(X_{i}, Y_{j})$ will mean the simultaneous substitution $\theta\sbx{X_{i}}{x_{-i}}{Y_{j}}{x_{j}}$ (see \cref{seg-simsub}). We can use this notation as soon as we know $\theta$ has free variables among $\crul{x_{-1}, \dots x_{-n}, x_{1}, \dots x_{m}}$; there is no need to show $T$-provability of F1, F2, F3 first.\end{rk} \bigskip

\newcommand{\mr}{\mathfrak{m}_{r}}
\newcommand{\ml}{\mathfrak{m}_{l}}
\begin{defns} Let $\ml \in \prehom(\mathfrak{A}, \mathfrak{B})$ and $\mr \in \prehom(\mathfrak{B}, \mathfrak{C})$ be premorphisms. The composition $\mr \circ \ml$ is defined according to the following table. \end{defns}

\begin{center}
\begin{tabular}{| c | c | c |}
\hline
Case & Form & Definition\\
\hline
\hline
$1 \circ 1$& $\alpha(X_{i}) \overset{\theta_{1}}{\rightarrow} \beta(Y_{j})\overset{\theta_{2}}{\rightarrow}\gamma(Z_{k})$ & \makecell{$\mr \circ \ml = \exists T_{j}\rak{\theta_{1}(x_{-i}, T_{j})\wedge \theta_{2}(T_{j}, x_{k})}$,\\ where $(T_{j})$ is some sequence of $m$ new\\ variables, in increasing index order} \\
\hline
$1 \circ 3$ & $s \overset{\theta_{1}}{\rightarrow}\beta(Y_{j})\overset{\theta_{2}}{\rightarrow}\gamma(Z_{k})$ & $\text{as above but without $x_{-i}$}$\\
\hline
$3 \circ 2$ & $\alpha(X_{i})\rightarrow s \overset{\theta_{2}}{\rightarrow} \gamma(Z_{k})$ & $\mr \circ \ml = \alpha(x_{-i})\wedge \theta_{2}(x_{k})$\\
\hline
$2 \circ 1$ & $\alpha(X_{i})\overset{\theta_{1}}{\rightarrow} \beta(Y_{j})\rightarrow s$ & $\mr \circ \ml = *_{\alpha s}$\\
\hline
$4 \circ 2$ & $\alpha(X_{i})\rightarrow s \rightarrow s'$ & $\mr \circ \ml = *_{\alpha s'}$\\
\hline
$2 \circ 3$ & $s \overset{\theta_{1}}{\rightarrow} \beta(Y_{j}) \rightarrow s'$ & $\mr \circ \ml = *_{ss'}$\\
\hline
$3\circ 4$ & $s \rightarrow s' \overset{\theta_{2}}{\rightarrow}\gamma(Z_{k})$ & $\mr \circ \ml = s \wedge \theta_{2}(x_{k})$\\
\hline
$4 \circ 4$ & $s \rightarrow s' \rightarrow s''$ & $\mr \circ \ml = *_{ss''}$\\
\hline
\end{tabular}
\end{center}
\bigskip

\begin{clm}\label{dec-prehomcomp} We have $\mr \circ \ml \in \prehom(\mathfrak{A}, \mathfrak{C})$. \end{clm}

\pf We go in order of the cases above; in each case it is clear what type of premorphism $\mr \circ \ml$ must be and hence what conditions we must check.

Case 1 ($1 \circ 1 = 1$). $\mr \circ \ml$ is a formula $\theta$ in an appropriate set of free variables. We have \[
\fone^{\theta} =  \forall x_{-i} \forall x_{k} \rak{\exists T_{j}\rak{\theta_{1}(x_{-i}, T_{j})\wedge \theta_{2}(T_{j}, x_{k})}\rightarrow \alpha(x_{-i})\wedge \gamma(x_{k})}.
\]

A general remark which will be the basis for almost all our syntactic arguments: To show $T$-provability of a sentence of the form $\forall \_\_ \phi$ such as $\fone^{\theta}$ above, it suffices to obtain\[
\Gamma, \neg E \Rightarrow \phi \tag{a}
\]
for $\Gamma$ some finite subset of $T$. This is because by assumption $T$ consists of sentences, so the LHS of (a) has no free variables, and we may apply the rule \frlS to obtain $\Gamma, \neg E\Rightarrow \forall \_\_ \phi$. We can then eliminate the $\neg E$ hypothesis by the rule \pc, as we always have $E \Rightarrow \forall \_\_ \phi$ by \cref{dec-existstononempty} in Chapter 1.

So to show $T$-provability of \fone$^{\theta}$, we will obtain\[
\Gamma, \neg E \Rightarrow \exists T_{j}\rak{\theta_{1}(x_{-i}, T_{j})\wedge \theta_{2}(T_{j}, x_{k})}\rightarrow \alpha(x_{-i})\wedge \gamma(x_{k}) \tag{*}
\]

A second remark on how our derivations will go: to derive (*) we obviously need some derivable sequents mentioning its free variables, with which to work. These are supplied by the assumed $T$-provability of \fone$^{\theta_{1}}$, \fone$^{\theta_{2}}$, ect (since $\theta_{1}$ and $\theta_{2}$ are given premorphisms). Indeed, applying the rule \frelim to the corresponding sequents, we have for example\[
\Gamma, \neg E \Rightarrow \theta_{1}(x_{-i}, T_{j}) \rightarrow \alpha(x_{-i})\wedge \beta(T_{j}) \tag{1}
\]
\[
\Gamma, \neg E \Rightarrow \theta_{2}(T_{j}, x_{k}) \rightarrow \beta(T_{j})\wedge \gamma(x_{k}) \tag{2}.
\]

Now, by \meta (1) and (2) give\[
\Gamma, \neg E, \theta_{1}(x_{-i}, T_{j}) \Rightarrow \alpha(x_{-i})\wedge \beta(T_{j}) \tag{1$\pr$}
\]
\[
\Gamma, \neg E, \theta_{2}(T_{j}, x_{k}) \Rightarrow \beta(T_{j})\wedge \gamma(x_{k}) \tag{2$\pr$}.
\]
We then proceed \begin{paracol}{2}
\qquad \qquad \qquad (1$\pr$)
\[\overline{\Gamma, \neg E, \theta_{1}(x_{-i}, T_{j}) \Rightarrow \alpha(x_{-i})} (\sdrop) \tag{3}
\]
\switchcolumn \qquad \qquad \qquad (2$\pr$)
\[\overline{\Gamma, \neg E, \theta_{2}(T_{j}, x_{k}) \Rightarrow \gamma(x_k)} (\sdrop) \tag{4}
\]
\end{paracol}\[(3), (4)
\]\[\overline{\Gamma, \neg E, \theta_{1}(x_{-i}, T_{j}) \wedge \theta_{2}(T_{j}, x_{k}) \Rightarrow \alpha(x_{-i}) \wedge \gamma(x_k)} \qquad (\andd) \tag{5}
\]\[
\overline{\Gamma, \neg E, \exists T_{j}\rak{\theta_{1}(x_{-i}, T_{j}) \wedge \theta_{2}(T_{j}, x_{k})} \Rightarrow \alpha(x_{-i}) \wedge \gamma(x_k)} \qquad (\eA) \tag{6}
\]\[
\overline{\Gamma, \neg E \Rightarrow \exists T_{j}\rak{\theta_{1}(x_{-i}, T_{j}) \wedge \theta_{2}(T_{j}, x_{k})} \rightarrow \alpha(x_{-i}) \wedge \gamma(x_k)}\qquad (\meta) \tag{*}
\]

and the $T$-provability of the first functionality sentence \fone$^{\mathfrak{m}_{r}\circ \mathfrak{m}_{l}}$ in Case 1 is shown.

Given the somewhat roundabout nature of this argument, a reasonable question is whether it wouldn't be better to omit universal quantifiers from the premorphism conditions altogether. The reason we don't do this is because of the freedom to choose variables granted by $\frelim$: we could have invoked (1) and (2) with any free variables we like, not just the $T$'s and $x_{-i}$/$x_{k}$'s. This will be necessary in more complex setups.

Rather than spelling out the ``setup reasoning'' every time, many of our demonstrations of $T$-provability of a universal sentence will begin simply by postulating a goal sequent $*$ as above, and proceed to derive $*$ using any \frelim's we need of known $T$-provable universal sentences. As the situations get more complex we must always be sure to indicate why $*$ leads quickly to $\Gamma \Rightarrow (\text{desired sentence})$.

We now proceed to show $T$-provability of F2 and F3 by this method, finishing Case 1 of \cref{dec-prehomcomp}.\\


F2: \[\Gamma, \neg E \Rightarrow \alpha(x_{-i})\rightarrow \exists x_{k}\exists T_{j}\rak{\theta_{1}(x_{-i}, T_{j})\wedge \theta_{2}(T_{j}, x_{k})}\tag{*}
\]

We have
\[
\Gamma, \neg E \Rightarrow \beta(T_{j})\rightarrow \exists x_{k}\rak{\theta_{2}(T_{j}, x_{k})} \qquad (\text{\frelim on $\ftwo^{\theta_{2}}$}) \tag{7}
\]
\[
\overline{\Gamma, \neg E, \beta(T_{j})\Rightarrow \exists x_{k}\rak{\theta_{2}(T_{j}, x_{k})}} \qquad (\meta) \tag{$7\pr$}
\]
\[
\overline{\Gamma, \neg E, \beta(T_{j}) \wedge \theta_{1}(x_{-i}, T_{j}) \Rightarrow \exists x_{k}\rak{\theta_{2}(T_{j}, x_{k})} \wedge \theta_{1}(x_{-i}, T_{j})} \qquad (\andd)\tag{8}
\]
\[
\overline{\Gamma, \neg E, \beta(T_{j}) \wedge \theta_{1}(x_{-i}, T_{j}) \Rightarrow \exists x_{k}\rak{\theta_{2}(T_{j}, x_{k}) \wedge \theta_{1}(x_{-i}, T_{j})}} \qquad (\ins) \tag{9}
\]
\[
\overline{\Gamma, \neg E, \beta(T_{j}), \theta_{1}(x_{-i}, T_{j}) \Rightarrow \exists x_{k}\rak{\theta_{2}(T_{j}, x_{k}) \wedge \theta_{1}(x_{-i}, T_{j})}} \qquad (\un) \tag{10}
\]

Now we also have\[
(1\pr \text{ from above}) \qquad
\]\[
\overline{\Gamma, \neg E, \theta_{1}(x_{-i}, T_{j}) \Rightarrow \beta(T_{j})} \qquad (\sdrop) \tag{11}
\]

Therefore \[
(11), (10)
\]\[
\overline{\Gamma, \neg E, \theta_{1}(x_{-i}, T_{j}) \Rightarrow \exists x_{k}\rak{\theta_{2}(T_{j}, x_{k}) \wedge \theta_{1}(x_{-i}, T_{j})}} \qquad (\ch) \tag{12}
\]
\[
\overline{\Gamma, \neg E, \exists T_{j}\rak{\theta_{1}(x_{-i}, T_{j})} \Rightarrow \exists T_{j}\exists x_{k}\rak{\theta_{2}(T_{j}, x_{k}) \wedge \theta_{1}(x_{-i}, T_{j})}} \qquad (\eD) \tag{13}
\]

\[
\qquad \Gamma, \neg E \Rightarrow \alpha(x_{-i}) \rightarrow \exists T_{j}\rak{\theta_{1}(x_{-i}, T_{j})} \qquad (\frelim \ftwo^{\theta_{1}}) \tag{14}
\]\[
\overline{\Gamma, \neg E, \alpha(x_{-i}) \Rightarrow \exists T_{j}\rak{\theta_{1}(x_{-i}, T_{j})}} \qquad (\meta) \tag{14$\pr$}
\]

\[
(13), (14\pr)
\]\[
\overline{\Gamma, \neg E, \alpha(x_{-i}) \Rightarrow \exists T_{j}\exists x_{k}\rak{\theta_{2}(T_{j}, x_{k}) \wedge \theta_{1}(x_{-i}, T_{j})}} \qquad (\ch) \tag{15}
\]\[
\overline{\Gamma, \neg E, \alpha(x_{-i}) \Rightarrow \exists x_{k}\exists T_{j}\rak{\theta_{2}(T_{j}, x_{k}) \wedge \theta_{1}(x_{-i}, T_{j})}} \qquad (\reord) \tag{16}
\]
\[
\overline{\Gamma, \neg E \Rightarrow \alpha(x_{-i}) \rightarrow \exists x_{k}\exists T_{j}\rak{\theta_{2}(T_{j}, x_{k}) \wedge \theta_{1}(x_{-i}, T_{j})}} \qquad (\meta)
\]

which is (*), up to commutativity of $\wedge$, which we will generally not bother about. \qtw

In the rest of our derivations, for compactness, we may perform multiple simple rule applications across one horizontal bar and not name all of them.\\

F3: \[
\Gamma, \neg E, \exists T_{j}\rak{\theta_{1}(x_{-i}, T_{j})\wedge \theta_{2}(T_{j}, x_{k})} \wedge \exists T_{j}\rak{\theta_{1}(x_{-i}, T_{j})\wedge \theta_{2}(T_{j}, x_{l + k})} \Rightarrow \bigwedge_{k} x_{k} = x_{l + k} \tag{*}
\]

Let $(L_{j})$ be a new increasing sequence. We have \[
\qquad \Gamma, \neg E \Rightarrow \theta_{1}(x_{-i}, T_{j})\wedge \theta_{1}(x_{-i}, L_{j}) \rightarrow \bigwedge_{j} T_{j} = L_{j} \qquad (\frelim \fth^{\theta_{1}}) \tag{1}
\]
\[
\overline{\Gamma, \neg E, \theta_{1}(x_{-i}, T_{j})\wedge \theta_{1}(x_{-i}, L_{j}) \Rightarrow \bigwedge_{j} T_{j} = L_{j}} \qquad \tag{1$\pr$}
\]

\[
\theta_{2}(L_{j}, x_{l + k}) \Rightarrow \theta_{2}(L_{j}, x_{l + k}) \qquad (\assm) \tag{2}
\]
\[
\overline{\theta_{2}(L_{j}, x_{l + k}), \bigwedge_{j} T_{j} = L_{j} \Rightarrow \theta_{2}(T_{j}, x_{l + k})} \qquad (\subs) \tag{3}
\]

\[
(1\pr), (3)
\]
\[
\overline{\Gamma, \neg E, \theta_{1}(x_{-i}, T_{j})\wedge \theta_{1}(x_{-i}, L_{j}), \theta_{2}(L_{j}, x_{l + k}) \Rightarrow \theta_{2}(T_{j}, x_{l + k})} \qquad (\ch) \tag{4}
\]
\[
\overline{\Gamma, \neg E, \theta_{1}(x_{-i}, T_{j})\wedge \theta_{1}(x_{-i}, L_{j}), \theta_{2}(L_{j}, x_{l + k}), \theta_{2}(T_{j}, x_{k})\Rightarrow \theta_{2}(T_{j}, x_{l + k})} \qquad (\ant) \tag{5}
\]
\[
\overline{\Gamma, \neg E, \theta_{1}(x_{-i}, T_{j}), \theta_{2}(T_{j}, x_{k}), \theta_{1}(x_{-i}, L_{j}), \theta_{2}(L_{j}, x_{l + k})\Rightarrow \theta_{2}(T_{j}, x_{l + k})} \qquad \tag{5\pr}
\]

We may certainly add $\theta_{2}(T_{j}, x_{k})$ to the RHS of (5\pr) as it is already on the LHS (formally we use \assm to obtain a new sequent $\theta_{2}(T_{j}, x_{k})\Rightarrow \theta_{2}(T_{j}, x_{k})$, add the other hypotheses which appear on the LHS of 5\pr to it by \ant, and apply \andd and \lo) to obtain\[
\Gamma, \neg E, \theta_{1}(x_{-i}, T_{j}), \theta_{2}(T_{j}, x_{k}), \theta_{1}(x_{-i}, L_{j}), \theta_{2}(L_{j}, x_{l + k})\Rightarrow \theta_{2}(T_{j}, x_{k})\wedge \theta_{2}(T_{j}, x_{l + k}) \qquad \tag{5\pr\pr}
\]

We have by \frelim on $\fth^{\theta_{2}}$\[
\Gamma, \neg E, \theta_{2}(T_{j}, x_{k})\wedge \theta_{2}(T_{j}, x_{l + k}) \Rightarrow \bigwedge_{k} x_{k} = x_{l + k} \tag{6}
\] and \ch ing with (5\pr\pr) we get\[
\Gamma, \neg E, \theta_{1}(x_{-i}, T_{j}), \theta_{2}(T_{j}, x_{k}), \theta_{1}(x_{-i}, L_{j}), \theta_{2}(L_{j}, x_{l + k})\Rightarrow \bigwedge_{k} x_{k} = x_{l + k} \tag{7}
\]
\[
\overline{\Gamma, \neg E, \theta_{1}(x_{-i}, T_{j})\wedge \theta_{2}(T_{j}, x_{k}), \theta_{1}(x_{-i}, L_{j}) \wedge \theta_{2}(L_{j}, x_{l + k})\Rightarrow \bigwedge_{k} x_{k} = x_{l + k}}\qquad (\un)
\]
\[
\overline{\Gamma, \neg E, \theta_{1}(x_{-i}, T_{j})\wedge \theta_{2}(T_{j}, x_{k}), \exists L_{j}\rak{\theta_{1}(x_{-i}, L_{j}) \wedge \theta_{2}(L_{j}, x_{l + k})}\Rightarrow \bigwedge_{k} x_{k} = x_{l + k}}\qquad (\eA)
\]
\[
\overline{\Gamma, \neg E, \exists T_{j}\rak{\theta_{1}(x_{-i}, T_{j})\wedge \theta_{2}(T_{j}, x_{k})}, \exists L_{j}\rak{\theta_{1}(x_{-i}, L_{j}) \wedge \theta_{2}(L_{j}, x_{l + k})}\Rightarrow \bigwedge_{k} x_{k} = x_{l + k}}\qquad (\eA)
\]
\[
\overline{\Gamma, \neg E, \exists T_{j}\rak{\theta_{1}(x_{-i}, T_{j})\wedge \theta_{2}(T_{j}, x_{k})}\wedge \exists L_{j}\rak{\theta_{1}(x_{-i}, L_{j}) \wedge \theta_{2}(L_{j}, x_{l + k})}\Rightarrow \bigwedge_{k} x_{k} = x_{l + k}}
\]
which is clearly (*) up to alpha-equivalence. \qth \qedone\\

Case 2 ($1 \circ 3 = 3$). The composition $\mr \circ \ml$ is given by the same formula as in Case 1, except that now $\theta_{1}$ has no negative free variables.\\

mF1:\[
\Gamma, \neg E, \exists T_{j} \rak{\theta_{1}(T_{j})\wedge \theta_{2}(T_{j}, x_{k})}\Rightarrow s \wedge \gamma(x_{k}) \tag{*}
\]
The derivation runs identically to that of F1 in Case 1; absence of the $x_{-i}$ does not affect anything.\qed\\


mF2 (note mF2 is not a universal sentence):\[
\Gamma, s \Rightarrow \exists x_{k} \rak{\exists T_{j}\rak{ \theta_{1}(T_{j})\wedge \theta_{2}(T_{j}, x_{k})} } \tag{*}
\]

By \frelim $\mfone^{\theta_{1}}$ and $\ftwo^{\theta_{2}}$ respectively, have \begin{paracol}{2}\[
\Gamma, \neg E, \theta_{1}(T_{j}) \Rightarrow \beta(T_{j})\tag{1}
\] \switchcolumn  \[
\Gamma, \neg E, \beta(T_{j})\Rightarrow \exists x_{k}\rak{\theta_{2}(T_{j}, x_{k})}\tag{2}
\] 
\end{paracol}
Chaining (1) and (2) we have\[
\Gamma,\neg E, \theta_{1}(T_{j}) \Rightarrow \exists x_{k}\rak{\theta_{2}(T_{j}, x_{k})}\tag{3}
\]
\[
\overline{\Gamma,\neg E, \theta_{1}(T_{j}) \Rightarrow \theta_{1}(T_{j}) \wedge \exists x_{k}\rak{\theta_{2}(T_{j}, x_{k})}} \tag{3\pr}
\]
\[
\qquad \qquad \overline{\Gamma,\neg E, \theta_{1}(T_{j}) \Rightarrow \exists x_{k}\rak{ \theta_{1}(T_{j}) \wedge \theta_{2}(T_{j}, x_{k})}}\qquad (\ins) \tag{3\pr\pr}
\]
\[
\qquad \qquad \overline{\Gamma, \neg E, \exists T_{j}\rak{\theta_{1}(T_{j})} \Rightarrow \exists T_{j}\exists x_{k}\rak{ \theta_{1}(T_{j}) \wedge \theta_{2}(T_{j}, x_{k})}}\qquad (\eD)\tag{4}
\]\[
\overline{\Gamma, \neg E, \exists T_{j}\rak{\theta_{1}(T_{j})} \Rightarrow \exists x_{k}\exists T_{j}\rak{ \theta_{1}(T_{j}) \wedge \theta_{2}(T_{j}, x_{k})}} \tag{4\pr}
\]

Now we may drop the $\neg E$ in (4\pr) because the other hypothesis $\exists T_{j}\rak{\theta_{1}(T_{j})}$, being an existential sentence, already implies it. Let this be (4\pr\pr).
By $T$-provability of $\mftwo^{\theta_{1}}$ we have \[
\Gamma, s \Rightarrow \exists T_{j}\theta_{1}(T_{j})
\]
and chaining this with (4\pr\pr) we obtain (*).\qed\\

mF3: Identical to F3 in Case 1; absence of $x_{-i} $ makes no difference. \qed \qed\\

Case 3 ($3 \circ 2 = 1$). We have $\mr \circ \ml = \alpha(x_{-i})\wedge \theta_{2}(x_{j})$.\\

F1:\[
\Gamma, \neg E, \alpha(x_{-i})\wedge \theta_{2}(x_{j}) \Rightarrow \alpha(x_{-i})\wedge \gamma(x_{k}) \tag{*}
\]
By \frelim on $\mfone^{\theta_{2}}$ we have\[
\Gamma, \neg E, \theta_{2}(x_{k})\Rightarrow \gamma(x_{k})\tag{2}
\]
and (*) follows immediately.\qed\\ 

F2: Since $\ml$ is a type 2 premorphism $\alpha \to s$, we have that $T$ proves $\forall x_{-i} \rak{\alpha(x_{-i})\rightarrow s}$. Therefore by \frelim and trivial manipulations we have\[
\Gamma, \neg E, \alpha(x_{-i}) \Rightarrow s\tag{3}
\]

By $T$-provability of $\mftwo^{\theta_{2}}$ we have\[
\Gamma, s\Rightarrow \exists x_{k}\theta_{2}(x_{k}) \tag{4}
\]
and chaining (3) and (4) gives \[
\Gamma, \neg E, \alpha(x_{-i}) \Rightarrow \exists x_{k}\rak{\theta_{2}(x_{k})}\tag{5}
\]\[
\overline{\Gamma, \neg E, \alpha(x_{-i}) \Rightarrow \alpha(x_{-i})\wedge \exists x_{k}\rak{\theta_{2}(x_{k})}}\tag{5\pr}
\]
The $x_{k}$ and $x_{-i}$ are distinct variables, so we may apply \ins to (5\pr) to obtain (*).\qed\\

F3:\[
\Gamma, \neg E, \alpha(x_{-i})\wedge \theta_{2}(x_{k}), \alpha(x_{-i})\wedge \theta_{2}(x_{k + l}) \Rightarrow \bigwedge_{k} x_{k} = x_{k + l}\tag{*}
\]

By \frelim on $\mfth^{\theta_{2}}$ we have\[
\Gamma, \neg E, \theta_{2}(x_{k}), \theta_{2}(x_{k + l}) \Rightarrow \bigwedge_{k} x_{k} = x_{k + l} \tag{6}
\]
and (*) follows immediately.\qed \qed\\

Case 4 ($2 \circ 1 = 2$). We just need to show that we have the formal object $*_{\alpha s}$, which by definition is just the claim that $T$ proves $\forall x_{-i} \rak{\alpha(x_{-i})\rightarrow s}$. So we need \[
\Gamma, \neg E, \alpha(x_{-i}) \Rightarrow s \tag{*}
\]
Since $*_{\beta s}$ exists, we have\[
\Gamma, \neg E, \beta(T_{j}) \Rightarrow s\qquad  \tag{1}
\]\[
\overline{\Gamma, \neg E, \exists T_{j}\beta(T_{j}) \Rightarrow s} \qquad (\eA)\tag{1\pr}
\]

We have
\[
\qquad \qquad \qquad \Gamma, \neg E, \theta_{1}(x_{-i}, T_{j}) \Rightarrow \beta(T_{j})\qquad (\frelim \text{on $\fone^{\theta_{1}}, \sdrop$}) \tag{2}
\]\[
\overline{\Gamma, \neg E, \exists T_{j} \theta_{1}(x_{-i}, T_{j}) \Rightarrow \exists T_{j}\beta(T_{j})}\qquad \qquad (\eD) \tag{3}
\]

By \frelim on $\ftwo^{\theta_{1}}$ have
\[
\Gamma, \neg E, \alpha(x_{-i}) \Rightarrow \exists T_{j} \theta_{1}(x_{-i}, T_{j}) \tag{4}
\]

\[
(4), (3)
\]\[
\overline{\Gamma, \neg E, \alpha(x_{-i}) \Rightarrow \exists T_{j} \beta(T_{j})} \qquad (\ch) \tag{5}
\]

and chaining (5) and (1\pr) we have (*).\qed\\

Case 5 ($4 \circ 2 = 2$). The setup is\[
\alpha(X_{i})\rightarrow s \rightarrow s'
\]

To show existence of $*_{\alpha s'}$ we need\[
\Gamma, \neg E, \alpha(x_{-i}) \Rightarrow s' \tag{*}
\]

This is immediate from\[
\Gamma, \neg E, \alpha(x_{-i}) \Rightarrow s 
\]

and\[
\Gamma, s \Rightarrow s'.
\]\qed\\

Case 6 ($2 \circ 3 = 4$). The setup is\[
s \overset{\theta_{1}}{\rightarrow} \beta(Y_{j}) \rightarrow s'
\]

To show existence of $*_{s s'}$ we need\[
\Gamma, s \Rightarrow s' \tag{*}
\]

Since $*_{\beta s'}$ exists we have \[
\Gamma, \neg E, \beta(x_{j}) \Rightarrow s' \qquad \qquad \tag{1}
\]\[
\overline{\Gamma, \neg E, \exists x_{j}\beta(x_{j}) \Rightarrow s'} \qquad (\eA) \tag{1\pr}
\]

We have by $T$-provability of $\mftwo^{\theta_{1}}$\[
\Gamma, s \Rightarrow \exists x_{j} \theta_{1}(x_{j}) \tag{2}
\]

and by \frelim $\mfone^{\theta_{1}}$
\[
\Gamma, \neg E, \theta_{1}(x_{j}) \Rightarrow \beta(x_{j}) \tag{3}
\]\[
\overline{\Gamma, \neg E,\exists x_{j} \theta_{1}(x_{j}) \Rightarrow \exists x_{j}\beta(x_{j})}\qquad (\eD) \tag{3\pr}
\]

\[
(2),(3\pr)\qquad \qquad
\]
\[
\qquad \qquad \overline{\Gamma, \neg E, s \Rightarrow\exists x_{j}\beta(x_{j})} \qquad (\ch)\tag{4}
\]

\[
(4),(1')\qquad \qquad
\]
\[
\qquad \qquad \overline{\Gamma, \neg E, s \Rightarrow s'} \qquad (\ch)\tag{5}
\]

Now in fact by (2) we must have $\Gamma, s \Rightarrow \neg E$, by \ch and \cref{dec-existstononempty} of Chapter 1. Therefore we can eliminate the $\neg E$ hypothesis from (5), to obtain *. \qed\\

Case 7 ($3 \circ 4 = 3$). The setup is\[
s \rightarrow s' \overset{\theta_{2}}{\rightarrow}\gamma(Z_{k})
\]
and we have $\mr\circ\ml = s \wedge \theta_{2}(x_{k})$. Thus we check\\

$\mfone:$\[
\Gamma, \neg E,  s \wedge \theta_{2}(x_{k}) \Rightarrow s\wedge\gamma(x_{k}) \tag{*}
\]
and the derivation is trivial from $\mfone^{\theta_{2}}$.\qed\\

$\mftwo:$\[
\Gamma, s \Rightarrow \exists x_{k} \rak{s \wedge \theta_{2}(x_{k})} \tag{*}
\] We have \[
\Gamma, s \Rightarrow s' \tag{1}\]
by existence of $*_{ss'}$. By $\mftwo^{\theta_{2}}$, also\[
\Gamma, s' \Rightarrow \exists x_{k}\rak{\theta_{2}(x_{k})} \tag{2}
\]
and chaining (1) and (2) gives \[
\Gamma, s \Rightarrow\exists x_{k}\rak{\theta_{2}(x_{k})}\tag{3}
\]
from which we obtain (*) easily as $s$ is a sentence.\qed\\

mF3:\[
\Gamma, \neg E, s\wedge \theta_{2}(x_{k}), s\wedge \theta_{2}(x_{k + l}) \Rightarrow \bigwedge_{k} x_{k} = x_{k + l} \tag{*}
\]
We have by $\mfth^{\theta_{2}}$\[
\Gamma, \neg E, \theta_{2}(x_{k}), \theta_{2}(x_{k + l}) \Rightarrow \bigwedge_{k} x_{k} = x_{k + l}
\]
so this is immediate. \qed\\

Case 8 ($4 \circ 4  = 4$): We have $\Gamma, s \Rightarrow s'$ and $\Gamma, s' \Rightarrow s''$, so $\Gamma, s\Rightarrow s''$ by \ch. \qed\\

This concludes the proof of \cref{dec-prehomcomp}. \qed\\

\section{Morphisms}
We now define the morphisms of our category.\\

\begin{defn}[nonsentence target] Let $\mathfrak{A}$ be a formula, and $\psi(Y_{j})$ a nonsentence. We set $\hom(\mathfrak{A}, \psi) = \prehom(\mathfrak{A}, \psi)/\sim$, where $\sim$ is the equivalence relation on $\prehom(\mathfrak{A}, \psi)$ defined as follows: $\theta \sim \theta '$ if and only if $\theta$, $\theta'$ are $T\cup \{\neg E\}$-provably equivalent, i.e. there is a derivable sequent\[
\Gamma, \neg E \Rightarrow \theta \leftrightarrow \theta'
\]
for some finite $\Gamma \subset T$. \end{defn}\bigskip

\begin{defn}[sentence target] Let $\mathfrak{A}$ be a formula, and $s$ a sentence. We set $\hom(\mathfrak{A}, s) = \prehom(\mathfrak{A}, s)$. \end{defn}\bigskip

\begin{rk}\label{dec-morsarethis} By \frelim, this is equivalent to stipulating the following: In the $\mathfrak{A}\to s$ (type 2 or 4) situation, a morphism is just a premorphism. When the target is a nonsentence (type 1 or 3 situation), a morphism is represented by a premorphism $\theta$, which is a formula as specified above; two premorphisms $\theta$, $\theta'$ represent the same morphism just when $T$ proves the sentence $\E^{\theta \leftrightarrow \theta'}:= \forall \text{(free variables of $\mathfrak{A}$ and $\psi$)}\rak{\theta \leftrightarrow \theta'}$.\end{rk}\bigskip

\begin{prop}\label{dec-fullreps} Let $\varphi(X_{i})\rightarrow \psi(Y_{j})$ be a type 1 or 3 morphism (so $(X_{i})$ possibly the empty sequence). Then there exists a representative $\Theta(x_{-i}, x_{j})$ (resp. $\Theta(x_{j})$) in which all the $x_{-i}$ and $x_{j}$ (resp. $x_{j}$) occur free.\end{prop}

\pf Take any representative $\theta$, and consider the formula $\theta \wedge \bigwedge_{i} x_{-i} = x_{-i}\wedge \bigwedge_{j} x_{j} = x_{j}$. We have by instances of \id and \andd the empty-LHS sequent\[
\Rightarrow \bigwedge_{i} x_{-i} = x_{-i}\wedge \bigwedge_{j} x_{j} = x_{j}\tag{1}
\]
and by \assm\[
\theta \Rightarrow \theta \tag{2}
\]

\[
(1)(2)
\]
\[\overline{\theta \Rightarrow \theta \wedge \bigwedge_{i} x_{-i} = x_{-i}\wedge \bigwedge_{j} x_{j} = x_{j}} \qquad \qquad (\andd)\tag{3}
\]

Conversely we certainly have \[
\theta \wedge \bigwedge_{i} x_{-i} = x_{-i} \wedge \bigwedge_{j} x_{j} = x_{j} \Rightarrow \theta
\]
by \assm and \sdrop. Therefore setting $\Theta$ equal to this formula, we have that $\theta$ and $\Theta$ are $\emptyset$-provably equivalent, hence certainly $T$-provably equivalent. It follows by the Equivalence Lemma that $T$ proves all the functionality sentences for $\Theta$, hence $\Theta$ is a premorphism $\varphi(X_{i})\rightarrow \psi(Y_{j})$. It represents the same morphism as $\theta$, since the formulas $\theta,\Theta$ are certainly also $T \cup \{\neg E\}$-provably equivalent.\qed\\

From now on, unqualified references to ``provable equivalence/implication" between formulas or ``provability'' of a formula shall refer to $S = T \cup \{\neg E\}$.\\ 

\begin{prop}\label{dec-descendants} The composition of premorphisms we have defined descends to morphisms.\end{prop}

\pf This is immediate from the Equivalence Lemma, since in each case the composition is either formal, or a well-defined construction sequence in terms of ``constant formulas'' $\alpha$, $s$, etc, and the input premorphisms $\theta_{1}$ and $\theta_{2}$. For example, in Case 3 this sequence is $\mathcal{E}_{\neg}(\mathcal{E}_{\vee}(\mathcal{E}_{\neg}(\alpha\sbb{x_{-i}}{X_{i}}), \mathcal{E}_{\neg}(\theta_{2})))$. The Equivalence Lemma says the provable equivalence class of this construction is invariant under varying either of the inputs within their provable equivalence classes. Hence the morphism represented by the composition is invariant under varying the premorphism representatives of either input morphism.\qed\\

In light of \cref{dec-descendants}, from now on in arguments we will mostly identify (non-formal) morphisms with their representatives, i.e. with formulas $\theta$. By \cref{dec-fullreps} we can always assume $\theta$ actually contains all possible free variables. \\

\subsection{Associativity}
\begin{prop}\label{dec-assoc} The composition of morphisms is associative. \end{prop}

\pf As will become our general pattern, we first do the case where all morphisms are of type 1 (to be called, throughout the paper, ``Case 0''); and then we indicate modifications to the argument in other cases. So let \begin{center}\begin{tikzcd}
\alpha(X_{i}) \arrow[r, "\theta_{1}"] & \beta(Y_{j}) \arrow[r, "\theta_{2}"] & \gamma(Z_{k}) \arrow[r, "\theta_{3}"]  & \delta(P_{w})
\end{tikzcd}\end{center}

$\theta_{3}\circ ( \theta_{2} \circ \theta_{1})$ is defined as \[
\exists M_{k}\rak{\exists T_{j}\rak{\theta_{1}(x_{-i}, T_{j})\wedge \theta_{2}(T_{j}, M_{k})}\wedge \theta_{3}(M_{k}, x_{w})} \tag{1}\]

$x_{w}$ and $M_{k}$ are taken distinct from $T_{j}$, so by \ins and \ex the subformula $\exists T_{j}\rak{\theta_{1}(x_{-i}, T_{j})\wedge \theta_{2}(T_{j}, M_{k})}\\ \wedge \theta_{3}(M_{k}, x_{w})$ is provably equivalent to $\exists T_{j}\rak{\theta_{1}(x_{-i}, T_{j})\wedge \theta_{2}(T_{j}, M_{k})\wedge \theta_{3}(M_{k}, x_{w})}$. Therefore by the Equivalence Lemma (1) is provably equivalent to \[
\exists M_{k}\exists T_{j}\rak{\theta_{1}(x_{-i}, T_{j})\wedge \theta_{2}(T_{j}, M_{k})\wedge \theta_{3}(M_{k}, x_{w})}\tag{1$\pr$}
\]

On the other hand we have\[
(\theta_{3}\circ  \theta_{2}) \circ \theta_{1} = \exists T_{j} \rak{\theta_{1}(x_{-i}, T_{j})\wedge \exists M_{k }\rak{\theta_{2}(T_{j}, M_{k})\wedge \theta_{3}(M_{k}, x_{w})}}\tag{2}
\]
and by similar reasoning the RHS of (2) is provably equivalent to \[
(\theta_{3}\circ  \theta_{2}) \circ \theta_{1} = \exists T_{j}\exists M_{k}\rak{\theta_{1}(x_{-i}, T_{j})\wedge \theta_{2}(T_{j}, M_{k})\wedge \theta_{3}(M_{k}, x_{w})}\tag{2$\pr$}\]
By \reord this is provably equivalent to (1$\pr$). Therefore by transitivity of provable equivalence, (1) is p.e. to the RHS of (2). Therefore they represent the same morphism $\alpha \to \delta$. \qedo\\ 

There are {\it a priori} 15 other cases, but we observe that if the last formula $\delta$ is a sentence there is nothing to check, since both compositions must be the formal object $*_{\alpha \delta}$. Therefore we have in fact 7 other cases to check, each where at least one of $\alpha$, $\beta$, $\gamma$ in 

\begin{center}
\begin{tikzcd}
\alpha \arrow[r]& \beta \arrow[r] & \gamma \arrow[r, "\theta_{3}"] & \delta(P_{w})
\end{tikzcd}
\end{center}

is a sentence. All will be fairly trivial. Since sentences later in the chain simplify things, we organize them as follows.\\

1. Suppose $\gamma = s''$ is a sentence.\\

1.1. Suppose $\alpha = s$ is a sentence.\\

1.1.1. Suppose $\beta = s'$ is a sentence. Then our diagram is 

\begin{center}
\begin{tikzcd}
s \arrow[r, "*_{1}"]& s' \arrow[r, "*_{2}"] & s'' \arrow[r, "\theta_{3}(x_{w})"] & \delta(P_{w})
\end{tikzcd}
\end{center}


We have $(\theta_{3}  \circ  *_{2})  \circ *_{1} = (s' \wedge \theta_{3})\circ *_{1} = s \wedge s' \wedge \theta_{3}(x_{w})$. On the other hand

\[\theta_{3}  \circ  (*_{2}  \circ *_{1}) = \theta_{3}\circ *_{s s''} = s\wedge \theta_{3}(x_{w})\]

These are trivially provably equivalent since $s$ $T$-provably implies $s'$, by existence of the type 4 morphism $s\to s'$. \qed\\

1.1.2. Suppose $\beta = \beta(Y_{j})$. Then our diagram is 

\begin{center}
\begin{tikzcd}
s \arrow[r, "\theta_{1}(x_{j})"]& \beta(Y_{j}) \arrow[r, "*_{2}"] & s'' \arrow[r, "\theta_{3}(x_{w})"] & \delta(P_{w})
\end{tikzcd}
\end{center}


We have $(\theta_{3}  \circ  *_{2})  \circ \theta_{1} = (\beta(x_{-j})\wedge \theta_{3}(x_{w}))\circ \theta_{1} = \exists T_{j}\rak{\theta_{1}(T_{j})\wedge \beta(T_{j})\wedge \theta_{3}(x_{w})}$. On the other hand

\[\theta_{3}  \circ  (*_{2}  \circ \theta_{1}) = \theta_{3}\circ *_{s s''} = s\wedge \theta_{3}(x_{w})\]
To show provable equivalence of these formulas, it suffices to obtain

\[
\Gamma, \neg E, s\wedge \theta_{3}(x_{w})\Rightarrow \exists T_{j}\rak{\theta_{1}(T_{j})\wedge \beta(T_{j})\wedge \theta_{3}(x_{w})} \tag{*a}
\]\[
\Gamma, \neg E, \exists T_{j}\rak{\theta_{1}(T_{j})\wedge \beta(T_{j})\wedge \theta_{3}(x_{w})} \Rightarrow s\wedge \theta_{3}(x_{w}) \tag{*b}
\]
By $T$-provability of $\mftwo^{\theta_{1}}$, we have \[
\Gamma, s\Rightarrow \exists T_{j} \theta_{1}(T_{j})\tag{1}
\]
By \frelim on $\mfone^{\theta_{1}}$, we have \[
\Gamma, \neg E, \theta_{1}(T_{j})\Rightarrow \beta(T_{j})
\]\[
\overline{\Gamma, \neg E, \theta_{1}(T_{j})\Rightarrow \theta_{1}(T_{j})\wedge \beta(T_{j})}
\]\[
\qquad \qquad \overline{\Gamma, \neg E, \exists T_{j} \theta_{1}(T_{j})\Rightarrow \exists T_{j}\rak{\theta_{1}(T_{j})\wedge \beta(T_{j})}} \qquad \qquad (\eD) \tag{2}
\]
and chaining (1) with this we have\[
\Gamma, s, \neg E \Rightarrow \exists T_{j}\rak{\theta_{1}(T_{j})\wedge \beta(T_{j})} \tag{3}
\]
By \id we have $\theta_{3}(x_{w})\Rightarrow \theta_{3}(x_{w})$, and \andd ing this with (3) we obtain\[
\Gamma, s, \neg E, \theta_{3}(x_{w}) \Rightarrow \exists T_{j}\rak{\theta_{1}(T_{j})\wedge \beta(T_{j})}\wedge\theta_{3}(x_{w}) \tag{4}
\]
$T_{j}$ are disjoint from $x_{w}$, so by \ins from (4) we obtain (*a).\\

Conversely, also by \frelim on $\mfone^{\theta_{1}}$ we have\[
\Gamma, \neg E, \theta_{1}(T_{j})\Rightarrow s
\]\[
\qquad \qquad \overline{\Gamma, \neg E, \theta_{1}(T_{j})\wedge \theta_{3}(x_{w})\Rightarrow s\wedge \theta_{3}(x_{w})}\qquad \qquad (\andd)
\]\[
\overline{\Gamma, \neg E, \beta(T_{j})\wedge \theta_{1}(T_{j})\wedge \theta_{3}(x_{w})\Rightarrow s\wedge \theta_{3}(x_{w})}
\]
and we get (*b) by \eA.\qed\\

1.2. Suppose $\alpha = \alpha(X_{i})$.\\

1.2.1. Suppose $\beta = s'$ is a sentence. Then our diagram is

\begin{center}
\begin{tikzcd}
\alpha(X_{i}) \arrow[r, "*_{1}"]& s' \arrow[r, "*_{2}"] & s'' \arrow[r, "\theta_{3}(x_{w})"] & \delta(P_{w})
\end{tikzcd}
\end{center}


We have $(\theta_{3}  \circ  *_{2})  \circ *_{1} = (s' \wedge \theta_{3}(x_{w}))\circ *_{1} = \alpha(x_{-i})\wedge s' \wedge \theta_{3}(x_{w})$. On the other hand

\[
\theta_{3}  \circ  (*_{2}  \circ *_{1}) = \theta_{3}\circ *_{\alpha s''} = \alpha(x_{-i}) \wedge \theta_{3}(x_{w})
\]

These are provably equivalent, since we have by \frelim on $\forall x_{-i}\rak{\alpha(x_{i})\rightarrow s'}$\[
\Gamma, \neg E \Rightarrow \alpha(x_{-i}) \rightarrow s' \qed
\]

1.2.2. Suppose $\beta = \beta(Y_{j})$. Then our diagram is

\begin{center}
\begin{tikzcd}
\alpha(X_{i}) \arrow[r, "{\theta_{1}(x_{-i}, x_{j})}"]& \beta(Y_{j}) \arrow[r, "*_{2}"] & s'' \arrow[r, "\theta_{3}(x_{w})"] & \delta(P_{w})
\end{tikzcd}
\end{center}

We have $(\theta_{3}  \circ  *_{2})  \circ \theta_{1} = (\beta(x_{-j})\wedge \theta_{3}(x_{w}))\circ \theta_{1} = \exists T_{j}\rak{\theta_{1}(x_{-i}, T_{j})\wedge \beta(T_{j})\wedge \theta_{3}(x_{w})}$. On the other hand

\[\theta_{3}  \circ  (*_{2}  \circ \theta_{1}) = *_{\alpha s''}\circ \theta_{3} = \alpha(x_{-i})\wedge \theta_{3}(x_{w})\]

The proof in this case is identical to 1.1.2; the addition of negative free variables in $\alpha(x_{-i})$ and $\theta_{1}(x_{-i}, x_{j})$ makes no difference for deriving the corresponding (*a) and (*b). \qed\\

2. Suppose $\gamma = \gamma(Z_{k})$.

2.1. Suppose $\beta = s'$ is a sentence.

2.1.1. Suppose $\alpha = s$ is a sentence. Then our diagram is

\begin{center}
\begin{tikzcd}
s \arrow[r, "*_{1}"]& s' \arrow[r, "\theta_{2}(x_{k})"] & \gamma(Z_{k}) \arrow[r, "{\theta_{3}(x_{-k}, x_{w})}"] & \delta(P_{w})
\end{tikzcd}
\end{center}

We have $(\theta_{3}  \circ  \theta_{2})  \circ *_{1} = (\exists T_{k}\rak{\theta_{2}(T_{k})\wedge \theta_{3}(T_{k}, x_{w})})\circ *_{1} = s\wedge \exists T_{k}\rak{\theta_{2}(T_{k})\wedge \theta_{3}(T_{k}, x_{w})}$. On the other hand

\[ \theta_{3}  \circ  (\theta_{2}  \circ *_{1}) = \theta_{3}\circ (s\wedge \theta_{2}(x_{k})) =  \exists T_{k}\rak{ s \wedge \theta_{2}(T_{k})\wedge \theta_{3}(T_{k}, x_{w})}\]

Provable equivalence of these is trivial by \ins, as $s$ is a sentence.\qed\\

2.1.2. Suppose $\alpha = \alpha(X_{i})$. Then our diagram is

\begin{center}
\begin{tikzcd}
\alpha(X_{i}) \arrow[r, "*_{1}"]& s' \arrow[r, "\theta_{2}(x_{k})"] & \gamma(Z_{k}) \arrow[r, "{\theta_{3}(x_{-k}, x_{w})}"] & \delta(P_{w})
\end{tikzcd}
\end{center}

We have $(\theta_{3}  \circ  \theta_{2})  \circ *_{1} = (\exists T_{k}\rak{\theta_{2}(T_{k})\wedge \theta_{3}(T_{k}, x_{w})}) \circ *_{1} = \alpha(x_{-i})\wedge \exists T_{k}\rak{\theta_{2}(T_{k})\wedge \theta_{3}(T_{k}, x_{w})}$. On the other hand

\[ \theta_{3}  \circ  (\theta_{2}  \circ *_{1}) = \theta_{3}\circ (\alpha(x_{-i})\wedge \theta_{2}(x_{k})) = \exists T_{k}\rak{\alpha(x_{-i}) \wedge \theta_{2}(T_{k})\wedge \theta_{3}(T_{k}, x_{w})}\]

$T_{k}$ does not occur free in $\alpha(x_{-i})$, so again we just use \ins.\qed\\

2.2. Suppose $\beta = \beta(Y_{j})$.

2.2.1 ($\alpha$ a nonsentence) is Case 0, which we have already shown. 2.2.2 ($\alpha = s$): the diagram is

\begin{center}
\begin{tikzcd}
s \arrow[r, "\theta_{1}(x_{j})"]& \beta(Y_{j}) \arrow[r, "{\theta_{2}(x_{-j}, x_{k})}"] & \gamma(Z_{k}) \arrow[r, "{\theta_{3}(x_{-k}, x_{w})}"] & \delta(P_{w})
\end{tikzcd}
\end{center}

There are no formal morphisms present, so the compositions are identical to in Case 0 except that $\theta_{1}$ is missing the free variables $x_{-i}$, which makes no difference to the proof. \qed\\

This concludes the proof of \cref{dec-assoc}. \qed\\

In the next part we establish a way to obtain certain morphisms we call {\it special morphisms}, which will be crucial to the rest of our constructions. As a corollary we will obtain identities and hence finish the construction of the category $\syn(T)$.

\chapter{Morphism Lemmas}\label{ch-lemmas}

\section{Special morphisms}
\newcommand{\tf}{\widetilde{f}}

\begin{clm}\label{dec-havespecial} Let $\psi(S_{i})$ and $\varphi(R_{j})$ be nonsentences, $\tf: \free(\varphi)\rightarrow \free(\psi)$ 
be a map, and suppose that $T$ proves the sentence
\[
\forall S_{i}\rak{\psi(S_{i})\rightarrow\varphi(\tf(R_{j}))}.
\]
Then if we let $f$ be the map of subsets of $\mathbb{N}$ such that $\tf(R_{j}) = S_{f(j)}$, the formula\[
S\tf(x_{-i}, x_{j}):= \psi(x_{-i})\wedge \bigwedge_{j}x_{j}= x_{-f(j)}
\]
defines a (type 1) morphism $\psi \to \varphi$, which we call the {\it $\tf$-special morphism $\psi \to \varphi$.} \end{clm}

\pf $S\tf$ has the right free variables, so we just need to show $T$-provability of the functionality sentences.

F1: \[\Gamma, \neg E, S\tf(x_{-i}, x_{j}) \Rightarrow \psi(x_{-i})\wedge \varphi(x_{j})\tag{*}
\]
We trivially have $S\tf(x_{-i}, x_{j})\Rightarrow \psi(x_{-i})$. By \frelim on the supposition (instantiating $S_{i}$ as $x_{-i}$), we have\[
\Gamma, \neg E, \psi(x_{-i})\Rightarrow \varphi(x_{-f(j)})\tag{1}
\]\[
\overline{\Gamma, \neg E, \psi(x_{-i})\Rightarrow \psi(x_{-i}) \wedge \varphi(x_{-f(j)}) }\tag{1\pr}
\]

(the instantiation on $\varphi$ works out as $\varphi(\tf(R_{j}))\sbb{x_{-i}}{S_{i}} := \varphi\sbb{S_{f(j)}}{R_{j}}\sbb{x_{-i}}{S_{i}} =  \varphi\sbb{x_{-f(j)}}{R_{j}}=: \varphi(x_{-f(j)})$. We will not rigorously justify the middle equals sign here, but it is intuitively clear.) 
Chaining the trivial sequent with (1\pr) we obtain\[
\Gamma, \neg E, S\tf(x_{-i}, x_{j}) \Rightarrow \psi(x_{-i}) \wedge \varphi(x_{-f(j)}) \tag{2}
\]

from which (*) follows by manipulations with \subs (using the definition of $S\tf$).\qed

F2: \[
\Gamma, \neg E, \psi(x_{-i})\Rightarrow \exists x_{j} S\tf(x_{-i},x_{j})\tag{*}\]

The existential claim has an ``obvious witness''. Rigorously, $\psi(x_{-i})\wedge \bigwedge_{j} x_{-f(j)} = x_{-f(j)}$ is the same formula as $S\tf \sbb{x_{-f(j)}}{x_{j}}$, so we have \[
\psi(x_{-i}), \bigwedge_{j} x_{-f(j)} = x_{-f(j)} \Rightarrow \psi(x_{-i})\wedge \bigwedge_{j} x_{-f(j)} = x_{-f(j)}  \qquad \qquad (\assm) \tag{1}
\]\[
\overline{\psi(x_{-i}), \bigwedge_{j} x_{-f(j)} = x_{-f(j)}, \neg E \Rightarrow \exists x_{j} S\tf} \qquad \qquad (\eS) \tag{2}
\]
We trivially have $\Rightarrow \bigwedge_{j} x_{-f(j)} = x_{-f(j)}$ by \id, and chaining this with (2) and adding in the extra hypotheses we obtain (*).\qed

F3: \[
\Gamma, \neg E, S\tf(x_{-i}, x_{j})\wedge S\tf(x_{-i}, x_{n + j}) \Rightarrow \bigwedge_{j}x_{j} = x_{n + j}\tag{*}
\]
This is derivable by transitivity of equality (a consequence of the rule \subs) since $S\tf$ asserts all output (positive) variables equal to some negative variables. \qed\\

This shows \cref{dec-havespecial}. \qed\\

\begin{prop}\label{dec-soutgoing} Let $\theta: \varphi(R_{j})\to \rho(W_{k})$ be type 1. Then we have\[
\theta\circ  S\tf= \psi(x_{-i}) \wedge \theta(x_{-f(j)}, x_{k}).
\]\end{prop}

\pf The composition is defined as $\exists T_{j}\rak{\psi(x_{-i})\wedge \bigwedge_{j}T_{j} = x_{-f(j)}\wedge \theta(T_{j}, x_k)}$. It suffices to show that this formula is provably equivalent to the claimed one, since then $\Gamma\cup \{\neg E\}$ must prove the \frelim s for functionality sentences for the claimed formula, proving them for this composition. We will thus conclude the claimed formula is a premorphism, as usual by \frlS and \pc. And then by the definition of our equivalence relation on Prehom it must represent the same morphism as the composition. This is what is meant by the ``equation'' of formulas.

Indeed, by the Equivalence Lemma and Subs, $\theta \circ S\tf$ is clearly provably equivalent to \[
\exists T_{j}\rak{\psi(x_{-i})\wedge \theta(x_{-f(j)}, x_k) \wedge \bigwedge_{j}T_{j} = x_{-f(j)}}\tag{1}\]
By \ins/ \ex (1) is provably equivalent to \[
\psi(x_{-i})\wedge \theta(x_{-f(j)}, x_k) \wedge \exists T_{j}\rak{\bigwedge_{j}T_{j} = x_{-f(j)}}\tag{1\pr}
\]
and this is provably equivalent to the claimed formula, as the last clause is vacuous (we have the $\neg E$ hypothesis).\qed\\

\newcommand{\ran}{\text{ran}}

\begin{prop}\label{dec-sincoming} Suppose $\tf$ is injective, and let $\theta: \rho(W_{k})\to \psi(S_{i})$ be type 1. Let $m_{1}, \dots$ be the natural numbers $\leq n$ which are {\it not} in $\ran(f)$, i.e. $\tf(R_{j})\neq S_{m_{l}}$ for any $j$, $l$. Then defining \[
Q_{i} = \begin{cases} x_{f^{-1}(i)} & i \in \ran(f) \\ T_{m_{l}} & i = m_{l}\end{cases}
\]we have \[
S\tf \circ \theta = \exists T_{m_{1}}\exists T_{m_{2}}\dots \theta(x_{-k}, Q_{i}).
\]
If there are no such $m_{l}$ ($f, \tf$ are bijective) then $Q_{i} = x_{f^{-1}(i)}$ for all $i$ and there is no quantification, i.e. we have $S\tf \circ \theta = \theta(x_{-k}, x_{f^{-1}(i)})$.\end{prop}

\pf The composition is defined as $S\tf \circ \theta = \exists T_{i}\rak{\theta(x_{-k}, T_{i})\wedge S\tf(T_{i}, x_{j})}$\[
= \exists T_{i}\rak{\theta(x_{-k}, T_{i})\wedge \psi(T_{i})\wedge\bigwedge_{j }x_{j}=T_{f(j)}}. \tag{1}
\]

We first consider the subformula\[
\theta(x_{-k}, T_{i})\wedge \psi(T_{i})\wedge\bigwedge_{j }x_{j}=T_{f(j)} \tag{2}
\]
Since by $\fone^{\theta}$ we know that $\Gamma \cup \{\neg E\}$ proves $\theta(x_{-k}, T_{i})\rightarrow \psi(T_{i})$, (2) is provably equivalent to \[
\theta(x_{-k}, T_{i})\wedge \bigwedge_{j }x_{j}=T_{f(j)}. \tag{2\pr}
\]

Clearly by \subs 
(2\pr) is provably equivalent to\[
\theta(x_{-k}, T_{i})\sbb{x_{1}}{T_{f(1)}}\dots\sbb{x_{m}}{T_{f(m)}} \wedge\bigwedge_{j }x_{j}=T_{f(j)}
\]\[
=\theta(x_{-k}, Q_{i}) \wedge\bigwedge_{j }x_{j}=T_{f(j)}.\tag{3}
\]

Therefore by the Equivalence Lemma (1) is provably equivalent to\[
\exists T_i\rak{\theta(x_{-k}, Q_{i}) \wedge\bigwedge_{j }x_{j}=T_{f(j)}}
\]
p.e. to \[
\exists T_{m_{1}}\exists T_{m_{2}}\dots \exists T_{f(1)}\dots \exists T_{f(m)} \rak{\theta(x_{-k}, Q_{i}) \wedge\bigwedge_{j }x_{j}=T_{f(j)}}\tag{4}
\]
by \reord. Now (4) has a subformula $\exists T_{f(1)}\dots \exists T_{f(m)} \rak{\theta(x_{-k}, Q_{i}) \wedge\bigwedge_{j }x_{j}=T_{f(j)}}$, in which no $Q_{i}$ is quantified; by \ins and \ex this subformula is provably equivalent to \[
\theta(x_{-k}, Q_{i}) \wedge \exists T_{f(1)}\dots \exists T_{f(m)} \rak{\bigwedge_{j }x_{j}=T_{f(j)}},
\]
which is p.e. to just \[\theta(x_{-k}, Q_{i}),\] since the last clause is vacuous.

Therefore by the Equivalence Lemma again, (4) is provably equivalent to\[
\exists T_{m_{1}} \exists T_{m_{2}}\dots \rak{\theta(x_{-k}, Q_{i})}
\]
and by transitivity of provable equivalence 
this formula is p.e. to (1) and hence a representative for the morphism $S\tf \circ \theta$, as claimed.

The argument needs no modification in the case of bijectivity, where there are just no $\exists T_{m_{l}}$'s sitting out in front. We also note that it goes through identically if $\theta$ is assumed type 3, i.e. $\rho = s$ and $\theta$ does not contain any negative free variables $x_{-k}$.\qed\\

\newcommand{\ti}{\widetilde{\iota}}
\section{Identities}
\subsection{Nonsentences}
For any nonsentence $\varphi(X_{i})$, letting $\ti: \free(\varphi)\to \free(\varphi)$ be the identity map, we clearly have an $\ti$-special morphism $I_{\varphi}: \varphi \to \varphi$. The associated index-function $\iota$ is also the identity on $\crul{1, \dots |\varphi|}$.

Now if $\theta: \gamma(Z_{k}) \to \varphi$ is any type 1 morphism, then by \cref{dec-sincoming} we have \[
I_\varphi \circ \theta = \theta(x_{-k}, x_i) = \theta
\]
as morphisms, and similarly $I_{\varphi}\circ \theta = \theta$ if $\theta$ is type 3. Therefore, $I_{\varphi}$ acts as an identity for postcomposing.\\

By \cref{dec-soutgoing}, if $\theta: \varphi \to \gamma(Z_{k})$ is type 1, we have \[
\theta\circ I_{\varphi} = \varphi(x_{-i})\wedge \theta(x_{-i}, x_{k})
\]
and this represents the same morphism $\varphi\to \gamma$ as the formula $\theta(x_{-i}, x_{k})$, since we can already conclude $\varphi(x_{-i})$ from $\theta(x_{-i}, x_{k})$, by $\fone^{\theta}$. Therefore composition with $I_{\varphi}$ also acts as an identity in this case.

Finally, if we have a type 2 morphism $\varphi \to s$, precomposition with $I_{\varphi}$ again does nothing, since there is a unique morphism $\varphi \to s$.\\

Therefore $I_{\varphi}$ is an identity morphism for $\varphi$ in the categorical sense. We note that $I_{\varphi}$ is given by the formula\[
\varphi(x_{-i})\wedge \bigwedge_{i} x_{i} = x_{-i}.
\]

\subsection{Sentences}
A sentence $s$ has a unique morphism to itself by construction, and postcomposing a map into $s$ with it does nothing, by the formal composition rules. Similarly for precomposing $s \to s$ with $s \to s'$. Finally, if we have a type 3 morphism $\theta: s \to \gamma(Z_{k})$, by definition the precomposition with $s \to s$ is type 3 represented by\[
s\wedge \theta(x_{k}).
\]
This is the same morphism as $\theta: s\to \gamma(Z_{k})$, because we already conclude $s$ from $\theta(x_{k})$ by $\mfone^{\theta}$.

Therefore $s \to s$ is a categorical identity for $s$.\\

This completes our construction of the (weak) syntactic category $\syn(T)$. \qed

\section{Terminology}
In general if $\ti:\free(\varphi)\to \free(\psi)$ is the identity map, an $\ti$-special morphism will be called a {\it special monomorphism}. The index function $\iota$ is also the identity of corresponding subsets of $\mathbb{N}$ in this case. Thus \cref{dec-sincoming} says that the precomposition with a morphism represented by $\theta$ is again represented by $\theta$. It follows that $S\ti$ is indeed a monomorphism. 

We will also call any map of sentences a special monomorphism. It is a mono (as is also any type 3 morphism) since by construction our category does not contain two parallel maps with sentence target.

A special monomorphism which is an iso will be called a special isomorphism.

\section{Isomorphisms}
\begin{lem}[Extraction]\label{dec-exlem}Let $\theta(x_{-i}, x_j)$ be $T$-functional (i.e. be a type 1 premorphism). Then for any formula $\gamma$, variables $(P_i)$ and $(Q_j)$, distinct variables $(K_{j})$, and new distinct variables $(W_j)$, we claim derivability of\[
\Gamma, \neg E \Rightarrow \theta(P_i, Q_j) \wedge \exists W_j \rak{\theta(P_i, W_j)\wedge \gamma \sbb{W_{j}}{K_{j}}} \rightarrow \gamma\sbb{Q_{j}}{K_{j}}.
\] \end{lem}

\pf We have \[
\Gamma, \neg E, \theta(P_i, Q_j), \theta(P_i, W_j) \Rightarrow \bigwedge_{j} Q_j = W_j \tag{1}\]
by \frelim on $\fth^{\theta}$.

Also by \subs we have \[
\gamma \sbb{W_{j}}{K_{j}}, \bigwedge_{j} Q_{j} = W_{j} \Rightarrow \gamma \sbb{W_{j}}{K_{j}}\sbb{Q_{j}}{W_{j}}\tag{2}
\]
and the RHS is $\gamma \sbb{Q_j}{K_j}$ as the $W_j$ are new. Therefore chaining (1) and (2) we have
\[
\Gamma, \neg E, \theta(P_i, Q_j), \theta(P_i, W_j), \gamma \sbb{W_{j}}{K_{j}} \Rightarrow \gamma \sbb{Q_j}{K_j}
\]
Again as $W_{j}$ are new by $\exists$A this gives
\[
\Gamma, \neg E, \theta(P_i, Q_j), \exists W_{j}\rak{\theta(P_i, W_j) \wedge \gamma \sbb{W_{j}}{K_{j}}} \Rightarrow \gamma \sbb{Q_j}{K_j}
\]
which gives the claimed sequent.\qed\\

\begin{prop}[Type 1 isomorphism condition, sufficient] \label{dec-t1isosuf} Let $\theta: \alpha(X_{i})\to \beta(Y_{j})$ be Type 1, and suppose the negated-variables formula $\widetilde{\theta}:=\theta\sbx{x_{i}}{x_{-i}}{x_{-j}}{x_{j}}=\theta(x_{i}, x_{-j})$ is a premorphism $\beta \to \alpha$. Then $\theta$ represents an isomorphism, and the inverse is represented by $\widetilde{\theta}$.\end{prop}

\pf Suppose $\widetilde{\theta}:\beta \rightarrow \alpha$ is a premorphism. The task is to establish $\widetilde{\theta}\circ \theta = I_{\alpha}$ and $\theta \circ \widetilde{\theta} = I_{\beta}$ as morphisms. But we observe that negating variables again (i.e. forming $\widetilde{\widetilde{\theta}}$) yields $\theta$ (literally, not just up to provable equivalence) and this is also a premorphism. It follows that we need only establish the first equality, and the second will also follow, by application of the argument to $\theta' = \widetilde{\theta}$.

For any $(P_{j})$, $(Q_{i})$, we have $\widetilde{\theta}(P_{j}, Q_{i}) \overset{\text{def}}{=} \theta\sbx{x_{i}}{x_{-i}}{x_{-j}}{x_{j}}\sbx{P_{j}}{x_{-j}}{Q_{i}}{x_{i}}$\[
=\theta\sbx{P_{j}}{x_{j}}{Q_{i}}{x_{-i}} \overset{\text{def}}{=} \theta(Q_{i}, P_{j}).
\]

Therefore $\widetilde{\theta} \circ \theta=\exists T_{j}\rak{\theta(x_{-i}, T_{j})\wedge \widetilde{\theta}(T_{j}, x_{i})}=\exists T_{j}\rak{\theta(x_{-i}, T_{j}) \wedge \theta(x_{i}, T_{j})}$. Therefore we need to obtain\[
\Gamma, \neg E, \exists T_{j}\rak{\theta(x_{-i}, T_{j}) \wedge \theta(x_{i}, T_{j})} \Rightarrow \alpha(x_{-i}) \wedge \bigwedge_{i} x_{i}=x_{-i} \tag{*a}
\]\[
\Gamma, \neg E, \alpha(x_{-i}) \wedge \bigwedge_{i} x_{i}=x_{-i} \Rightarrow \exists T_{j}\rak{\theta(x_{-i}, T_{j}) \wedge \theta(x_{i}, T_{j})} \tag{*b}
\]

By the above we have \[\fth^{\widetilde{\theta}}
=\forall x_{-j} \forall x_{i} \forall x_{n+i}\rak{\widetilde{\theta}(x_{-j}, x_{i})\wedge\widetilde{\theta}(x_{-j}, x_{n+i}) \rightarrow \bigwedge_{i} x_{i}=x_{n+i}} \]
\[=\forall x_{-j} \forall x_{i} \forall x_{n+i}\rak{\theta(x_{i}, x_{-j})\wedge\theta(x_{n+i}, {x_{-j}}) \rightarrow \bigwedge_{i} x_{i}=x_{n+i}}. \]

By \frelim on this we have\[
\Gamma, \neg E, \theta(x_{-i}, T_{j}) \wedge \theta(x_{i}, T_{j}) \Rightarrow \bigwedge_{i} x_{-i}=x_{i}.\tag{1}
\]

By $\fone^{\theta}$ we have \[
\Gamma, \neg E, \theta(x_{-i}, T_{j})\Rightarrow \alpha(x_{-i}) \tag{2}
\]
Applying $\andd $ and then $\eA$, we obtain (*a).\\

Conversely, by $\ftwo^{\theta}$, we have\[
\Gamma, \neg E, \alpha(x_{-i}) \Rightarrow \exists T_{j}\rak{\theta(x_{-i}, T_{j})} \tag{3}
\]

An instance of \subs is \[
\Gamma, \neg E, \theta(x_{-i}, T_{j}), \bigwedge_{i} x_{i} = x_{-i}\Rightarrow \theta(x_{i}, T_{j})
\]\[
\overline{\Gamma, \neg E, \theta(x_{-i}, T_{j}), \bigwedge_{i} x_{i} = x_{-i}\Rightarrow\theta(x_{-i}, T_{j}) \wedge\theta(x_{i}, T_{j})}
\]\[
\qquad \overline{\Gamma, \neg E, \exists T_{j} \rak{\theta(x_{-i}, T_{j})}, \bigwedge_{i} x_{i} = x_{-i}\Rightarrow\exists T_{j}\rak{\theta(x_{-i}, T_{j}) \wedge\theta(x_{i}, T_{j})}}\qquad  (\eD) \tag{4}
\]
and chaining (3) and (4) we obtain (*b).\qed\\

\begin{prop}[Type 1 isomorphism condition, necessary] \label{dec-t1isonec} Let $\theta: \alpha(X_{i}) \rightarrow \beta(Y_{j})$ represent an isomorphism. Then its inverse is represented by $\widetilde{\theta}$.\end{prop}

\pf Let $\iota(x_{-j}, x_{i})$ represent the inverse of $\theta$. Similarly to in the proof of \cref{dec-soutgoing}, it suffices to show $\widetilde{\theta}$ and $\iota$ are provably equivalent, i.e. to obtain \[
\Gamma, \neg E, \widetilde{\theta} \Rightarrow \iota(x_{-j}, x_{i}) \tag{*a}
\]
\[
\Gamma, \neg E, \iota({x}_{-j} ,x_{i}) \Rightarrow \widetilde{\theta}, \tag{*b}\]


By \frelim on $\E^{I_{\alpha}\rightarrow\iota\circ\theta}$ (instantiating $x_{-i}$, $x_{i}$ both as $x_{i}$) we have \[
\Gamma,\neg E, \alpha(x_{i}) \wedge \bigwedge_{i} x_{i}=x_{i} \Rightarrow \exists M_{j}\rak{\theta(x_{i}, M_{j})\wedge\iota(M_{j},x_{i})} \tag{1}
\]

By $\fone^{\theta}$ we have
\[\Gamma,\neg E, \theta(x_{i}, x_{-j}) \Rightarrow \alpha(x_{i})\]\[
\overline{\Gamma,\neg E, \theta(x_{i}, x_{-j}) \Rightarrow \alpha(x_{i}) \wedge\bigwedge_{i} x_{i}=x_{i}}\tag{2}
\]

Chaining (2) and (1) we have \[
\Gamma, \neg E, \theta(x_{i}, x_{-j}) \Rightarrow \exists M_{j}\rak{\theta(x_{i}, M_{j})\wedge\iota(M_{j}, x_{i})}\tag{3}
\]\[
\overline{\Gamma, \neg E, \theta(x_{i}, x_{-j}) \Rightarrow\theta(x_{i}, x_{-j})\wedge \exists M_{j}\rak{\theta(x_{i}, M_{j})\wedge\iota(M_{j}, x_{i})}}\tag{3\pr}
\]
Now as $\theta$ is functional, (3\pr) can be chained with the output of the Extraction Lemma (taking the formula $\gamma$ there to be $\iota$) to give\[
\Gamma, \neg E,\theta(x_{i}, x_{-j}) \Rightarrow\iota(x_{-j}, x_{i})
\]
This is (*a) by the definition of $\widetilde{\theta}$.

Conversely, for (*b), we apply the same trick in the opposite direction. Namely, by \frelim on $\E^{I_{\beta} \rightarrow \theta\circ \iota}$ we have
\[
\Gamma,\neg E, \beta(x_{-j}) \wedge \bigwedge_{j} x_{-j}=x_{-j} \Rightarrow \exists T_{i}\rak{\iota(x_{-j}, T_{i})\wedge\iota(T_{i},x_{-j})} \tag{1}
\]
We have $\Gamma, \neg E, \iota(x_{-j}, x_{i}) \Rightarrow \beta(x_{-j})$ by $\fone^{\iota}$ and thus\[
 \Gamma,\neg E, \iota(x_{-j},x_{i})\Rightarrow\beta(x_{-j})\wedge\bigwedge_{j} x_{-j}=x_{-j}\tag{2}
 \]
 Chaining (2) and (1) we have \[
 \Gamma,\neg E, \iota(x_{-j}, x_{i}) \Rightarrow \exists M_{i}\rak{\iota (x_{-j}, M_{i}) \wedge \theta(M_{i}, x_{-j})}
 \]\[
 \overline{\Gamma,\neg E, \iota(x_{-j}, x_{i}) \Rightarrow \iota(x_{-j}, x_{i})\wedge \exists M_{i}\rak{\iota (x_{-j}, M_{i}) \wedge \theta(M_{i}, x_{-j})}} \tag{3}
 \]
whence\[
\Gamma, \neg E,\iota(x_{-j}, x_{i}) \Rightarrow \theta(x_{i}, x_{-j})
\] by Extraction as $\iota$ is functional. This is (*b). \qed\\

\begin{rk}\label{dec-specialisonice} In the case of a special isomorphism of nonsentences, the postcomposition formula \cref{dec-soutgoing} is \[
\theta \circ S\ti  = \psi(x_{-i}) \wedge \theta(x_{-i}, x_{k})
\] and we in fact may drop the $\psi(x_{-i})$ clause. This is because by \cref{dec-t1isonec}, we can represent $(S\ti)^{-1}$ by the negated-variables formula\[
\psi(x_{i})\wedge \bigwedge_{i}x_{-i} = x_{i}.
\]
By $\fone^{\theta}$, from $\theta(x_{-i}, x_{k})$ alone we can conclude $\varphi(x_{-i})$, and hence $\exists x_{i} (S\ti)^{-1}$ by $\ftwo^{(S\ti)^{-1}}$. But then given the form of $(S\ti)^{-1}$ above we clearly have $\psi(x_{-i})$, essentially by the same argument as when proving the Extraction Lemma.\qed \\

The effect is that composing with a special {\it iso} on {\it either} side acts trivially on the representing formula (while only {\it postcomposition with} a special {\it mono} in general acts trivially). Simpler reasoning than the above shows this for special isomorphisms of sentences, without appealing to the foregoing isomorphism characterizations.

Thus special isomorphic objects are ``the same'' in a sense beyond that of simple isomorphism: the natural correspondence of their maps {\it preserves any associated defining formulas}.\end{rk}

\bigskip

\begin{rk} \label{dec-specialisosare}
It follows easily from the above that two objects are special isomorphic iff \begin{itemize}
\item they are $T\cup \crul{\neg E}$-provably equivalent nonsentences in exactly the same free variables, or
\item they are $T$-provably equivalent sentences.
\end{itemize} \end{rk}

\chapter{Prospect}\label{ch-prospect}

To effect the reduction of G\"{o}del's Theorem to Theorem 4 in (\cite{lu}, Lecture 6), we need to show first of all that the category $\syn(T)$ we have constructed is a {\it small consistent coherent category}. This will be the content of Chapters 5-8. We give the relevant definitions from (\cite{lu}, Lecture 5) below.\\

\begin{defn}
A category $\mathcal{C}$ is {\it coherent} if it satisfies the following.\begin{enumerate}
\item $\mathcal{C}$ admits finite limits, including a final object
\item Every morphism in $\mathcal{C}$ factors into an {\it effective epimorphism} followed by a monomorphism
\item For every $X\in \text{Ob}(\mathcal{C})$, the natural poset $\sub(X)$ of subobjects of $X$ is in fact a {\it join semilattice}
\item The pullback of an effective epimorphism (sometimes abbreviated e.e.) along any morphism in $\mathcal{C}$, is again an effective epimorphism 
\item For every morphism $X \to Y$ in $\mathcal{C}$, the induced {\it pullback map} $\sub(Y)\to \sub(X)$ is a homomorphism of join semilattices.
\end{enumerate}
\end{defn}

\bigskip

And, a coherent $\mathcal{C}$ is {\it consistent} (see \cite{lu}, Lecture 6) if there is a nontrivial subobject of the final object.\\

\begin{rk}\label{dec-smol} Smallness of our category $\syn(T)$ is immediate since it certainly has no more morphisms than there are formulas of $L$, and these form a set.\end{rk}
\bigskip

\begin{rk}\label{dec-final} It is also immediate that $\syn(T)$ admits a final object: the validity $S: = \forall x \rak{x = x}$. Indeed, we have

\[
\overline{\Rightarrow x_{0} = x_{0}} \qquad \qquad (\id)
\]\[
\overline{\Rightarrow \forall x_{0}\rak{x_{0} = x_{0}}} \qquad \qquad (\frlS).
\]

Thus, for any nonsentence $\varphi$, $T$ certainly proves $\forall \free(\varphi)\rak{\varphi \to S}$ (it proves $S$ unconditionally, and we can add $\varphi$ as an antecedent by \ant, then apply \meta and \frlS). We thus have a (necessarily unique) formal morphism $\varphi \to S$. We also clearly have unique $s\to S$ for any sentence $s$.\qed \end{rk}
\bigskip

\begin{rk}\label{dec-consistent}
Further, if we assume $\syn(T)$ satisfies the properties 1-5, then its consistency as a coherent category is immediate from the assumed consistency of $T$ as a first-order theory. Indeed, $\sub(S)$ being trivial ($=$containing only a single element, necessarily the equivalence class of $S \to S$ itself) would mean that $T$ proved every sentence.

For any $s$ is an obvious subobject of $S$, and since this would be equal to $S\to S$, we would have that $s$ and $S$ are isomorphic. Then according to our terminology in the last Chapter they must be special isomorphic. Then by \cref{dec-specialisosare} there, $s$ and $S$ are $T$-provably equivalent sentences, but $S$ is $\emptyset$-provable, so this just means that $T$ proves $s$.\qed
\end{rk}
\bigskip

In light of these remarks, our goal for chapters 5-8 is to verify the coherence properties 1-5 for $\syn(T)$. We shall imitate the development of \cite{lu}, Lectures 1-5. That is, we shall\begin{enumerate}[label = (\alph*)]
\item establish existence of finite limits. To do this, in light of \cref{dec-final} it suffices by abstract nonsense to obtain binary pullbacks. We construct these in the next Chapter.
\item construct the morphism factorization for property 2. This is Chapter 6.
\item analyze $\sub(X)$ for general $X$, showing, as in \cite{lu}, that it is in fact a Boolean lattice. We shall obtain an explicit description of the pullback map which lets us check that it is a Boolean lattice homomorphism. This takes care of properties 3 and 5, and is done in Chapter 7.
\item use part (b) to get an easy characterization of e.e.'s, which will allow us to verify their stability under pullback (property 4), using our explicit description of the map. This is Chapter 8.
\end{enumerate}

\chapter{Fiber Products}\label{ch-fps}

We establish the existence of fiber products in $\syn(T)$ in this part. We do it first in Case 0, i.e. we construct a pullback of two type 1 morphisms with common target.

We will always set things up such that all bound variables are obvious (i.e. we never quantify over variables which may occur in other places than where they are explicitly written).

\newcommand{\tf}{\widetilde{f}}
\section{Case 0}
Let $\theta_{\alpha}: \alpha(X_{i}) \rightarrow \gamma(Z_{k})$ and $\theta_{\beta}: \beta(Y_{j}) \rightarrow \gamma(Z_{k})$ be type 1, and without loss of generality assume the representing premorphisms actually contain all possible free variables. Then we define the formula
$$\FP=\exists W_{k}\rak{\theta_{\alpha}(T_{i}, W_{k})\wedge \theta_\beta (T_{n+ j}, W_{k})}$$
where $(T_{s})_{1\leq s \leq n + m}$ is an {\it arbitrarily} chosen sequence 
increasing in index, and $(W_{k})$ is any sequence of distinct variables distinct from all $T_{s}$.\\

We let $\tf_{\alpha}:\free(\alpha)\to \free(\FP)$ be the injection which sends $X_{i}$ to $T_{i}$. The index function $f_{\alpha}$ is $\crul{1, \dots n}\hookrightarrow \crul{1, \dots m +n}$, by the choice of $(T_{s})$ to be increasing. We have $\alpha(\tf_{\alpha}(X_{i})) = \alpha(T_{i})$, and clearly by $\fone^{\theta_{\alpha}}$ we have $T$-derivability of\[
\forall T_{s} \rak{\FP \to \alpha(T_{i})}.
\]
Therefore we have an $\tf_{\alpha}$-special morphism $\FP \to \alpha$, which we call $\pi_{\alpha}$.\\

Similarly letting $Y_{j}\overset{\tf_{\beta}}{\mapsto} T_{n + j}$, we have $\tf_{\beta}$-special $\pi_{\beta}: \FP \to \beta$.\\

By \cref{dec-soutgoing} of Morphism Lemmas we have\[
\theta_{\alpha}\circ \pi_{\alpha} = \FP(x_{-s})\wedge \theta_{\alpha}(x_{-f_{\alpha}(i)}, x_{k}) = \FP(x_{-s})\wedge \theta_{\alpha}(x_{-i}, x_{k})
\]
\[
=\FP\sbx{x_{-i}}{T_{i}}{x_{-n-j}}{T_{n+j}}\wedge\theta_{\alpha}(x_{-i}, x_{k}).
\]
This is equal to\[
\exists W_{k}\rak{\theta_{\alpha}(x_{-i}, W_{k})\wedge \theta_\beta (x_{-n -j}, W_{k})}\wedge \theta_{\alpha}(x_{-i}, x_{k}).
\]
For $\theta_{\beta}\circ \pi_{\beta}$, the index map is $j \overset{f_{\beta}}{\mapsto} n + j$, so by \cref{dec-soutgoing} \[
\theta_{\beta}\circ \pi_{\beta} = \exists W_{k}\rak{\theta_{\alpha}(x_{-i}, W_{k})\wedge \theta_\beta (x_{-n -j}, W_{k})}\wedge \theta_{\beta}(x_{-n-j}, x_{k})
\]
By the Extraction Lemma these are easily seen to be provably equivalent, hence they represent the same morphism $\FP \to \gamma$ and we have a commutative square

\begin{center}\begin{tikzcd}
\FP \arrow[r, "\pi_{\alpha}"]\arrow[d, "\pi_{\beta}"] & \alpha(X_{i}) \arrow[d, "\theta_{\alpha}"]\\
\beta(Y_{j})\arrow[r, "\theta_{\beta}"] & \gamma(Z_{k}).
\end{tikzcd}\end{center}\bigskip

We shall establish that this square is a pullback, by checking the universal property.\\

So suppose we are given a commutative square \begin{center}\begin{tikzcd}
\delta \arrow[r, "\tau_{\alpha}"]\arrow[d, "\tau_{\beta}"] & \alpha(X_{i}) \arrow[d, "\theta_{\alpha}"]\\
\beta(Y_{j})\arrow[r, "\theta_{\beta}"] & \gamma.
\end{tikzcd}\end{center}
We must show existence and uniqueness of a solution $\eta: \delta \to \FP$ to the factoring problem.

\section{Subcase 0.0}
We first suppose that $\delta = \delta(S_q)$ is a nonsentence. Then we claim \[\eta:=\tau_{\alpha}(x_{-q}, x_{i}) \wedge \tau_{\beta}(x_{-q}, x_{n+j})\] works, i.e. this formula represents the unique map $\delta \rightarrow \FP$ which factors $\tau_{\alpha}$ through $\pi_{\alpha}$ and $\tau_{\beta}$ through $\pi_{\beta}$. First we must check it represents a map $\delta \rightarrow \FP$.\\

Indeed, $\eta$ has the right free variables. Modifying our previous notation in two special cases, for any sequences $(A_{q})$, $(B_{i})$, and $(C_{j})$ we shall denote the simultaneous substitution of all $A_{q}$ for $x_{-q}$, $B_{i}$ for $x_{i}$, and $C_{j}$ for $x_{n + j}$ in $\eta$ by $\eta(A_{q}, B_{i}, C_{j})$. Similarly $\FP(B_{i}, C_{j})$ shall denote the simultaneous substitution of all $B_{i}$ for $T_{i}$ and $C_{j}$ for $T_{n + j}$ in $\FP$.\\

$\fone^{\eta}$: it suffices to derive
\[\Gamma, \neg E, \eta(x_{-q}, x_{i}, x_{n + j}) \Rightarrow  \delta(x_{-q}) \wedge \FP(x_{i}, x_{n + j})\tag{*}\]

By \frelim on $\fone^{\tau_{\alpha}}$, we have \[\Gamma, \neg E, \tau_{\alpha}(x_{-q}, x_{i}) \Rightarrow \delta(x_{-q})\]\[
\qquad \qquad \overline{\Gamma, \neg E, \eta(x_{-q}, x_{i}, x_{n + j}) \Rightarrow \delta(x_{-q})} \qquad \qquad (\ant) \tag{1}.
\]

Hence by \andd it will suffice to get \[
\Gamma, \neg E, \eta \Rightarrow \FP(x_{i}, x_{n + j})  
\]\[
=\Gamma, \neg E, \tau_{\alpha}(x_{-q}, x_{i})\wedge \tau_{\beta}(x_{-q}, x_{n + j}) \Rightarrow \exists P_{k}\rak{\theta_{\alpha}(x_{-i},P_{k})\wedge  \theta_{\beta}(x_{n + j}, P_{k})}\tag{**}    
\] by alpha-equivalence, where $(P_{k})$ are new and distinct, in particular not overlapping with $(x_{-q})$.\\

By \frelim on $\ftwo^{\theta_{\alpha} \circ \tau_{\alpha}}$, we have\[
\Gamma, \neg E, \delta(x_{-q})\Rightarrow \exists P_{k}\rak{(\theta_{\alpha}\circ\tau_{\alpha})(x_{-q,}, P_{k})}\tag{2}\]
\[(1) (2)\]
\[\qquad \qquad \overline{\Gamma, \neg E, \eta(x_{-q}, x_{i}, x_{n + j}) \Rightarrow \exists P_{k}\rak{(\theta_{\alpha} \circ \tau_{\alpha})(x_{-q}, P_{k})}} \qquad \qquad (\ch)\tag{3}\]

Now by \frelim on $\E^{\theta_{\alpha} \circ \tau_{\alpha} \rightarrow \theta_{\beta}\circ \tau_{\beta}}$ we have

$$\Gamma, \neg E, (\theta_{\alpha} \circ \tau_{\alpha})(x_{-q}, P_{k}) \Rightarrow (\theta_{\beta} \circ \tau_{\beta})(x_{-q}, P_{k})$$
\[\overline{\Gamma, \neg E, (\theta_{\alpha}\circ \tau_{\alpha})(x_{-q}, P_{k}) \Rightarrow (\theta_{\alpha}{\circ} \tau_{\alpha})(x_{-q}, P_{k}) \wedge (\theta_{\beta} \circ\tau_{\beta})(x_{-q}, P_{k})}
\]\[
\overline{\Gamma, \neg E, \exists P_{k}\rak{(\theta_{\alpha}\circ \tau_{\alpha})(x_{-q}, P_{k})} \Rightarrow \exists P_{k}\rak{(\theta_{\alpha}{\circ} \tau_{\alpha})(x_{-q}, P_{k}) \wedge (\theta_{\beta} \circ\tau_{\beta})(x_{-q}, P_{k})}}\hspace{2mm} (\eD)\tag{4}
\]

\[(3)(4)
\]\[
\overline{\Gamma, \neg E, \eta(x_{-q}, x_{i}, x_{n + j})\Rightarrow \exists P_{k}\rak{(\theta_{\alpha}{\circ} \tau_{\alpha})(x_{-q}, P_{k}) \wedge (\theta_{\beta} \circ\tau_{\beta})(x_{-q}, P_{k})}}\hspace{2mm}(\ch)\tag{5}
\]

This is almost what we want.\\

Now, we claim \[
\Gamma, \tau_{\alpha}(x_{-q}, x_{i}),(\theta_{\alpha} \circ \tau_{\alpha})(x_{-q}, P_{k}) \Rightarrow \theta_\alpha (x_{i}, P_{k}) \tag{A}\]

Indeed, this is simply (taking $M_{i}$ unused) \[\Gamma, \tau_{\alpha}(x_{-q}, x_{i}), \exists M_{i}\rak{\tau_{\alpha}(x_{-q}, M_{i}) \wedge \theta_{\alpha}(M_{i}, P_{k})} \Rightarrow \theta_{\alpha}(x_{i}, P_{k})
\]

and we have this by the Extraction Lemma (with $\gamma = \theta_{\alpha}$) as $\tau_{\alpha}$ is functional. \qcl\\

Similarly we have\[
\Gamma, \tau_{\beta}(x_{-q}, x_{n+j}), (\theta_{\beta}{\circ} \tau_{\beta})(x_{-q}, P_{k}) \Rightarrow \theta_{\beta}(x_{n+j}, P_{k})\tag{B}
\]

\[(A),(B)\]
\[\overline{\Gamma, \eta, \theta_{\alpha}{\circ} \tau_{\alpha}(x_{-q}, P_{k}) \wedge \theta_{\beta} \circ \tau_\beta(x_{-q}, P_{k}) \Rightarrow \theta_{\alpha}(x_{i}, P_{k}) \wedge \theta_{\beta}(x_{n+j}, P_{k})} \qquad (\andd)
\]\[
\overline{\Gamma, \eta, \exists P_{k}\rak{(\theta_{\alpha} \circ \tau_{\alpha})(x_{-q}, P_{k}) \wedge(\theta_{\beta} \circ \tau_{\beta})(x_{-q}, P_{k})} \Rightarrow\exists P_{k} \rak{\theta_{\alpha}(x_{i}, P_{k}) \wedge \theta_{\beta}(x_{n+j}, P_{k})}}\hspace{1mm} (\eD) \hspace{1mm} (6)
\]

Chaining (5) and (6) we obtain (**). \qo\\

F2:\[
\Gamma, \neg E, \delta(x_{-q}) \Rightarrow \exists x_{i} \exists x_{n+j}\rak{\tau_{\alpha}(x_{-q}, x_{i}) \wedge \tau_{\beta}(x_{-q}, x_{n+j})} \tag{*}
\]

We have \begin{paracol}{2}\[
\Gamma, \neg E, \delta(x_{-q})\Rightarrow \exists x_{i} \tau_{\alpha}(x_{-q}, x_{i})\tag{1}\]
\switchcolumn
\[
\Gamma, \delta(x_{-q}) \Rightarrow \exists x_{n+j} \tau_{\beta}(x_{-q}, x_{n+j})\tag{2}\]
\end{paracol}

by \frelim applied to $\ftwo^{\tau_{\alpha}}, \ftwo^{\tau_{\beta}}$ respectively.

\[(1)(2)\]
\[\overline{\Gamma,\neg E, \delta(x_{-q}) \Rightarrow \exists x_{i} \rak{ \tau_{\alpha}(x_{-q}, x_{i})} \wedge \exists x_{n+j}\rak{ \tau_{\beta}(x_{-q}, x_{n+j})}} \qquad (\andrr)\tag{3}\]

The RHS of (3) is provably equivalent to the RHS of (*), by \ins/ \ex and the Equivalence Lemma. \qtw\\

F3: by \frelim it suffices to take $(P_{s})$, $(Q_{s})$ unused variables and derive \[
\Gamma, \neg E, \eta(x_{-q}, P_{i}, P_{n+j})\wedge \eta(x_{-q}, Q_{i}, Q_{n+j})\Rightarrow \bigwedge_{1 \leq s \leq n+m} P_s = Q_s \tag{*}
\]

But this is (technically after \un and rearrangement)
\[\Gamma, \tau_{\alpha}(x_{-q}, P_{i}), \tau_{\alpha}(x_{-q}, Q_{i}), \tau_{\beta}(x_{-q}, P_{n+j}), \tau_{\beta}(x_{-q}, Q_{n+j}) \Rightarrow \bigwedge_{s} P_{s}=Q_{s}\tag{*\pr}
\]
and we may write the RHS as $\bigwedge_i P_{i} = Q_{i} \wedge \bigwedge_{j} P_{n + j}=Q_{n+j}$. Thus (*\pr) is just obtained by
using \andd on the appropriate \frelim s of $\fth^{\tau_{\alpha}}, \fth^{\tau_{\beta}}$.\qth\\

This concludes the proof that $\eta = \tau_{\alpha}(x_{-q}, x_{i}) \wedge \tau_{\beta}(x_{-q}, x_{n+j})$ represents a morphism $\delta \rightarrow \FP$ (in the $\delta \neq \text{sentence}$ Subcase 0.0).

\section{Verification of solution}
We continue in this subcase and verify that $\eta$ represents the unique factoring map. To check both the solution of $\eta$ and its uniqueness as a morphism, we would like some characterization of those $\eta'\in \prehom(\delta, \FP)$ such that $\pi_{\alpha}\circ \eta' = \tau_{\alpha}$ and $\pi_{\beta}\circ \eta' = \tau_{\beta}$ simultaneously. Fortunately this is easy:\\

\begin{clm}\label{dec-c1c2} Let $\eta': \delta \to \FP$. Then $\pi_{\alpha}\circ \eta^{\prime}=\tau_{\alpha}$ and $\pi_{\beta} \circ \eta^{\prime}=\tau_{\beta}$ simultaneously hold iff we can derive

\[\Gamma, \neg E \Rightarrow\exists x_{i} \rak{\eta^{\prime}(x_{-q}, x_i, x_{n+j})} \leftrightarrow \tau_{\beta}(x_{-q}, x_{n+j}) \tag{C1}\]

\[\Gamma, \neg E \Rightarrow \exists x_{n+j}\rak{\eta^{\prime}(x_{-q}, x_{i}, x_{n+j})} \leftrightarrow \tau_{\alpha}(x_{-q}, x_{i}). \tag{C2}\]
\end{clm}

\pf By \cref{dec-sincoming} of Special Morphisms, we have\[
\pi_{\beta}\circ \eta\pr = \exists x_i \rak{\eta^{\prime}(x_{-q}, x_i, x_{n+j})}\tag{1}
\]\[
\pi_{\alpha}\circ \eta^{\prime} = \exists x_{n+j}\rak{\eta^{\prime}(x_{-q}, x_{i}, x_{n+j})}\tag{2}
\]
since, in the case of (1), the map $\tf_{\beta}: \free(\beta)\to \free(\FP)$ misses exactly the $T_i$; in the case of (2) $\tf_{\alpha}:\free(\alpha)\to \free(\FP)$ misses the $T_{n + j}$.

The definition of (1) representing the same morphism $\delta \to \FP$ as $\tau_{\beta}$, is that (C1) should be derivable; similarly for (2) and (C2). \qed\\

Now we can easily do the verification. Indeed, to check $\tau_{\alpha} (x_{-q}, x_{i})\wedge \tau_{\beta}(x_{-q}, x_{n+j})$ is a solution, it suffices to obtain (C1) and (C2). The forwards half of (C1) in this case is\[
\Gamma, \neg E, \exists x_i \rak{\tau_\alpha(x_{-q}, x_i)\wedge \tau_{\beta}(x_{-q}, x_{n+j})} \Rightarrow \tau_{\beta}(x_{-q}, x_{n+j})
\] the derivation of which is trivial, without use of any background theory. Conversely, by \frelim on $\fone^{\tau_{\beta}}$ we have
\[
\Gamma, \neg E, \tau_\beta(x_{-q}, x_{n+j})\Rightarrow \delta(x_{-q})\tag{2}
\] and by \frelim on $\ftwo^{\tau_{\alpha}}$ we have 
\[
\Gamma, \neg E, \delta (x_{-q})\Rightarrow \exists x_i\tau_\alpha(x_{-q}, x_i)\tag{3}
\]
\[
(2)(3)
\]\[
\overline{\Gamma, \neg E, \tau_\beta(x_{-q}, x_{n+j})\Rightarrow\exists x_i\tau_\alpha(x_{-q}, x_i)}\qquad (\ch)
\]\[
\overline{\Gamma,\neg E, \tau_\beta(x_{-q}, x_{n+j})\Rightarrow\exists x_i\rak{\tau_\alpha(x_{-q}, x_i)}\wedge \tau_\beta(x_{-q}, x_{n+j})}
\]\[
\qquad \qquad \overline{\Gamma, \tau_\beta(x_{-q}, x_{n+j})\Rightarrow\exists x_i \rak{\tau_\alpha(x_{-q}, x_i)\wedge \tau_\beta(x_{-q}, x_{n+j})}}\qquad (\ins)
\]
This is the other half of (C1). The derivation of (C2) is exactly analogous.\qed

\subsection{Uniqueness}\label{seg-uniqueness}
Now let $\eta' : \delta \rightarrow \FP$ be another solution, i.e. (by virtue of \cref{dec-c1c2}) another formula such that we have (C1) and (C2). Suppose we want to show $T$ proves $\E^{\eta\to \eta\pr}$. Let $(M_{s})$ be new. By alpha-equivalence it is easy to see that we derive 
\[
\Gamma, \neg E \Rightarrow \exists M_i\eta (x_{-q}, M_i, x_{n+j}) \rightarrow \exists M_i \eta^{\prime}(x_{-q}, M_i, x_{n+j}) \tag{1}
\]
\[
\Gamma, \neg E \Rightarrow \exists M_{n+j}\eta (x_{-q}, x_{i}, M_{n+j}) \rightarrow \exists M_{n+j} \eta^{\prime}(x_{-q}, x_{i}, M_{n+j}) \tag{2}
\]

since we derive that both existential conditions in (1) are provably equivalent to the condition $\tau_\beta(x_{-q}, x_{n+j})$, and both in (2) to $\tau_\alpha(x_{-q}, x_{i})$.\\

Certainly also
\[
\Gamma, \lnot E, \eta(x_{-q}, x_i, x_{n+j})\Rightarrow \exists M_i \eta (x_{-q}, M_i, x_{n+j}) \tag{a}
\]
\[
\Gamma, \neg E, \eta(x_{-q}, x_{i}, x_{n+j}) \Rightarrow \exists M_{n+j} \eta(x_{-q}, x_{i}, M_{n+j})\tag{b}
\] by \eS, and we obtain \[
\Gamma, \neg E, \eta(x_{-q}, x_{i}, x_{n+j})\Rightarrow \exists M_i\eta^{\prime}(x_{-q},M_i, x_{n+j})\tag{a\pr}
\]
\[
\Gamma, \neg E, \eta(x_{-q}, x_{i}, x_{n+j})\Rightarrow\exists M_{n+j}\eta^{\prime}(x_{-q},x_{i}, M_{n+j})\tag{b\pr}
\]

by chaining with (1) and (2) respectively. Now we have
\[
\Gamma, \neg E, \eta^{\prime}(x_{-q}, M_i, x_{n+j}), \eta^{\prime}(x_{-q}, x_{i}, M_{n+j})\Rightarrow\bigwedge_{i} M_i=x_{i}\wedge\bigwedge_{j} x_{n+j}=M_{n+j}\tag{3}
\]

by \frelim on $\fth^{\eta'}$. By \assm and \subs we have
\[
\eta'(x_{-q}, M_i, x_{n + j})\wedge \bigwedge_{i} M_i = x_{i} \Rightarrow \eta'(x_{-q}, x_{i}, x_{n + j})\tag{c}
\]
\[
\eta'(x_{-q}, M_i, x_{n + j})\wedge \bigwedge_{j} M_{n+j} = x_{n+j} \Rightarrow \eta'(x_{-q}, x_{i}, x_{n + j})\tag{d}
\]

\[(c) (d)
\]
\[
\overline{\eta'(x_{-q}, M_i, x_{n + j})\wedge \eta'(x_{-q}, x_{i}, M_{n + j}), \bigwedge_{i} M_i = x_{i}\wedge \bigwedge_{j}M_{n + j} =  x_{n+j}\Rightarrow \eta'(x_{-q}, x_{i}, x_{n+j})} (\andd) \tag{4}
\]

and by chaining (3) and (4) we have 
\[
\Gamma, \neg E, \eta'(x_{-q}, M_i, x_{n + j}),\eta'(x_{-q}, x_{i}, M_{n + j}) \Rightarrow \eta'(x_{-q}, x_{i}, x_{n+j}).
\]By \eA twice this gives \[
\Gamma, \exists M_i \eta'(x_{-q}, M_i, x_{n + j}),\exists M_{n+j}\eta'(x_{-q}, x_{i}, M_{n + j}) \Rightarrow \eta'(x_{-q}, x_{i}, x_{n+j}).
\] Chaining with (a\pr) and (b\pr) we have \[
\Gamma, \neg E,  \eta(x_{-q}, x_{i}, x_{n+j}) \Rightarrow \eta'(x_{-q}, x_{i}, x_{n+j})
\]
which shows $\E^{\eta \to \eta'}$. $\E^{\eta' \to \eta}$ is derived exactly symmetrically; this argument only invoked functionality sequents for $\eta$ and $\eta'$, and the sequents (1) and (2), for which we have corresponding reverses.

Therefore we have shown $\eta = \eta'$ as morphisms $\delta \to \FP$. This concludes Subcase 0.0. \qed\\

\section{Bootstrap}
We will now argue abstractly that we have in fact solved all of Case 0, i.e. constructed a genuine fiber product in $\syn(T)$ of the two type 1 morphisms $\theta_{\alpha}$ and $\theta_{\beta}$. But to do so we will need the final isomorphism criterion, for the dual type 2/type 3 maps. We develop this here.\\

\subsection{Isomorphisms of sentences with nonsentences}
\begin{prop}\label{dec-t3iso} Let $s \overset{\theta(x_{i})}{\longrightarrow} \varphi(X_{i})$ be type 3. Then it is an isomorphism if and only if $T$ derives the sentence \[
\forall x_{i}\rak{\theta(x_{i})\leftrightarrow \varphi(x_{i})};
\] in other words iff $\theta$ and $\varphi\sbb{x_{i}}{X_{i}}$ are provably equivalent. \end{prop}

\pf Suppose it is an isomorphism; we then need to obtain\[
\Gamma, \neg E \Rightarrow \theta(x_{-i})\leftrightarrow \varphi(x_{-i}). \tag{*}
\]
$\Gamma, \neg E \Rightarrow \theta(x_{-i})\rightarrow \varphi(x_{-i})$ is immediate from $\mfone^{\theta}$. As $s \to \varphi$ is an isomorphism there must be a (formal) map in the other direction, 
and the composition $\varphi \rightarrow s \overset{\theta}{\rightarrow} \varphi$ is $I_{\varphi}$, which referencing the definitions means we have\[
\varphi(x_{-i})\wedge \theta(x_{i})= \varphi(x_{-i})\wedge \bigwedge_{i} x_{i} = x_{-i} \qquad \qquad \text{as morphisms $\varphi \to \varphi$.}
\]

By \frelim on $\E^{\text{RHS}\to \text{LHS}}$, letting $(M_{i})$ be unused variables, we have\[
\Gamma, \neg E, \varphi(x_{-i}), \bigwedge_{i}M_{i} = x_{-i} \Rightarrow\theta(M_{i}),\tag{1}
\]\[
\overline{\Gamma, \neg E, \varphi(x_{-i})\wedge\bigwedge_{i} M_{i} = x_{-i}\Rightarrow\theta(M_{i})\wedge\bigwedge_{i} M_{i} = x_{-i}}.\tag{1\pr}
\]
By \subs we have \[
\theta(M_{i})\wedge\bigwedge_{i} M_{i} = x_{-i} \Rightarrow \theta(x_{-i}) \tag{2}
\]Chaining (1\pr) and (2) we have\[
\Gamma, \neg E, \varphi(x_{-i}), \bigwedge_{i}M_{i} = x_{-i} \Rightarrow \theta(x_{-i})\tag{3}
\]\[
\qquad \overline{\Gamma, \neg E, \varphi(x_{-i}), \exists M_{i}\rak{\bigwedge_{i}M_{i} = x_{-i}} \Rightarrow \theta(x_{-i})}\qquad (\eA)\tag{3\pr}
\]
The last hypothesis on the LHS in (3\pr) is vacuous given $\neg E$, so we can drop it to obtain $\Gamma, \neg E, \varphi(x_{-i}) \Rightarrow \theta(x_{-i})$. This combines with the easy direction to give (*).\qed\\

Conversely, suppose $\theta$ and $\varphi(x_{i})$ are provably equivalent. By $\mfone^{\theta}$, we have\[
\Gamma, \neg E, \theta(x_{i})\Rightarrow s\]
\[\overline{\Gamma, \neg E, \varphi(x_{i})\Rightarrow s}
\] Therefore we have a type 2 morphism $\varphi \to s$. The composite $s \overset{\theta}{\rightarrow}\varphi\rightarrow s$ is forced to be the unique formal morphism $s\to s$, which is $I_{s}$. For the other order, we need to show that\[
\varphi(x_{-i})\wedge \theta(x_{i})= I_{\varphi}=\varphi(x_{-i})\wedge \bigwedge_{i} x_{i} = x_{-i} \qquad \qquad \text{as morphisms $\varphi \to \varphi$.}
\]\\

So we need
\[
\Gamma, \neg E, \varphi(x_{-i})\wedge \theta(x_{i})\Rightarrow \varphi(x_{-i})\wedge \bigwedge_{i} x_{i} = x_{-i} \tag{*a}
\]
By provable equivalence, it suffices to change both $\varphi$'s to $\theta$'s:\[
\Gamma, \neg E, \theta(x_{-i})\wedge \theta(x_{i})\Rightarrow \theta(x_{-i})\wedge \bigwedge_{i} x_{i} = x_{-i} \tag{*a\pr}
\]
This is now an immediate consequence of $\mfth^{\theta}$.\\

We also need
\[
\Gamma, \neg E,  \varphi(x_{-i})\wedge \bigwedge_{i} x_{i} = x_{-i}\Rightarrow \varphi(x_{-i})\wedge \theta(x_{i}) \tag{*b}
\]

This is immediate by provable equivalence and \subs. This concludes the proof of \cref{dec-t3iso}.\qed\\

\begin{cor}\label{dec-isosingular}Let $s \overset{\theta}{\longrightarrow} \varphi(X_{i})$ be type 3. Then for any sequence $(T_{i})$ of distinct variables, we have an isomorphism\[
s \overset{\theta}{\longrightarrow} \theta(T_{i}).
\]\end{cor}

\pf We must first establish $\theta: s \to \theta(T_{i})$ is a type 3 morphism. It has the right free variables. Derivability of $\mfone^{\theta}$ for this claim is trivial, since we have $\Gamma, \neg E, \theta \Rightarrow s$ by $\mfone^{\theta}$ for $\theta: s \to \varphi(X_{i})$, and the other required conclusion is tautological. $\mftwo^{\theta}$ does not mention the target object and hence is identical with regard to the claims $\theta: s\to \varphi(X_{i})$ and $\theta: s \to \theta(T_{i})$. Similarly for $\mfth$ as it only mentions the formula defining the map.\\

That this morphism is then an iso is immediate from \cref{dec-t3iso}.\qed\\

\begin{cor}\label{dec-t2iso} Let $\varphi(X_{i})\rightarrow s$ be type 2. Then it is an isomorphism if and only if $\varphi(x_{i})\in \prehom(s, \varphi(X_{i}))$, in which case $\varphi(x_{i})$ represents the (type 3) inverse.\end{cor}

\pf If it is an isomorphism, it must have an inverse type 3 iso $s \to \varphi(X_{i})$, and by Prop 1 $\varphi(x_{i})$ is a representing premorphism for this. Conversely, if $\varphi: s\to \varphi(X_{i})$ is a premorphism, by Prop 1 it represents an isomorphism. There is only one possibility for the inverse iso, namely the given formal type 2 morphism $\varphi(X_{i}) \to s$. \qed\\

\subsection{Subcase 0.1}
Suppose now $\delta = s$, and we are given type 3 maps $s \overset{\tau	_{\alpha}(x_{i})}{\rightarrow} \alpha(X_{i})$, $s\overset{\tau_{\beta}(x_{j})}{\rightarrow}\beta(Y_{j})$ such that $\theta_{\alpha}\circ \tau_{\alpha} = \theta_{\beta}\circ \tau_{\beta}$. We shall reduce to the Subcase 0.0 situation.

Indeed, by \cref{dec-isosingular}, $s$ is isomorphic to a nonsentence $\delta'$ by an isomorphism $\iota: s\to \delta'$ (we can take $\delta'$ and $\iota$ to both be $\tau_{\alpha}(x_{i})$, say.)

Then certainly $\theta_{\alpha} \circ \tau_{\alpha}\circ \iota^{-1} = \theta_{\beta} \circ \tau_{\beta}\circ \iota^{-1}$. Therefore by Subcase 0.0, there is a unique $\widetilde{\eta}:\delta' \to \FP$ such that we have
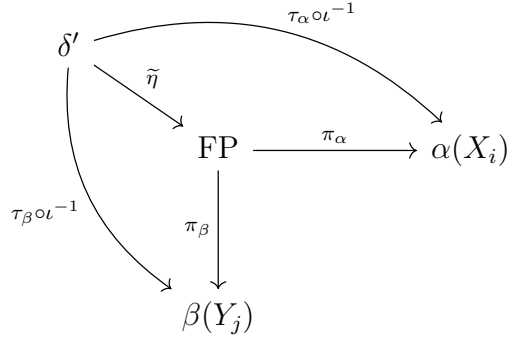
\begin{center}\begin{tikzcd}
\delta' \arrow[rd, "\widetilde{\eta}"] \arrow[rrrd, "\tau_{\alpha} \circ \iota^{-1}", bend left] \arrow[rddd, "\tau_{\beta}\circ \iota^{-1}"', bend right] &                                                       &  &             \\
                                                                                                                                                    & \FP \arrow[dd, "\pi_\beta"'] \arrow[rr, "\pi_\alpha"] &  & \alpha(X_i) \\
                                                                                                                                                    &                                                       &  &             \\
                                                                                                                                                    & \beta(Y_j)                                            &  &            
\end{tikzcd}\end{center}

It follows by pure abstract nonsense that the unique $\eta: s \to \FP$ such that $\pi_{\alpha} \circ \eta = \tau_{\alpha}$ and $\pi_{\beta} \circ \eta = \tau_{\beta}$, is $\eta:=\widetilde{\eta} \circ \iota$. Thus we have solved the factorization problem in the subcase that $\delta = s$ is a sentence. This completes the full Case 0, i.e. the verification of the universal property of the pullback for $(\FP, \pi_{\alpha}, \pi_{\beta})$.\qed

\section{Other cases}
Out of the five other possible combinations of maps to fiber, three (where the common target is a sentence) will essentially consist of producing an ordinary product, and will be more or less trivial. The other two will be reduced to Case 0 easily. We begin with the main product case.

\subsection{Case 1}\label{seg-oversentence}
Suppose we wish to construct a pullback of\begin{center}\begin{tikzcd}
\alpha(X_{i})\arrow[rd] & & \beta(Y_{j})\arrow[ld]\\
& s&
\end{tikzcd}\end{center}

We let again $(T_{s})_{1 \leq s \leq m+n}$ be any variables in increasing order, and consider \[
\PD:= \alpha(T_{i})\wedge \beta(T_{n +j}).
\]Similarly to in Case 0, we have obvious special morphisms \[\pi_{\alpha}: \PD \to \alpha,
\]\[\pi_{\beta}:\PD \to \beta.
\]
We shall show a pullback diagram is 
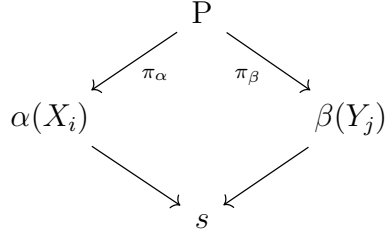
\begin{center}\begin{tikzcd}
& \PD \arrow[ld, "\pi_{\alpha}"] \arrow[rd, "\pi_{\beta}"']& \\
\alpha(X_{i})\arrow[rd] & & \beta(Y_{j})\arrow[ld]\\
& s&
\end{tikzcd}\end{center}

\pf We have commutativity of the diamond trivially. And in general, because any composition with sentence target is formal, given an arbitrary $\delta$ and maps $m_{\alpha}:\delta \to \alpha$ and $m_{\beta}:\delta \to \beta$, we have a similar commutative diamond. Hence, to check the universal property in this case, we must exactly check that 
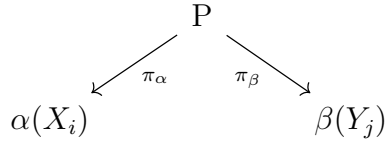
\begin{center}\begin{tikzcd}
& \PD \arrow[ld, "\pi_{\alpha}"] \arrow[rd, "\pi_{\beta}"']& \\
\alpha(X_{i})& & \beta(Y_{j})
\end{tikzcd}\end{center} is a {\it product} diagram.\\

\subsubsection{Subcase 1.0}
As in Case 0, we first solve the factoring problem for 
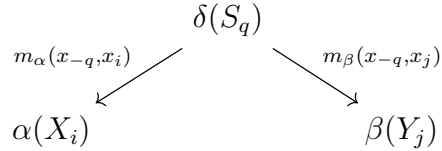
\begin{center}\begin{tikzcd}
& \delta(S_q) \arrow[ld, "{m_{\alpha}(x_{-q}, x_{i})}"'] \arrow[rd, "{m_{\beta}(x_{-q}, x_{j})}"]& \\
\alpha(X_{i})& & \beta(Y_{j})
\end{tikzcd}\end{center} 
and then use this to bootstrap the subcase in which $\delta$ is a sentence.\\

Indeed, given the above diagram we claim\[
m_{\alpha}(x_{-q}, x_{i})\wedge m_{\beta}(x_{-q}, x_{n + j})
\]
represents the unique $\eta: \delta \to \PD$ such that $\pi_{\alpha}\circ \eta = m_{\alpha}$ and $\pi_{\beta}\circ \eta = m_{\beta}$. First we must check $\eta: \delta\to \alpha(T_{i})\wedge \beta(T_{n + j})$.\\

F1: \[
\Gamma, \neg E, m_{\alpha}(x_{-q}, x_{i})\wedge m_{\beta}(x_{-q}, x_{n + j})\Rightarrow \delta(x_{-q}) \wedge \alpha(x_{i})\wedge \beta(x_{n+j}) \tag{*}
\]
This is immediate from $\fone^{m_{\alpha}}$ and $\fone^{m_{\beta}}$. \qed\\

F2: \[
\Gamma, \neg E, \delta(x_{-q}) \Rightarrow \exists x_{i} \exists x_{n + j}\rak{m_{\alpha}(x_{-q}, x_{i})\wedge m_{\beta}(x_{-q}, x_{n + j})} \tag{*}
\]
By \ins, \ex and Equivalence Lemma, the RHS is provably equivalent to $\exists x_{i}\rak{m_{\alpha}(x_{-q}, x_{i})}\wedge \exists x_{n+j}\rak{m_{\beta}(x_{-q}, x_{n+j})}$, so this is immediate from $\ftwo^{m_{\alpha}}$ and $\ftwo^{m_{\beta}}$. \qed\\

F3: \[
\Gamma, \neg E, m_{\alpha}(x_{-q}, x_{i}), m_{\beta}(x_{-q}, x_{n + j}), m_{\alpha}(x_{-q}, T_{i}), m_{\beta}(x_{-q}, T_{n + j})\Rightarrow \bigwedge_{s} x_{s} = T_{s} \tag{*}
\]
Immediate from $\fth^{m_{\alpha}}$ and $\fth^{m_{\beta}}$.\qed\\

Now, as in Case 0, for any $\eta': \delta \to \PD $ we have by the precomposition formulas \[
\pi_{\alpha} \circ \eta' = \exists M_{n + j} \rak{\eta'(x_{-q}, x_{i}, M_{n + j})}\tag{A}
\]\[
\pi_{\beta} \circ \eta' = \exists M_{i} \rak{\eta'(x_{-q}, M_{i}, x_{n + j})}. \tag{B}
\]
Therefore to check $\eta$ is a solution, it suffices to obtain\[
\Gamma, \neg E\Rightarrow  \exists M_{n + j} \rak{m_{\alpha}(x_{-q}, x_{i})\wedge m_{\beta}(x_{-q}, M_{n + j})
} \leftrightarrow m_{\alpha}(x_{-q}, x_{i}) \tag{*a}
\]\[
\Gamma, \neg E\Rightarrow  \exists M_{i} \rak{m_{\alpha}(x_{-q}, M_{i})\wedge m_{\beta}(x_{-q}, x_{n + j})
} \leftrightarrow m_{\beta}(x_{-q}, x_{n+j}) \tag{*b}
\]
The derivations are obvious, making use only of \ins and \ex and $\fone^{m_{\alpha}}$, $\ftwo^{m_{\beta}}$ (respectively $\fone^{m_{\beta}}$ and $\ftwo^{m_{\alpha}}$ for (*b)).\qed\\

Lastly for this main subcase, we must check uniqueness, which will follow if we show that any $\eta_{1}$ and $\eta_{2}$ satisfying both (A) and (B) are provably equivalent. The argument for this is a nearly verbatim copy of \cref{seg-uniqueness}.\qed

\subsubsection{Subcase 1.1 (bootstrap)}
Now suppose we are given \begin{center}
\begin{tikzcd}
& s\arrow[rd, "{m_{\alpha}(x_{i})}"]\arrow[ld, "{m_{\beta}(x_{j})}"'] &\\
\alpha(X_{i}) & & \beta(Y_{j})
\end{tikzcd}
\end{center}

As in the bootstrap for Case 0, to show that the factoring problem has unique solution, it suffices to obtain any isomorphism of $s$ with a nonsentence. We have this again by \cref{dec-isosingular} and either of the type 3 maps $m_{\alpha}$, $m_{\beta}$. \qed\\

This concludes the verification of universal property of the product, and hence the proof that our claimed diamond is a pullback. \qed

\subsection{Case 2}\label{seg-sentoversent}

Suppose we wish to construct a pullback for \begin{center}
\begin{tikzcd}
\alpha(X_{i}) \arrow[rd]& &s' \arrow[ld]\\
& s &
\end{tikzcd}
\end{center}

We claim it is given by \begin{center}
\begin{tikzcd}
& \alpha(X_{i})\wedge s' \arrow[ld]\arrow[rd]&\\
\alpha(X_{i}) \arrow[rd]& &s' \arrow[ld]\\
& s &
\end{tikzcd}
\end{center}
where the upper-left map is a special mono.\\

\pf As in Case 1, the claim is just the assertion that \begin{center}
\begin{tikzcd}
& \alpha(X_{i})\wedge s' \arrow[ld]\arrow[rd]&\\
\alpha(X_{i}) & &s' \\
\end{tikzcd}
\end{center}
is a product.
Therefore consider $\theta(x_{-q}, x_{i}): \delta(S_q) \to \alpha(X_{i})$, where $\delta$ also has a type 2 map to $s'$. We claim $\theta: \delta\to \alpha\wedge s'$ also, and that this uniquely solves the factoring problem.

Indeed, we obtain F1 easily from the given $\fone^{\theta}$ and the fact that we have\[
\Gamma, \neg E, \delta(x_{-q})\Rightarrow s',
\] by existence of the type 2 $\delta \to s'$. F2 and F3 do not mention the target formula and hence are immediate. The composition relations hold trivially. Finally, we have uniqueness as $\alpha \wedge s'\to \alpha$ is a monomorphism. \qed

If $\delta$ is a sentence, the given map to $\alpha(X_{i})$ again gives rise to an isomorphism with a nonsentence via \cref{dec-isosingular}, and we have the bootstrap reduction. \qed

\subsection{Case 3}
Suppose we wish to construct a pullback for 
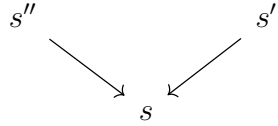
\begin{center} \begin{tikzcd}
s'' \arrow[rd]& &s' \arrow[ld]\\
& s &
\end{tikzcd}
\end{center}

It is just $s' \wedge s''$. Any object with (necessarily formal) maps to $s'$ and $s''$ must have a formal map to this, and the composition relations and uniqueness are automatic. \qed

\subsection{Case 4}
Suppose we wish to {\it show existence of} a pullback for 
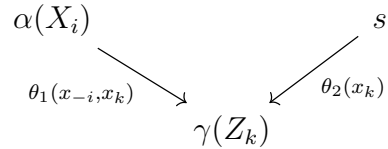
\begin{center} \begin{tikzcd}
\alpha(X_{i}) \arrow[rd, "{\theta_{1}(x_{-i}, x_{k})}"']& &s \arrow[ld, "\theta_{2}(x_{k})"]\\
& \gamma(Z_{k}) &
\end{tikzcd}
\end{center}
(we will not need this case in practice, so we don't analyze it in detail). By \cref{dec-isosingular}, we have an isomorphism $\widetilde{\theta_{2}}: s \to \theta_{2}(Z_{k})$. We invoke Case 0 to obtain a pullback of 
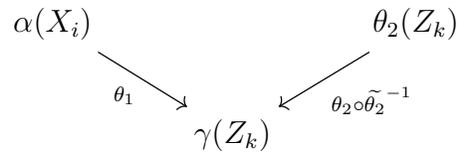
\begin{center} \begin{tikzcd}
\alpha(X_{i}) \arrow[rd, "{\theta_{1}}"']& &\theta_{2}(Z_{k}) \arrow[ld, "\theta_{2} \circ \widetilde{\theta_{2}}^{-1}"]\\
& \gamma(Z_{k}) &
\end{tikzcd}
\end{center}
and by abstract nonsense we get a pullback of the original maps by postcomposing the right projector with $\widetilde{\theta_{2}}^{-1}$.\qed\\

\subsection{Case 5}
Suppose we wish to show existence of a pullback for 
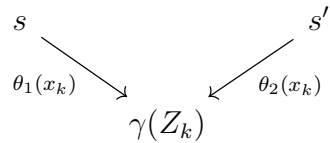
\begin{center} \begin{tikzcd}
s \arrow[rd, "{\theta_{1}(x_{k})}"']& &s' \arrow[ld, "\theta_{2}(x_{k})"]\\
& \gamma(Z_{k}) &
\end{tikzcd}
\end{center}
We argue exactly as above but using isomorphisms for both legs. \qed\\

\chapter{Factorization of Morphisms}\label{ch-factorization}

In a category with fiber products, an {\it effective epimorphism} is a morphism $g$ which is a coequalizer of the projection maps $\pi_{1}$, $\pi_{2}: g\times g \rightrightarrows \text{dom(}g)$. This is a stronger condition than being an epi. It is a fact which generalizes the surjective/injective factorization in \sett, that if $f$ can be factored as $h \circ g$ where $h$ is a monomorphism and $g$ is an effective epimorphism, then that factorization is unique up to unique isomorphism (\cite{lu}, Lecture 3).

Here, we construct such a factorization for arbitrary $f$ in $\syn(T)$. Its uniqueness will give us certain ``standardized forms'' which will be of use in the next Chapter. We do the factorization first in Case 0 ($f$ of type 1), and then for the other cases, which will be easy.\\

We fix the new ``index-limit pair'' $(w, c)$.

\section{Case 0}
\subsection{Preliminaries and definitions}

Let $\theta: \alpha(X_{i}) \rightarrow \gamma(Z_{k})$ be type 1, and let $\im \theta$ be the object \[\exists M_{i} \theta (M_{i}, Z_{k}) \]
where $(M_{i})$ are distinct unused variables.\\

\begin{clm}\label{dec-imagereceives}We have $\theta:\alpha \rightarrow \im \theta$.\end{clm}

\pf F1: we must derive \[
\Gamma, \neg E, \theta(x_{-i}, x_{k}) \Rightarrow \alpha (x_{-i}) \wedge \exists M_{i} \theta(M_{i}, x_{k})\tag{*}\]

By F1 for $\theta:\alpha\to \gamma$, we have 
\[\Gamma, \neg E, \theta(x_{-i}, x_{k}) \Rightarrow\alpha (x_{-i})\tag{1}
\]
and trivially \[
\neg E, \theta(x_{-i}, x_{k}) \Rightarrow \exists M_{i} \theta(M_{i}, x_{k})\tag{2}\]
whence we get (*) by \andd.\qo

F2, F3: these sentences are independent of the target, so we have them as $\theta:\alpha \to \gamma$. \qcl\\

\begin{clm} \label{dec-toimageisee}$\theta: \alpha\to \im \theta$ is an effective epimorphism.\end{clm}

\pf As we have constructed it, the fiber product of $\theta$ with itself is an object $\FP(T_{i}, T_{n+i})$ in ascending variables $T_{1}, \ldots T_{2n}$, given by 

\[\FP=\exists W_{k}\rak{\theta(T_{i}, W_{k}) \wedge \theta(T_{n+i}, W_{k})}.\]

The projections $\FP \rightrightarrows \alpha$ are given by
\[\pi_{1}(x_{-i}, x_{-n-i}, x_{i})=\FP(x_{-i}, x_{-n-i}) \wedge \bigwedge_{i} x_{i} = x_{-i}\]
\[\pi_{2}(x_{-i}, x_{-n-i}, x_{i})=\FP(x_{-i},x_{-n-i}) \wedge \bigwedge_{i} x_{i}=x_{-n-i}.\]

We must show that $\theta: \alpha \rightarrow \im \theta$ is a coequalizer of $\pi_{1}$ and $\pi_{2}$. By the fiber product construction, we know $\theta\circ \pi_{1} = \theta\circ \pi_{2}$, so this makes sense.\\

To check the universal property, we consider an arbitrary $u: \alpha \to \chi$ such that $u \circ \pi_{1} = u \circ \pi_{2}$.\\

First suppose $\chi = s$ is a sentence. Then uniqueness of a factoring $\eta: \im \theta \to \chi$ is automatic given existence. We have existence by the following: as $u : \alpha \to s$, we have\[
\Gamma, \neg E, \alpha(M_{i}) \Rightarrow s
\]\[\overline{\Gamma, \neg E, \exists M_{i} \alpha(M_{i}) \Rightarrow s}.\tag{1}
\]By $\fone^{\theta}$,\[
\Gamma, \neg E, \theta(M_{i}, Z_{k})\Rightarrow \alpha(M_{i})
\]\[
\overline{\Gamma, \neg E, \exists M_{i}\theta(M_{i}, Z_{k})\Rightarrow \exists M_{i}\alpha(M_{i})} \qquad (\eD)\tag{2}
\]
Chaining (2) and (1) gives \[
\Gamma, \neg E, \im \theta \Rightarrow s
\]
so we have a formal $\eta: \im \theta \to s$. We must have $\eta\circ \theta = u$ as there is nothing else for the composition to be.\qed\\

Now suppose $\chi = \chi(Q_{w})$. To verify the existence and uniqueness in this main subcase, we need a characterization of the maps $u$ to be factored.\\

\begin{lem}\label{dec-mapslemma} Let $u : \alpha \to \chi(Q_{w})$ be such that $u \circ \pi_{1}= u \circ \pi_{2}$. 
Then we can derive\[
\Gamma \Rightarrow \forall x_{-i} \forall x_{-n - i} \rak{\FP(x_{-i}, x_{-n - i}) \rightarrow \exists x_{w}\rak{u(x_{-i}, x_{w})\wedge u(x_{-n-i}, x_{w})}} \tag{*}
\]
Set-theoretically this says that such a $u$ respects the fibers of $\theta$.\\
\end{lem}

\pf By \cref{dec-soutgoing} of Special Morphisms, we have $u \circ \pi_{1} = \FP(x_{-i}, x_{-n-i})\wedge u(x_{-i}, x_w)$, and $u \circ \pi_{2} = \FP(x_{-i}, x_{-n-i})\wedge u(x_{-n-i}, x_w)$. Thus if these are equal as morphisms $\FP \to \chi$, we have \[\Gamma, \neg E,  \FP(x_{-i}, x_{-n-i})\wedge u(x_{-i}, x_w)\Rightarrow \FP(x_{-i}, x_{-n-i})\wedge u(x_{-n-i}, x_w)\tag{1}
\]\[
\overline{\Gamma, \neg E,  \FP(x_{-i}, x_{-n-i})\wedge u(x_{-i}, x_w) \Rightarrow  u(x_{-i}, N_w) \wedge u(x_{-n-i}, x_{w})} \tag{1\pr}
\]
\[
\overline{\Gamma, \neg E, \exists x_{w}(u\circ \pi_{1})(x_{-i}, x_{-n-i}, x_{w}) \Rightarrow \exists x_{w}\rak{u(x_{-i}, x_{w}) \wedge u(x_{-n-i}, x_{w})}} \qquad (\eD)\tag{1\pr\pr}
\]

Now by $\ftwo^{u \circ \pi_{1}}$, we have\[
\Gamma, \neg E, \FP(x_{-i}, x_{-n-i}) \Rightarrow \exists x_{w}(u\circ \pi_{1})(x_{-i}, x_{-n-i}, x_{w}). \tag{2}
\]

\[(2), (1)
\]\[
\overline{\Gamma, \neg E, \FP(x_{-i}, x_{-n-i}) \Rightarrow \exists x_{w}\rak{u(x_{-i}, x_{w})\wedge u(x_{-n-i}, x_{w})}} \qquad (\ch)
\]

which gives (*) by \meta, \frlS, and eliminating $\neg E$.\qed\\

Now we can solve the factoring problem in this subcase. Let $u: \alpha \to \chi(Q_{w})$ be a map such that $u \circ \pi_{1} = u \circ \pi_{2}$. We claim the unique map $\eta: \im \theta \to \chi$ such that $ u = \eta \circ \theta$ is defined by\[
\eta = \exists R_{i}\rak{\theta(R_{i}, x_{-k}) \wedge u(R_{i}, x_{w})}
\] where $(R_{i})$ are unused variables.

\subsection{Functionality}
Verification of the first two functionality sequents for the claim $ \eta: \im \theta\rightarrow \chi$ is straightforward. Obtaining $\fth^\eta$, or checking that the map is well-defined, requires \cref{dec-mapslemma}.\\

F3: extend the sequence of $R$'s by letting $(R_{n+i})$ be additional unused variables. It suffices to get\[
\Gamma , \neg E, \exists R_{i}\rak{\theta(R_{i}, x_{-k}) \wedge u(R_{i}, x_{w})}, \exists R_{n+i}\rak{\theta(R_{n+i}, x_{-k}) \wedge u(R_{n+i}, x_{c+w})}\Rightarrow \bigwedge_{w} x_{w}=x_{c+w}\tag{*}\]
Trivially we have
\[
\Gamma,\neg E,\theta(R_{i},x_{-k}),\theta(R_{n+i},x_{-k})\Rightarrow\exists W_{k}\rak{\theta(R_{i},W_{k})\wedge\theta(R_{n+i},W_{k})}\tag{1}
\] and the right-hand side is $\FP(R_{i}, R_{n+i})$. 
Let $(N_{w})$ be unused variables. By \frelim on the output of \cref{dec-mapslemma} (taking the bound $x_{-i}$ there to $R_{i}$, $x_{-n-i}$ to $R_{n+i}$) we have
\[
\Gamma,\neg E,\theta(R_{i},x_{-k}),\theta(R_{n+i},x_{-k})\Rightarrow\exists N_{w}\rak{u(R_{i},N_{w})\wedge u(R_{n+i},N_{w})}\tag{2}.
\] 

Now by two \frelim's of $\fth^{u}$ \andd ed together, we have
\[
\Gamma,\neg E, u(R_{i}, N_{w}), u(R_{n+i}, N_{w}),u(R_{i}, x_{w}), u(R_{n+i}, x_{c+w}) \Rightarrow \bigwedge_{w} N_{w}=x_{w} \wedge \bigwedge_{w} N_{w} =x_{c+w}\tag{3}
\] and by transitivity of equality (formally, a manipulation with \subs) this yields
\[
\Gamma, \neg E, u(R_{i},N_{w}), u(R_{n+i},N_{w}), u(R_{i},x_{w}), u(R_{n+i},x_{c+w})\Rightarrow\bigwedge_{w} x_{w}=x_{c+w}. \tag{4}
\] By \eA (4) gives
\[
\Gamma, \neg E, \exists N_{w}\rak{u(R_{i}, N_{w}) \wedge u(R_{n+i}, N_{w})},u(R_{i}, x_{w}), u(R_{n+i}, x_{c+w})\Rightarrow\bigwedge_{w} x_{w}=x_{c+w}.\tag{4\pr}
\]

Chaining (2) and (4\pr) we obtain\[
\neg E,\Gamma,\theta(R_{i},x_{-k}),\theta(R_{n+i},x_{-k}),u(R_{i},x_{w}),u(R_{n+i},x_{c+w})\Rightarrow\bigwedge_{w} x_{w}=x_{c+w}.
\]
Rearranging and applying \eA twice we obtain $(*)$. \qth

\subsection{Solution}
Now we must check this $\eta: \im \theta \to \chi$ indeed satisfies $\eta \circ \theta = u$. The composition on the left hand side is\[
\exists W_{k}\rak{\theta(x_{-i}, W_{k})\wedge \exists R_{i} \rak{\theta(R_{i}, W_{k})\wedge u(R_{i}, x_{w})}} \tag{1}
\]
and we must show this is provably equivalent ($T \cup \{\neg E\}$-provably equivalent) to $u(x_{-i}, x_{w})$. It is clearly provably equivalent to \[
\exists R_{i}\rak{u(R_{i}, x_{w})\wedge \exists W_{k}\rak{\theta(x_{-i}, W_{k})\wedge \theta(R_{i}, W_{k})}}.\tag{1\pr}
\]

To obtain \[
\Gamma, \neg E, (1\pr) \Rightarrow u(x_{-i}, x_{w}) \tag{*a}
\]
we note that (1\pr) has a subformula \[
\exists W_{k}\rak{\theta(x_{-i}, W_{k})\wedge \theta(R_{i}, W_{k})} \tag{2}
\]

which is $\FP(x_{-i}, R_{i})$. By \frelim on the output of \cref{dec-mapslemma}, this provably implies $\exists N_{w}\rak{u(x_{-i}, N_{w})\wedge u(R_{i}, N_{w})}$.
By the Implication Lemma (\cref{dec-implemma} in Chapter 1; clearly we also have a forward direction for $\wedge$) we conclude that $(1\pr)$ provably implies
\[
\exists R_{i}\rak{u(R_{i}, x_{w})\wedge \exists N_{w}\rak{u(x_{-i}, N_{w})\wedge u(R_{i}, N_{w})}}.\tag{2}
\] By Extraction and then the Implication Lemma, we see that (2) provably implies\[
\exists R_{i}\rak{ u(x_{-i}, x_{w})} \tag{3}
\](this is an application of Extraction where the formula $\gamma = u(x_{-i}, N_{w})$). (3) is provably equivalent to just $u(x_{-i}, x_{w})$, so we have the provable implication (*a).\\

Conversely for \[
\Gamma, \neg E, u(x_{-i}, x_{w}) \Rightarrow (1\pr)\tag{*b}
\] there is an obvious witness for $R_{i}$, namely $x_{-i}$. Rigorously, we have \[
\Gamma, \neg E, u(x_{-i}, x_{w}) \Rightarrow u(x_{-i}, x_{w}) \wedge \exists W_{k}\rak{\theta(x_{-i}, W_{k})\wedge \theta(x_{-i}, W_{k})} \tag{4}
\] by $\fone^{u}$ and $\ftwo^{\theta}$; and the RHS is equal to \[
\paren{u(R_{i}, x_{w})\wedge \exists W_{k}\rak{\theta(R_{i}, W_{k})\wedge \theta(R_{-i}, W_{k})}}\sbb{x_{-i}}{R_{i}}
\] so we can apply \eS to (4) to get (*b). \qed\\

This shows $\eta \circ \theta = u$ as morphisms. \qed\\

\subsection{Uniqueness}
Now suppose $\eta': \im \theta \rightarrow \chi(Q_{w})$ is such that $\eta' \circ \theta = u$. We shall show derivability of $\E^{\eta \to \eta'}$ and $\E^{\eta' \to \eta}$. For the former we must derive
\[
\Gamma, \neg E, \exists R_{i}\rak{\theta(R_{i},x_{-k})\wedge u(R_{i},x_{w})}\Rightarrow\eta'(x_{-k},x_{w}).\tag{*}
\]

As $\eta'\circ \theta = u$, we have by \frelim on $\E^{\eta' \circ \theta \to u}$ (choosing to instantiate the $x_{w}$ in the usual form of this sentence as $x_{c + w}$)
\[
\Gamma, \neg E, \exists W_k\rak{\theta(R_{i},W_k)\wedge\eta'(W_k,x_{c+w})}\Rightarrow u(R_{i},x_{c+w})\qquad \qquad \tag{1}
\]\[
\qquad \overline{\Gamma, \neg E, \theta(R_{i},x_{-k}),\eta'(x_{-k},x_{c + w})\Rightarrow u(R_{i},x_{c + w})} \qquad (\eelim)\tag{1\pr}
\]

By \frelim on $\fth^{u}$, we have\[
\Gamma, \neg E, u(R_{i}, x_{w}), u(R_{i}, x_{c + w})\Rightarrow \bigwedge_{w} x_{w} = x_{c + w}\qquad \qquad \qquad \qquad\tag{3}
\]
\[
\overline{\Gamma, \neg E, \theta(R_{i},x_{-k}),\eta'(x_{-k},x_{c + w}), u(R_{i},x_{w})\Rightarrow\bigwedge_{w} x_{w}=x_{c + w}} \qquad (\ch \text{ with 1\pr})\tag{3\pr} 
\] and by chaining with instances of \subs this gives
\[
\Gamma, \neg E, \theta(R_{i},x_{-k}),\eta'(x_{-k},x_{c + w}), u(R_{i},x_{w}) \Rightarrow \eta'(x_{-k}, x_{w})\qquad \qquad \qquad \qquad \tag{3\pr\pr}
\]\[
\overline{\Gamma, \neg E, u(R_{i},x_{w}),\theta(R_{i},x_{-k}),\exists x_{c + w} \eta'(x_{-k},x_{c + w})\Rightarrow \eta'(x_{-k}, x_{w})}  \qquad \qquad(\eA) \tag{3 \pr \pr \pr}
\]

Now by $\fone^{\theta} $ and $\ftwo^{\eta'}$, we have\[
\Gamma, \neg E, \theta(R_{i}, x_{-k}) \Rightarrow \exists x_{c + w}\eta'(x_{-k},x_{c + w})\tag{2}
\]

\[(2)(3''')
\]\[
\overline{\Gamma, \neg E, u(R_{i},x_{w}),\theta(R_{i},x_{-k}) \Rightarrow \eta'(x_{-k}, x_{w})}\qquad (\ch)
\]\[
\overline{\Gamma, \neg E, \exists R_{i}\rak{\theta(R_{i},x_{-k})\wedge u(R_{i},x_{w})} \Rightarrow \eta'(x_{-k}, x_{w})},\qquad (\eA)
\] 

which is (*).\qed \\

For the converse sentence $\E^{\eta' \to \eta}$, we need
\[
\Gamma, \neg E, \eta'(x_{-k},x_{w})\Rightarrow\exists R_{i}\rak{\theta(R_{i},x_{-k})\wedge u(R_{i},x_{w})} \tag{*}
\]

By \frelim and then \eelim on $\E^{\eta' \circ \theta \to u}$, as for (1) and (1\pr) (but this time instantiating the bound variables $x_{w}$ as themselves in the \frelim step), we have
\[
\Gamma, \neg E,\theta(R_{i},x_{-k}),\eta'(x_{-k},x_{w})\Rightarrow u(R_{i},x_{w})
\]\[
\overline{\Gamma, \neg E,\theta(R_{i},x_{-k}),\eta'(x_{-k},x_{w})\Rightarrow u(R_{i},x_{w})\wedge\theta(R_{i},x_{-k})}
\]\[
\qquad \overline{\Gamma, \neg E,\exists R_{i}\rak{\theta(R_{i},x_{-k})}, \eta'(x_{-k},x_{w})\Rightarrow\exists R_{i}\rak{u(R_{i},x_{w})\wedge\theta( R_{i}, x_{-k})}} \qquad \qquad (\eD)
\]

Now by $\fone^{\eta'}$ we have 
\[
\Gamma, \neg E, \eta'(x_{-k},x_{w})\Rightarrow\exists R_{i}\theta(R_{i},x_{-k}).
\] 
Chaining this with the above, we reduce to only the $\eta'(x_{-k},x_{w})$ hypothesis on the LHS and thus have (*).\qed\\

Thus we conclude any such $\eta'$ is provably equivalent to the $\eta$ we have given, and uniqueness holds in this subcase.\qed\\

This concludes the verification of the universal property of the coequalizer for $\theta: \alpha \to \im \theta$, and hence the proof of \cref{dec-toimageisee}. \qed\\

\subsection{Factorization}
We now have our Case 0 factorization result. Namely, by the results above, for any $\theta: \alpha(X_{i})\to \gamma (Z_{k})$ we have an effective epimorphism $\theta: \alpha(X_{i}) \to (\im \theta)(Z_{k})$. By the definition of $\im \theta$, we have an obvious special monomorphism $\im \theta\to \gamma$. The composite of these is the original morphism $\theta$, since postcomposition with special monomorphisms preserves representing formulas. \qed\\

We note that the factorization here is in fact into an e.e. followed by a {\it special} monomorphism. We shall define the intermediate object $\im f$ in the factorization for other types of morphisms $f$, in such a way that this continues to hold.

\section{Other cases}
For a type 4 morphism $s\to s'$, we simply take $s\to s \to s'$. The first map is automatically e.e. as it is an identity, and the second is a special mono as are all maps of sentences.\qed\\

For a type 3 morphism $s \overset{\theta(x_{i})}{\rightarrow} \varphi(X_{i})$, we take $s \overset{\theta(x_{i})}{\rightarrow} \theta(X_{i}) \rightarrow \varphi(X_{i})$. By \cref{dec-isosingular} in Chapter 4, the first map is an isomorphism, hence automatically e.e. We clearly have the second special monomorphism (by $\mfone^{\theta}$).\qed \\

For type 2 $f: \varphi(X_{i})\to s$, we take $\im f$ to be the sentence $\exists X_{i} \varphi$. Then we indeed have a type 2 (formal) map $\varphi \to \im f$, as well as formal $\im f \to s$, and the composite can only be the original map. The type 4 map is a special mono, and what we must show is that $\varphi(X_{i})\overset{p}{\rightarrow} \exists X_{i} \varphi(X_{i})$ is e.e., for any $\varphi(X_{i})$.

That $p$ equalizes the projections $\pi_{1}$, $\pi_{2}$ from its self-pullback, is immediate as it is type 2. Moreover if $\varphi(X_{i})\to s'$ is any other type 2 map out of $\varphi$, it too equalizes $\pi_{1}, \pi_{2}$. By \eA we clearly have a type 4 map $\exists X_{i} \varphi(X_{i}) \to s'$, which is unique and must factorize $\varphi(X_{i}) \to s'$. So the universal property of the coequalizer for $p$ is satisfied for the trivial case of $\chi = s'$, as in the proof of \cref{dec-toimageisee}.

Now we argue, somewhat oddly, as follows. If there are no type 1 maps $u: \varphi(X_{i})\to \chi(Q_{w})$ such that $u \circ \pi_{1} = u \circ \pi_{2}$ at all, then the verification of the universal property of the coequalizer is done. Otherwise, let $u = u(x_{-i},x_{w})$ be such a map (for some $\chi$). We claim a diagram \begin{center}\begin{tikzcd}
&\varphi(X_{i})\arrow[ld, "u"]\arrow[rd, "p"]&\\
\im u\arrow[rr, "\iota"] & & \exists X_{i}\varphi
\end{tikzcd}\end{center}
where $\iota$ is a type 2 {\it isomorphism}. By Case 0, the left leg $u$ is a coequalizer of maps $\pi_{1}^{u}$ and $\pi_{2}^{u}$ (not shown), so it will follow that $p$ is also a coequalizer of these, and hence a {\it regular epimorphism} (\tkref nLab). Since our category has pullbacks, this will show that $p$ is an effective epimorphism (\tkref NLAB).\\

{\it Proof of diagram:} Indeed, we have existence of the type 2 bottom leg, by the universal property of $u: \varphi(X_{i})\to \im u$. By \cref{dec-t2iso} in the last Chapter, to show it is an isomorphism we need to show that \[
(\im u)(x_{w}) = \exists M_{i} u(M_{i}, x_{w}) \in \prehom(\exists X_{i}\varphi, \im u).
\]

Deriving $\mfone$ and $\mftwo$ is straightforward. For $\mfth$ we need (by alpha-equivalence)\[
\Gamma, \neg E, \exists x_{-i}u(x_{-i}, x_{w}), \exists x_{-n-i} u(x_{-n-i}, x_{c + w}) \Rightarrow \bigwedge_{w} x_{w} = x_{c + w}. \tag{*}
\]

According to the construction in Fiber Products, the domain $p \times p$ of $\pi_{1}$, $\pi_{2}$ can be taken to be formula $\FP = \varphi(X_{i}) \wedge \varphi(T_{i})$. Therefore by the postcomposition formulas we have\[
u \circ \pi_{1} = \FP(x_{-i}, x_{-n-i}) \wedge u(x_{-i}, x_{w}) = \varphi(x_{-i})\wedge \varphi(x_{-n-i}) \wedge u(x_{-i}, x_{w})
\]\[
u \circ \pi_{2} = \FP(x_{-i}, x_{-n-i}) \wedge u(x_{-n-i}, x_{w}) = \varphi(x_{-i})\wedge \varphi(x_{-n-i}) \wedge u(x_{-n-i}, x_{w})
\]
and by $\fone^{u}$ these can be simplified to
\[
u \circ \pi_{1} =   \varphi(x_{-n-i}) \wedge u(x_{-i}, x_{w})
\]\[
u \circ \pi_{2} = \varphi(x_{-i})\wedge u(x_{-n-i}, x_{w}).
\]

As these formulas are assumed to be provably equivalent, we have \[
\Gamma, \neg E, \varphi(x_{-n-i}), u(x_{-i}, x_{w}) \Rightarrow u(x_{-n-i}, x_{w}). \tag{1}
\]

Now by $\fth^{u}$ we have\[
\Gamma, \neg E, u(x_{-n-i}, x_{w}), u(x_{-n-i}, x_{c + w})\Rightarrow \bigwedge_{w} x_{w} = x_{c + w}. \tag{2}
\]

\[ (1) (2)
\]\[
\overline{\Gamma, \neg E, \varphi(x_{-n-i}), u(x_{-i}, x_{w}), u(x_{-n-i}, x_{c + w})\Rightarrow \bigwedge_{w} x_{w} = x_{c + w}} \qquad (\ch)\tag{3}
\]

We may drop the $\varphi(x_{-n-i})$ hypothesis in (3) since $u(x_{-n-i}, x_{c + w})$ already implies it by $\fone^{u}$; thus we have\[
\Gamma, \neg E, u(x_{-i}, x_{w}), u(x_{-n-i}, x_{c + w})\Rightarrow \bigwedge_{w} x_{w} = x_{c + w}. \tag{3\pr}
\]
Applying \eA twice we obtain (*). \qed\\ 

This concludes the proof of the claimed diagram, and hence as per the above discussion the proof that $\varphi(X_{i}) \to \exists X_{i} \varphi$ is e.e. \qed\\

\section{Summary}\label{seg-factsum}
We briefly summarize the results of this part. \begin{enumerate}
\item If $\alpha(X_{i})\overset{\theta}{\rightarrow} \gamma(Z_{k})$ is type 1, then it factors into an effective epi followed by a special mono, as \[
\alpha(X_{i}) \overset{\theta}{\rightarrow}\exists M_{i}\theta(M_{i}, Z_{k})\rightarrow \gamma(Z_{k}).
\]
\item If $\varphi(X_{i}) \to s$ is type 2, then it factors into an effective epi followed by a special mono, as\[
\varphi(X_{i})\rightarrow \exists X_{i} \varphi(X_{i})\rightarrow s.
\]
\item If $s \overset{\theta}{\rightarrow} \varphi(X_{i})$ is type 3, then it factors into an effective epi followed by a special mono, as\[
s \overset{\theta}{\rightarrow} \theta(X_{i})\rightarrow \varphi(X_{i}).
\]
\item If $s \rightarrow s'$ is type 4, then it factors into an effective epi followed by a special mono, as\[
s \rightarrow s \rightarrow s'.\]
\end{enumerate}
In all cases we call the intermediate object the {\it image} of the morphism.

\chapter{$\sub(X)$}\label{ch-subx}

If $X$ is an object in a small category, we can denote by $\sub(X)$ the set of equivalence classes of monomorphisms into $X$, where $Y \to X$ and $Y' \to X$ are identified iff there is an isomorphism $\iota$ making\begin{center}
\begin{tikzcd}
Y \arrow[rd] \arrow[rr, "\iota"] & & Y' \arrow[ld]\\
& X&
\end{tikzcd}
\end{center} commute (\cite{lu}, Lecture 3).

\newcommand{\set}{\textbf{Set}}
It is then easily checked that defining $Y \to X \leq Z \to X$ to hold iff there exists a map $Y\to Z$ commuting the triangle, we obtain a poset $(\sub(X), \leq)$. (Asymmetry of $\leq$ uses the fact that such a map is necessarily unique, and is a monomorphism.)

In our category $\syn(T)$, this poset is in fact a Boolean lattice (=complemented distributive lattice with least and greatest elements) for every object $X$, with l.u.b.'s and g.l.b.'s given by logical $\vee$ and $\wedge$ respectively; and pulling back along an arbitrary morphism $f: A\to X$ defines a structure-preserving map $\sub(X) \to \sub(A)$. We view this as an analogue of the statement (in $\set$) that inverse images commute with unions and intersections. We shall develop this here.

As usual we first consider the case where $X = \Pi(X_{i})$ is a nonsentence.\\

For the remainder of our paper if $V$ is a nonempty finite set of variables, $E^{V}$ will denote the formula $X_{1} = X_{1}\wedge \dots \wedge X_{n} = X_{n}$, where $X_{1}, \dots X_{n}$ is an enumeration of the elements of $V$ in increasing index order.

\section{Case 0}
The first observation is that morphism factorization gives us a standardized form for elements of $\sub(\Pi)$: they are all represented by special monomorphisms $\varphi(X_{i}) \rightarrowtail \Pi(X_{i})$. Indeed, if $\omega \overset{m}{\rightarrow} \Pi(X_{i})$ is any monomorphism, we have two factorizations \begin{center}
\begin{tikzcd}
& \omega \arrow[dl]\arrow[dr, "\text{id}"]&\\
(\im m)(X_{i})\arrow[dr, tail] & & \omega \arrow[dl, "m"]\\
&\Pi(X_{i})&
\end{tikzcd}
\end{center}
into an effective epi followed by a mono, so by uniqueness there must be an isomorphism $\im m \cong \omega$ such that the diagram commutes. This says that $m$ is equal as a subobject to the tailed special monomorphism.\\

Thus to specify an element of $\sub(\Pi)$, we will usually just specify a special monomorphism representative, that is, a formula $\varphi(X_{i})$ such that $T$ proves the sentence $\forall X_{i}\rak{\varphi \to \Pi }$. A defining formula for the special mono itself is always given by
\[
\iota_{\varphi} := \varphi(x_{-i})\wedge \bigwedge_{i} x_{i} = x_{-i}.
\]
Also note that this is a defining formula for any special mono out of $\varphi$, including the identity $I_{\varphi}$.\\

\begin{lem}\label{dec-lequal}Let $\varphi$, $\psi$ represent elements of $\sub(\Pi)$. Then we have $\varphi \leq \psi$ if and only if $T$ proves\[
\forall X_{i}\rak{\varphi \rightarrow \psi},
\]i.e. iff $\varphi$ provably implies $\psi$.\end{lem}
\pf Suppose $T$ proves this sentence. Then by definition we have a special mono $m: \varphi \to \psi$. We have commutativity of \begin{center}
\begin{tikzcd}
\varphi \arrow[rd] \arrow[rr, "m"] & & \psi \arrow[ld]\\
& \Pi&
\end{tikzcd}
\end{center}
because $m$ is defined by the same formula $\iota_{\varphi}$ as the left leg, and the right leg transfers this formula to the composite as it is a special monomorphism. Therefore, we have $\varphi \leq \psi$.

Conversely, suppose $\varphi \leq \psi$ holds; then by definition we have a map $m$ such that the above commutes. As special monos transfer incoming formulas, this implies $m$ is defined by the formula $\iota_{\varphi}$. But then by \frelim on $\fone^{m}$ we have\[
\Gamma, \neg E, \varphi(x_{-i})\wedge \bigwedge_{i} x_{i} = x_{-i} \Rightarrow \psi(x_{i})\qquad \qquad
\]\[
\overline{\Gamma, \neg E, \varphi(x_{-i}), \bigwedge_{i} x_{i} = x_{-i} \Rightarrow \psi(x_{-i})}\qquad (\subs)
\]\[
\overline{\Gamma, \neg E, \varphi(x_{-i}), \exists x_{i}\rak{\bigwedge_{i} x_{i} = x_{-i}} \Rightarrow \psi(x_{-i})}\qquad (\eA)
\]
and the last hypothesis is vacuous given $\neg E$, and can be dropped to obtain\[
\Gamma, \neg E, \varphi(x_{-i}) \Rightarrow \psi(x_{-i}).\qed
\]

\begin{cor}\label{dec-subobjs-provequiv} Let $\varphi$, $\psi$ be formulas admitting special monomorphisms to $\Pi$. Then they represent the same subobject of $\Pi$ iff they are provably equivalent. \end{cor}
\pf Straightforward from \cref{dec-lequal}.\qed\\

Combined with our earlier remarks, this shows that we may identify $\sub(\Pi)$ with the set of such formulas up to provable equivalence.

\subsection{Lattice structure}\label{seg-latticestruct}
\begin{clm}\label{dec-join} Under the above identification, any $\varphi, \psi\in \sub(\Pi)$ admit the least upper bound $\varphi \vee \psi$. \end{clm}

\pf $\varphi\vee \psi$ represents an element of $\sub(\Pi)$ by the left introduction rule for $\vee$:\begin{align*}
&\Gamma, \neg E, \varphi \Rightarrow \Pi\\
&\Gamma, \neg E, \psi \Rightarrow \Pi\\
&\overline{\Gamma, \neg E, \varphi \vee \psi \Rightarrow \Pi} \qquad (\oorl)
\end{align*}

By the right introduction rule (and \cref{dec-lequal}), we see $\varphi \leq \varphi \vee \psi$ and $\psi \leq \varphi \vee \psi$.

Now suppose $\varphi \leq \theta$, $\psi \leq \theta$. Then exactly as above we obtain the sequent $\Gamma, \neg E, \varphi \vee \psi \Rightarrow \theta$, and we therefore have $\varphi \vee \psi \leq \theta$. Hence, $\varphi \vee \psi$ is least.\qed\\

\begin{clm}\label{dec-meet}Any $\varphi, \psi$ admit the greatest lower bound $\varphi \wedge \psi$.\end{clm}

\pf Symmetric.\qed\\

\begin{clm}\label{dec-least} $\sub(\Pi)$ has the least element $0_{\Pi}:=\neg \E^{\free(\Pi)}$. \end{clm}
\pf By the rules \expl and \id, we have sequents with $0_{\Pi}$ as sole hypothesis and with arbitrary conclusion. It certainly follows that $0_{\Pi}$ represents a subobject under the above identification, and also that it is $\leq$ anything.\qed\\

\begin{clm}\label{dec-greatest}$\sub(\Pi)$ has the greatest element $\Pi$. \end{clm}
\pf Clearly a subobject, and by assumption for any $\varphi \in \sub(\Pi)$ we have a special mono $\varphi \rightarrow \Pi$, so we have $\varphi \leq \Pi$ by \cref{dec-lequal}. \qed\\

\begin{clm}\label{dec-distributive} The lattice $\sub(\Pi)$ is distributive. \end{clm}
\pf Since the join of $\alpha$ and $\beta$ is given by $\alpha \vee \beta$ and the meet by $\alpha \wedge \beta$, one of these is immediate from the provable equivalence of $\alpha \wedge (\beta \vee \gamma) $ with $(\alpha \wedge \beta)\vee (\alpha \wedge \gamma)$, and similarly for the other way. \qed\\

\begin{clm}\label{dec-comp} Any $\varphi \in \sub(\Pi)$ admits the complement $\Pi \wedge \neg \varphi$.\end{clm}

\pf We have that $\varphi \vee (\Pi \wedge \neg \varphi)$ is equal as a subobject to $(\varphi \vee \Pi)\wedge (\varphi\vee \neg \varphi) = \varphi \vee \Pi = \Pi$. And $\varphi \wedge (\Pi \wedge \neg \varphi)$ proves anything, by \expl, so it is $\leq 0_\Pi$, so it is equal to it. \qed\\

\section{Case 1}\label{seg-subasentence}
Next we consider the poset $\sub(s)$, $s$ a sentence. Again, factorization gives us that every subobject is represented by a special monomorphism, which in this case is a type 4 map $s' \to s$. If $s''\to s$ is any other type 4 map, by definition we have a map $s' \to s''$ iff $s'$ $T$-provably implies $s''$, and in this case the triangle automatically commutes. Thus, we see that the relation $\leq$ on subobjects is given by $T$-provable implication (whereas it was $T \cup \crul{\neg E}$-provable implication in Case 0), and $\sub(s)$ is identified with the set of sentences which $T$-provably imply $s$, mod $T$-provable equivalence.

All the lattice structure claims are identical, and are proven simply without using \frelim and hence without the $\neg E$ hypothesis in the sequents. The only change is that the least element $0_{s}$ is now represented by the contradictory sentence $\exists x \rak{\neg x = x}$, for which we again have arbitrary-RHS sequents. \qed

\section{Pullback homomorphism of Boolean lattices}
If $f: A\to X$ is any morphism in a small category with pullbacks, it can be checked that the operation of forming a pullback along $f$ yields a well-defined map of sets $f^{-1}: \sub(X) \to \sub(A)$. (Monomorphisms pull back to monomorphisms, and pullbacks are unique up to an isomorphism {\it over the two domain objects}. If we vary a mono to $X$ by any source iso before pulling it back, by abstract nonsense we need not change the pullback along $f$.) In fact, $f^{-1}$ is a map of posets (\cite{lu}, Lecture 4, p.4). 

In $\syn(T)$, this map also respects all the lattice structure. To show this, it suffices to check preservation of meets, joins, and greatest and least elements.

\subsection{Type 1}
First we consider the case of a type 1 morphism $f$, so let $\theta(x_{-j}, x_{i}): \alpha(Y_{j})\to \Pi(X_{i})$. We would like to find a special mono representative in $\sub(\alpha)$ for the projector $\pi: \FP \to \alpha$, where $\FP = \theta \times \iota_{\varphi}$ is the fiber product of $\theta$ with some special mono $\varphi(X_{i})\rightarrowtail \Pi(X_{i})$. By our earlier factorization trick, we can simply take $\im \pi \rightarrowtail \alpha$. We know $\pi$ is defined by the formula $\FP(x_{-i}, x_{-n-j})\wedge \bigwedge_{j} x_{j} = x_{-n-j}$, so by the constructions in the last Chapter $\im \pi$ is the formula\[
\exists M_{i} \exists M_{n + j}\rak{\pi(M_{i}, M_{n + j}, Y_{j})}
\]\[
=\exists M_{i} \exists M_{n + j} \rak{\FP(M_{i}, M_{n + j})\wedge \bigwedge_{j}Y_{j} = M_{n + j}},
\]where $(M_{s})$ are unused variables. This is provably equivalent to (and hence the same subobject of $\alpha(Y_{j})$ as)
\[\exists M_{i} \exists M_{n + j} \rak{\FP(M_{i}, Y_{j})\wedge \bigwedge_{j}Y_{j} = M_{n + j}}
\]\[\equiv\exists M_{i}\rak{ \FP(M_{i}, Y_{j})\wedge \exists M_{n + j} \rak{\bigwedge_{j}Y_{j} = M_{n + j}}} 
\]\[\label{eq-impi}
\equiv \exists M_{i}\FP(M_{i}, Y_{j}),\tag{1}
\] as the second clause is vacuous given $\neg E$. Writing out (1), we find it is\[
\exists M_{i}\rak{\exists W_{i}\rak{\iota_{\varphi}(M_{i}, W_{i})\wedge \theta(Y_{j}, W_{i})}}
\]
\[=\exists M_{i} \exists W_{i}\rak{\varphi(M_{i})\wedge \bigwedge_{i} W_{i} = M_{i}\wedge \theta(Y_{j}, W_{i})}.
\]By a similar manipulation to above this is $\equiv$ to\[
\exists M_{i}\rak{\varphi(M_{i})\wedge \theta(Y_{j}, M_{i})}.
\]\\

\begin{rk} \label{dec-pullbackrep}This is the special mono representative we shall use for $f^{-1}(\varphi)$. Thinking set-theoretically, we see it should indeed define the preimage of the subset $\varphi$ under $\theta$. We can now straightforwardly check the structure preservation.\end{rk}

\bigskip
\bigskip

For join, we need the provable equivalence \[\exists M_{i}\rak{\varphi(M_{i})\wedge \theta(Y_{j}, M_{i})} \vee \exists M_{i}\rak{\psi(M_{i})\wedge \theta(Y_{j}, M_{i})}\equiv \exists M_{i}\rak{(\varphi(M_{i})\vee \psi(M_{i}))\wedge \theta(Y_{j}, M_{i})}\tag{*a}\] for any $\varphi$, $\psi$ in $\sub(\Pi)$. Indeed, by distributivity the RHS is p.e. to $$\exists M_{i}\rak{(\varphi(M_{i})\wedge \theta(Y_{j}, M_{i})) \vee (\psi(M_{i})\wedge \theta(Y_{j}, M_{i}))}$$
and this is p.e. to the LHS by \sep and \unsep.\qed\\

For meet, we need \[\exists M_{i}\rak{\varphi(M_{i})\wedge \theta(Y_{j}, M_{i})} \wedge \exists M_{i}\rak{\psi(M_{i})\wedge \theta(Y_{j}, M_{i})}\equiv \exists M_{i}\rak{\varphi(M_{i})\wedge \psi(M_{i})\wedge \theta(Y_{j}, M_{i})}.\tag{*b}\] Indeed, \[\Gamma, \neg E, (\RHS) \Rightarrow (\LHS)\]
is immediate by \ex and the Implication Lemma. Conversely, as $\theta$ is functional, by the Extraction Lemma we have\[
\Gamma, \neg E, \theta(Y_{j}, M_{i}), \exists W_{i}\rak{\varphi(W_{i})\wedge \theta(Y_{j}, W_{i})}\Rightarrow \varphi(M_{i})\tag{1}
\]\[
\Gamma, \neg E, \theta(Y_{j}, M_{i}), \exists W_{i}\rak{\psi(W_{i})\wedge \theta(Y_{j}, W_{i})}\Rightarrow \psi(M_{i})\tag{2}
\]

\[
(1) (2)
\]\[
\overline{\Gamma, \neg E, \theta(Y_{j}, M_{i}),  \exists W_{i}\rak{\varphi(W_{i})\wedge \theta(Y_{j}, W_{i})}\wedge \exists W_{i}\rak{\psi(W_{i})\wedge \theta(Y_{j}, W_{i})}\Rightarrow \varphi(M_{i})\wedge \psi(M_{i})}\qquad (\andd)
\]\[
\overline{\Gamma, \neg E, \theta(Y_{j}, M_{i}),  \exists W_{i}\rak{\varphi(W_{i})\wedge \theta(Y_{j}, W_{i})}\wedge \exists W_{i}\rak{\psi(W_{i})\wedge \theta(Y_{j}, W_{i})}\Rightarrow \varphi(M_{i})\wedge \psi(M_{i})\wedge \theta(Y_{j}, M_{i})}
\]\[
\overline{\Gamma, \neg E, \exists M_{i}\theta(Y_{j}, M_{i}),  \exists W_{i}\rak{\varphi(W_{i})\wedge \theta(Y_{j}, W_{i})}\wedge \exists W_{i}\rak{\psi(W_{i})\wedge \theta(Y_{j}, W_{i})}\Rightarrow (\RHS)}\qquad (\eD)
\]

We may drop the  $\exists M_{i}\theta(Y_{j}, M_{i})$ hypothesis, as the others clearly already imply it, and we get $\Gamma, \neg E, (\LHS) \Rightarrow (\RHS)$. \qed\\

For preservation of greatest element, we need\[
\alpha(Y_{j})\equiv \exists M_{i}\rak{\Pi(M_{i})\wedge \theta(Y_{j}, M_{i})} \tag{*c}\]
By $\fone^{\theta}$ we have \[
\Gamma, \neg E, \theta(Y_{j}, M_{i}) \Rightarrow \alpha(Y_{j})\wedge \Pi(M_{i})\tag{1}
\]\[
\overline{\Gamma, \neg E, \theta(Y_{j}, M_{i}) \Rightarrow \Pi(M_{i})\wedge \theta(Y_{j}, M_{i})} \tag{1\pr}
\]\[
\overline{\Gamma, \neg E, \exists M_{i}\theta(Y_{j}, M_{i}) \Rightarrow \exists M_{i}\rak{\Pi(M_{i})\wedge \theta(Y_{j}, M_{i})}}\qquad (\eD)\tag{1\pr \pr}
\]and by $\ftwo^{\theta}$ we have \[
\Gamma, \neg E, \alpha(Y_{j})\Rightarrow \exists M_{i }\theta(Y_{j}, M_{i})
\] so applying \ch we obtain $\Gamma, \neg E, (\LHS)\Rightarrow (\RHS)$. Conversely, (1) also gives\[
\Gamma, \neg E, \theta(Y_{j}, M_{i}) \Rightarrow \alpha(Y_{j})
\]\[
\overline{\Gamma, \neg E, \Pi(M_{i})\wedge \theta(Y_{j}, M_{i}) \Rightarrow \alpha(Y_{j})}
\]\[
\qquad \qquad \qquad \overline{\Gamma, \neg E, (\RHS) \Rightarrow (\LHS)}\qquad (\eA)
\]so we have the provable equivalence (*c). \qed\\

For preservation of least element, we need \[
\neg \E^{\free(\alpha)} \equiv \exists M_{i} \rak{\neg \E^{\free(\Pi)}(M_{i})\wedge \theta(Y_{j}, M_{i})} \tag{*d}
\] It suffices to show that the RHS is also contradictory, so that we have explosion given it as hypothesis. Indeed, let $(P_{i})$ be unused variables; then we have\[
\Rightarrow \E^{\free(\Pi)}(M_{i}) \qquad (\id)
\]\[
\overline{ \E^{\free(\Pi)}(P_{i}) \Rightarrow \E^{\free(\Pi)}(M_{i})} \qquad (\ant)
\]\[
\overline{ \neg \E^{\free(\Pi)}(M_{i}) \Rightarrow \neg \E^{\free(\Pi)}(P_{i})} \qquad (\cp)
\]\[
\overline{ \neg \E^{\free(\Pi)}(M_{i})\wedge \theta(Y_{j}, M_{i}) \Rightarrow \neg \E^{\free(\Pi)}(P_{i})} \qquad (\ant)
\]\[
\overline{ \exists M_{i}\rak{\neg \E^{\free(\Pi)}(M_{i})\wedge \theta(Y_{j}, M_{i})} \Rightarrow \neg \E^{\free(\Pi)}(P_{i})} \qquad (\eA)
\]\[
\overline{\E^{\free(\Pi)}(P_{i}) \Rightarrow \neg (\RHS)}\qquad (\cp)
\]\[
\qquad \qquad \qquad \qquad \qquad \overline{\Rightarrow \neg (\RHS)} \qquad (\ch \text{ with $\Rightarrow \E^{\free(\Pi)}(P_{i})$}).
\] \qed

This concludes the proof that $\varphi \mapsto f^{-1}(\varphi)$ is a homomorphism of Boolean lattices for type 1 maps $f$. \qed

\subsection{Type 3}
Now suppose the map we want to pull back along is type 3, i.e. we are given $f: s \overset{\theta(x_{i})}{\rightarrow} \Pi(X_{i})$ and need to consider $f^{-1}: \sub(\Pi) \to \sub(s)$. We can deal with this abstractly, because by last Chapter, we can factor $f$ into an isomorphism $g$ followed by the special mono $h: \theta(X_{i})\rightarrow \Pi(X_{i})$. Then by the type 1 case, $h^{-1}: \sub(\Pi)\to \sub(\theta)$ is a lattice homomorphism, and by abstract nonsense the original pullback map is given by taking $h^{-1}$ and then composing with the inverse of $g$. Any iso induces an isomorphism of posets via composition, so this exhibits the original pullback map $f^{-1}$ as a composition of lattice homomorphisms. Hence it is a lattice homomorphism.\qed

\subsection{Type 2}
Now suppose we want to pull back along a type 2 map $\alpha(Y_{j}) \to s$. According to the construction of \cref{seg-sentoversent}, a fiber product with $s'\in \sub(s)$ is given by $\alpha(Y_{j})\wedge s'$, equipped with the evident special monomorphism to $\alpha(Y_{j})$. We easily check preservation of \begin{itemize}
\item Joins: $\alpha(Y_{j})\wedge(s' \vee s'')$ is p.e. to $(\alpha(Y_{j})\wedge s')\vee (\alpha(Y_{j})\wedge s'')$ by distributivity
\item Meets: $\alpha(Y_{j})\wedge s' \wedge s''$ is p.e. to $\alpha(Y_{j})\wedge s'\wedge \alpha(Y_{j}) \wedge s''$
\item Least elements: $\alpha(Y_{j})\wedge 0_{s}$ is certainly contradictory hence triggers explosion, just as $0_{\alpha}$ does, so these are p.e.
\item Greatest elements: $\alpha(Y_{j})\wedge s$ is p.e. to just $\alpha(Y_{j})$, since by existence of the type 2 map $\alpha(Y_{j})\to s$ we derive \[\Gamma, \neg E , \alpha(Y_{j})\Rightarrow s.
\]
\end{itemize}\qed

\subsection{Type 4}
Lastly, if the map to pull back along is type 4, everything works identically with the exception of not using \frelim and $\neg E$ in the last bullet above. \qed

\section{Status}
We have now verified coherence properties 1, 2, 3 and 5 (cf Chapter 4) for $\syn(T)$.

\chapter{Stability}\label{ch-stability}

The last axiom of a coherent category we must verify is that the collection of effective epis is stable under pullback. \\

\section{Effective epis revisited}\label{seg-eere}
With our factorization result, we can obtain an easy characterization of when an arbitrary morphism is an e.e. Namely, the property of a morphism being e.e. is preserved by isomorphisms of the target (e.e.'s are the same as regular epimorphisms, and an isomorphism of the target preserves the fact that the relevant diagram is colimit) so if the special monomorphism $\im f \rightarrowtail \text{cod} f$ is in fact an iso, then $f$ is an e.e. Conversely, if $f$ is an e.e. then we obtain that this morphism is an iso by uniqueness of factorization. So we conclude\\

 \begin{cor}\label{dec-eesy} In $\syn(T)$, the morphism $f:X \to Y$ is an e.e. if and only if the formulas $\im f$ and $Y$ are special-isomorphic, 
 in other words $T/T\cup\crul{\neg E}$-provably equivalent. \end{cor}

Of course, to conclude $f$ is an e.e. we only need really to check that the target formula $T/T\cup\crul{\neg E}$-provably implies $\im f$, since the other implication is automatic by construction. We will use this criterion to check that, if $f$ is an e.e. of types 1, 2, 3, or 4, any pullback of it is also e.e.\\

We let $\equiv$ denote special isomorphism.

\section{Type 1}
Let $\theta(x_{-i}, x_{k}): \alpha \to \gamma$ represent a type 1 effective epimorphism, and suppose we want to calculate the image $\im \pi$ of the pullback along some type 1 $\omega(x_{-j}, x_{k}): \beta \to \gamma$. By the calculation leading to \hyperref[eq-impi]{(\ref*{eq-impi}) in Chapter 7}, we see that it is special isomorphic to \[
\exists T_{i}\rak{\FP(T_{i,}Y_{j})}
\]\[
=\exists T_{i}\exists W_{k}\rak{\theta(T_{i}, W_{k})\wedge \omega(Y_{j}, W_{k})}.
\]
So, to check the pullback $\pi$ is an e.e., we must obtain\[
\Gamma, \neg E, \beta(Y_{j})\Rightarrow \exists T_{i}\exists W_{k}\rak{\theta(T_{i}, W_{k})\wedge \omega(Y_{j}, W_{k})}, \tag{*}
\] given that we have \[
\Gamma, \neg E, \gamma(W_{k})\Rightarrow \exists T_{i}\rak{\theta(T_{i}, W_{k})}. \tag{1}
\]
Indeed, from (1) by \andd we have\[
\Gamma, \neg E, \omega(Y_{j}, W_{k})\wedge \gamma(W_{k})\Rightarrow \omega(Y_{j}, W_{k})\wedge \exists T_{i}\rak{\theta(T_{i}, W_{k})}\qquad
\]\[
\qquad \overline{\Gamma, \neg E, \exists W_{k}\rak{\omega(Y_{j}, W_{k})\wedge \gamma(W_{k})}\Rightarrow \exists W_{k}\rak{\omega(Y_{j}, W_{k})\wedge \exists T_{i}\rak{\theta(T_{i}, W_{k})}}}\qquad (\eD)\tag{1\pr}
\]

By combination of $\ftwo^{\omega}$ and $\fone^{\omega}$, we clearly have \[
\Gamma, \neg E, \beta(Y_{j})\Rightarrow \exists W_{k}\rak{\omega(Y_{j}, W_{k})\wedge \gamma(W_{k})} \tag{2}
\]\[
(2),(1\pr)
\]\[
\overline{\Gamma, \neg E, \beta(Y_{j})\Rightarrow \exists W_{k}\rak{\omega(Y_{j}, W_{k})\wedge \exists T_{i}\rak{\theta(T_{i}, W_{k})}}} \qquad (\ch)
\] whence we get (*) by \ex and \reord.\qed

Now, if the map to pull $\theta$ back along is some type 3 $\omega(x_{k}):s \to \gamma$, we reduce to the above case by factoring $\omega$, and using again the fact that e.e.'s are preserved by isomorphisms of the target. Thus we conclude any pullback of $\theta$ is an e.e.\qed

\section{Type 3}
If $\theta(x_{k}): s\to \gamma$ is a type 3 effective epimorphism, the result follows from the above case, by the following. We can pass to a type 1 $\widetilde{\theta}$ by an isomorphism of the source, and this preserves the property of being e.e. by the equivalence of regular and effective. Thus by the above case we have an e.e. pullback of $\widetilde{\theta}$ along anything. But this is also a pullback of $\theta$ along the same thing. \qed

\section{Type 2}
Suppose $\varphi(X_{i}) \to s$ is a type 2 e.e., and we want to calculate $\im \pi$ for $\pi$ the pullback along some other type 2 $\psi(Y_{j})\to s$. According to \cref{seg-oversentence}, we have\[
\FP  = \varphi(T_{i})\wedge \psi(T_{n+j})
\]
and $\im \pi$ is\[
\exists T_{i} \exists T_{n + j} \rak{\FP \wedge \bigwedge_{j}Y_{j} = T_{n+j}}
\]
\[
=\exists T_{i} \exists T_{n + j} \rak{ \varphi(T_{i}) \wedge\psi(T_{n+j})\wedge \bigwedge_{j}Y_{j} = T_{n+j}}
\]
\[
\equiv \exists T_{i} \rak{ \varphi(T_{i}) \wedge\psi(Y_{j})\wedge \exists T_{n+j}\rak{\bigwedge_{j}Y_{j} = T_{n+j}}}
\]
\[
\equiv \exists T_{i} \rak{ \varphi(T_{i}) \wedge\psi(Y_{j})}.
\] 

Therefore to show $\pi$ is an e.e. we need \[
\Gamma, \neg E, \psi(Y_{j})\Rightarrow  \exists T_{i} \rak{ \varphi(T_{i}) \wedge\psi(Y_{j})}. \tag{*}
\]
Clearly it suffices by \ins to obtain\[
\Gamma, \neg E, \psi(Y_{j})\Rightarrow \exists T_{i}\rak{\varphi(T_{i})},
\]
that is (by the factorization constructions)\[
\Gamma, \neg E, \psi(Y_{j})\Rightarrow \im(\varphi \to s). \tag{*\pr}
\]
But indeed, by existence of the map $\psi\to s$ and \frelim we have\[
\Gamma, \neg E, \psi(Y_{j})\Rightarrow s \tag{1}
\]
and as $\varphi\to s$ is an e.e., also
\[\Gamma, s \Rightarrow \im(\varphi \to s). \tag{2}
\] by \cref{dec-eesy}. Chaining (1) and (2) we obtain (*\pr). \qed\\

For pulling back along a type 4 $s' \to s$, we have by \cref{seg-sentoversent} \[
\FP = \varphi(X_{i})\wedge s'
\]
and $\pi: \FP \to s'$ is type 2. Therefore it is an e.e. if we have\[
\Gamma, s' \Rightarrow \exists X_{i}\rak{\varphi(X_{i})\wedge s'}.\tag{*}
\]
As in the first subcase it suffices to obtain \[
\Gamma, s' \Rightarrow \im \varphi. \tag{*\pr}
\] But we have $\Gamma, s' \Rightarrow s$, and we apply chain with (2) as before. \qed

\section{Type 4}
Since the image of $s \to s'$ is just $s$ by definition, by the remarks in \cref{seg-eere} if this map is an e.e. then it is an isomorphism. But isomorphisms pull back to isomorphisms, which are e.e. \qed

\chapter{Semantics}\label{ch-semantics}

\newcommand{\set}{\textbf{Set}}
\newcommand{\calm}{\mathcal{M}}
\newcommand{\tf}{\widetilde{f}}
\newcommand{\tg}{\widetilde{g}}
\newcommand{\qstp}{$\qed_{(*\pr)}$}

\section{Recap}\label{seg-recap}
The syntactic category $\syn(T)$\footnote{In our primary reference what we are calling $\syn(T)$ is called the {\it weak syntactic category} $\syn_{0}(T)$. Thus our notation is somewhat unfortunate but it seemed senseless to carry over Lurie's as we will not get anywhere near the object, the pretopos completion, for which $\syn(T)$ is reserved in his notes.} of a first-order theory $T$ over a language $L_{T}$ containing only predicate symbols, is a category whose objects are the formulas of $L_{T}$ and morphisms are the ``$T$-provable definable maps'' between them. It has a structure axiomatized by the notion of a {\it coherent category}, which is defined by factorization, limit and subobject properties approximating those of the category of sets.

The purpose of the present paper is to show that G\"{o}del's Completeness Theorem from classical first-order logic can in principle be viewed as arising from this structure of $\syn(T)$ for general consistent $T$. This is to be done via the following\\

{\bf Theorem} (Deligne, stated in \cite{lu}). Let $\mathcal{C}$ be a small consistent coherent category. Then there exists a functor $F: \mathcal{C}\rightarrow \set$ such that\begin{enumerate}
\item $F$ preserves finite limits, including final objects
\item $F$ carries effective epimorphisms to surjections of sets
\item For every $X \in \mathcal{C}$, $F$ induces a homomorphism $\sub(X) \to \sub(F(X))$ of the join semilattices of subobjects $\sub(X)$, $\sub(F(X))$.\\
\end{enumerate}

In the foregoing chapters we constructed $\syn(T)$ and derived its coherence on the basis of first-order proof theory. Since G\"{o}del's theorem links syntax and semantics, this is a necessary step if we are to truly derive G\"{o}del from the above theorem. Here we undertake the second step: given our definition of $\syn(T)$ and, by Deligne's Theorem, a functor $F:\syn(T)\to \sett$ satisfying the above conditions, can we concretely describe a model of $T$ in the traditional sense? As elaborated in the notes \cite*{lu} upon which our project is based, in categorical logic, rather than viewing a model as a single set, one essentially identifies it with the more ``unbiased'' data of a functor satisfying the above conditions (and one identifies elementary embeddings with natural transformations). Nevertheless we do not find it entirely trivial to show this viewpoint is equivalent to the traditional one.

\section{Semantic setup}\label{seg-semsetup}
In proving the GCT we will make the following assumptions about the standard semantic definitions. We sometimes write $|P|$ for $\text{arity}(P)$, and $M^{V}$ for the set of functions from $V \to M$. If we are dealing with some function $f: V_{0}\to V_{1}$, we denote the corresponding ``currying map'' by $\Lambda: M^{V_{1}}\to M^{V_{0}}$. It is just given by precomposition with $f$.\\

\bigskip

\begin{uass}\label{dec-structurespecdby}
If $M$ is a nonempty set, then a structure $\calm$ with underlying set $M$ can be specified by a choice of
\begin{enumerate}
\item sets $P^{\calm}\subset M^{|P|}$ for each predicate symbol $P$ of $L_{T}$ with positive arity
\item sets $P^{\calm}\subset U$ for each $P$ of arity 0, where $U$ is a fixed but arbitrary singleton.
\end{enumerate}\end{uass}
\bigskip

\begin{uass} If $V$ is a nonempty finite set of variables and $\calm$ is a structure with nonempty underlying set $M$, then $M^{V}$ is identified, via the natural correspondence, with the set of usual assignments in $M$ modulo agreement on $V$. If $\varphi$ is a formula such that $\free(\varphi) = V$ exactly, then the usual recursive definition of satisfaction of $\varphi$ by an $\calm$-assignment, translates in the obvious way under this trivial identification to a definition of the relation $\lambda \models \varphi$ for $\lambda \in M^{V}$.
\end{uass}

\bigskip

\begin{defns} \hypertarget{semanticdefs}{} If $\calm$ is a structure specified via \cref*{dec-structurespecdby}, we define sets\begin{itemize}
\item $\calm(\varphi):= \crul{\lambda \in M^{\free(\varphi)}| \lambda \models \varphi}$, for each nonsentence $\varphi$
\item $\calm(s):= U/\emptyset$ (according as $\calm \models s$ or not), for each sentence $s$.
\end{itemize}
\end{defns}

\bigskip

\begin{uass}\label{dec-missuchthat}If $\calm$ is a structure specified via \cref*{dec-structurespecdby}, then we have \begin{enumerate}[label = (\alph*)]
\item \hypertarget{suchthat-a}{}for $|P| = n > 0$, $\calm(Px_{i_{1}}\dots x_{i_{n}}) = \crul{\lambda \in M^{\crul{x_{i_{1}}, \dots x_{i_{n}}}} | (\lambda(x_{i_{1}}), \dots \lambda (x_{i_{n}}))\in P^{\calm}}$

\item \hypertarget{suchthat-b}{}for $P$ of arity 0, $\calm(P) = P^{\calm}$

\item \hypertarget{suchthat-c}{}for $\alpha $ equal to $x_{i} = x_{i}$, $\calm(\alpha) = M^{\crul{x_{i}}}$

\item \hypertarget{suchthat-d}{}for $\alpha $ equal to $ x_{i} = x_{j}$, $i \neq j$, $\calm(\alpha) = \text{(diagonal)}\subset M^{\crul{x_{i}, x_{j}}}$

\item \hypertarget{suchthat-e}{}for nonsentences $\alpha$ and $\beta$, $\calm(\alpha \vee \beta) = \calm(\alpha) \times M^{\free(\beta)\setminus \free(\alpha)}\cup \calm(\beta) \times M^{\free(\alpha)\setminus \free(\beta)}$; where $\times$ between subsets of $M^{A}$ and $M^{B}$ for disjoint $A$ and $B$ denotes the natural product subset of $M^{A \cup B}$

\item \hypertarget{suchthat-f}{}for nonsentence $\alpha$ and sentence $s$, $\calm(\alpha \vee s) = \calm(s \vee \alpha) = \begin{cases}
\calm(\alpha) & \calm(s) = \emptyset\\ M^{\free(\alpha)} & \text{otherwise}
\end{cases}$

\item \hypertarget{suchthat-g}{}for sentences $s$ and $s'$, $\calm(s \vee s') = \begin{cases} \emptyset & \text{both of $\calm(s)$, $\calm(s')$ are $\emptyset$}\\ U &\text{otherwise} \end{cases}$

\item \hypertarget{suchthat-h}{}for nonsentence $\alpha$, $\calm(\neg \alpha) = $ complement of $\calm(\alpha)$ in $M^{\free(\alpha)}$

\item \hypertarget{suchthat-i}{}for sentence $s$, $\calm(\neg s) = $ complement of $\calm(s)$ in $U$

\item \hypertarget{suchthat-j}{}if $X\in \free(\psi)$ and $X$ is not the only free variable, \[\calm(\exists X \psi) = \crul{\lambda \in M^{\free(\psi)\setminus \crul{X}} | \lambda = \Lambda(\widetilde{\lambda}) \text{ for some $\widetilde{\lambda}$ in $\calm(\psi)$}},\] where the map for currying is inclusion $\free(\psi)\setminus \crul{X}  \hookrightarrow \free(\psi)$

\item \hypertarget{suchthat-k}{}if $X\in \free(\psi)$ is the only free variable, $\calm(\exists X \psi) = \emptyset/U$ (according as $\calm(\psi) = \emptyset$ or not)

\item	\hypertarget{suchthat-l}{}if $X\not\in \free(\psi)$, $\calm(\exists X \psi) = \calm(\psi).$
\end{enumerate}
\end{uass}

\section{Lemmas}
\begin{prop}[Substitution] \label{dec-choosypully}
Let $\varphi(X_{i})$ be a nonsentence, $V = \{Y_{j}\}$ be a finite set of variables, and $\tf: \free(\varphi) \to V$ be a map. Then we have a pullback square
\begin{center}
\begin{tikzcd}
\varphi(\tf(X_i)) \wedge E^{V}\arrow[d] \arrow[rr] & & E^{V} \arrow[d] \\
\varphi(X_i) \arrow[rr] & & E^{\free(\varphi)}          
\end{tikzcd}
\end{center}
where the horizontal maps are special monomorphisms and the vertical maps are $\tf$-special morphisms.\end{prop}

\newcommand{\ti}{\widetilde{i}}
\newcommand{\tit}{\widetilde{\iota_{t}}}
\newcommand{\til}{\widetilde{\iota_{l}}}

\pf We clearly have $S\tf: E^{V} \to E^{\free(\varphi)}$, as well as a special mono $S\ti: \varphi\to E^{\free(\varphi)}$. We transform our constructed pullback of these into the above square, by isomorphisms of the top-left corner.

So choose $(T_{s})_{1\leq s\leq n+m}$ to be increasing such that $T_{n + j} = Y_{j}$. By \cref{ch-fps} a pullback of $S\tf$ and $S\ti$ is given by\[
\FP := \exists W_{i}\rak{S\ti(T_{i}, W_{i})\wedge S\tf(T_{n+j}, W_{i})}
\]\[
= \exists W_{i}\rak{\varphi(T_{i})\wedge \bigwedge_{i} W_{i} = T_{i}\wedge E^{V}(T_{n+j})\wedge \bigwedge_{i} W_{i} = T_{n + f(i)}}
\]\[
= \exists W_{i}\rak{\varphi(T_{i})\wedge \bigwedge_{i} W_{i} = T_{i}\wedge E^{V}\wedge \bigwedge_{i} W_{i} = Y_{f(i)}},
\]
together with maps $\pi_{\varphi}$, $\pi_{E^{V}}$. The former of these is $\til$-special, where $X_{i} \overset{\til}{\mapsto} T_{i}$, and the latter is $\tit$-special, where $\tit: V \hookrightarrow \free(\FP)$ is inclusion. We see that the index maps for these are given by $\iota_{t}(j) = n + j$, $\iota_{l}(i) = i$. Therefore we have\[
\pi_{E^{V}} = \FP(x_{-i}, x_{-n -j}) \wedge \bigwedge_{j} x_{j} = x_{-n-j}
\]\[
\pi_{\varphi} = \FP(x_{-i}, x_{-n -j}) \wedge \bigwedge_{i} x_{i} = x_{-i}.
\]

Step 1: It is easy to see that the formula $\FP$ is provably equivalent to\[
\varphi(T_{i}) \wedge E^{V} \wedge \bigwedge_{i} T_{i} = Y_{f(i)}\tag{1}
\] which is p.e. to \[
\varphi(Y_{f(i)}) \wedge E^{V} \wedge \bigwedge_{i} T_{i} = Y_{f(i)}. \tag{1\pr}
\]
(Going from (1) to (1\pr) is reversible, even if $f$ is not injective, since \[
\varphi\sbb{T_{i}}{X_{i}}\sbb{Y_{f(i)}}{T_{i}}, \bigwedge_{i} T_{i} = Y_{f(i)}\Rightarrow \varphi\sbb{T_{i}}{X_{i}}
\]
is a valid instance of \subs.) 

All these formulas share the same free variables. Therefore denoting (1\pr) by $\FP'$, we may pass by a special isomorphism from our constructed pullback to
\begin{center}
\begin{tikzcd}
\FP' \arrow[d, "\pi_{\varphi}"] \arrow[rr, "\pi_{E^{V}}"] & & E^{V} \arrow[d] \\
\varphi \arrow[rr] & & E^{\free(\varphi)}          
\end{tikzcd}
\end{center}
and keep the same representing formulas for $\pi_{\varphi}$ and $\pi_{E^{V}}$ (although for convenience let us now consider the $\FP$ clauses in them to be replaced by $\FP'$).\\
\newcommand{\te}{\widetilde{\varepsilon}}

Step 2: Let our desired upper-left corner $\varphi(Y_{f(i)})\wedge E^{V}$ be $\Theta(Y_{j})$. Then clearly we have an $\te$-special map $\FP' \to \Theta$, where $\free(\Theta)\overset{\te}{\hookrightarrow} \free(\FP')$ is inclusion. We claim this map is an isomorphism.

Indeed, the extra free variables $T_{i}$ which we are trying to get rid of were chosen to have lower index than the $Y_{j}$, so the index map $\varepsilon$ acts as $\varepsilon(j) = n+j$. Therefore we have\[
S\te = \FP'(x_{-i}, x_{-n-j})\wedge \bigwedge_{j} x_{j} = x_{-n-j}.
\] By \cref{dec-t1isosuf} of \cref*{ch-lemmas} it suffices to negate the variables of this representative and show the result is a premorphism $\Theta \to \FP'$.

Writing this formula out we find it is equal to\[
\FP'(x_{i}, x_{n+ j})\wedge \bigwedge_{j} x_{-j} = x_{n + j} 
\]\[
= \varphi(x_{n + f(i)})\wedge E^{V}(x_{n + j})\wedge \bigwedge_{i} x_{i} = x_{n + f(i)} \wedge \bigwedge_{j} x_{-j} = x_{n + j} \tag{2}.
\]
Now by \subs we have provable equivalence of the formulas\[
\bigwedge_{i} x_{i} = x_{n + f(i)} \wedge \bigwedge_{j} x_{-j} = x_{n + j}, \qquad \qquad \bigwedge_{i} x_{i} = x_{-f(i)} \wedge \bigwedge_{j} x_{-j} = x_{n + j}
\] and it follows by more use of \subs 
 that a provably equivalent formula to (2) is\[
\varphi(x_{-f(i)})\wedge E^{V}(x_{-j})\wedge \bigwedge_{i} x_{i} = x_{-f(i)} \wedge \bigwedge_{j} x_{-j} = x_{n + j} \tag{2\pr}.
\]

Now, the $\FP'(x_{i}, x_{n+ j})$ (target) conclusion in F1 is built into the formula (2). On the other hand the $\Theta$ half is built into (2\pr). Therefore we have \fone.

$\ftwo$ is immediate from the (2\pr) form. So is \fth, as all positive variables are explicitly assigned to negative ones. Therefore, we have verified the conditions for a premorphism $\Theta \to \FP'$ for this common provable equivalence class.\qed \\

Step 3: We show applying $S\te: \FP' \to \Theta$ takes the pullback square we have at the end of Step 1 to the desired square. This consists of verifying composition relations for the top and left legs.

For the top leg, we must check that $(\text{Evident special mono $\Theta \to E^{V}$})\circ S\te = \pi_{E^{V}}$. Since postcomposition with special monos preserves incoming formula, this amounts to checking provable equivalence of the formulas $S\te$ and $\pi_{E^{V}}$. But they are identical. \qed

For the other, we let $\omega(x_{-j}, x_{i})$ be the left leg of the desired square, i.e. the evident $\tf$-special morphism $\Theta \to \varphi$. We must check $\omega \circ S\te = \pi_{\varphi}$. We shall use the formula for {\it precomposition with} a special map (\cref{dec-soutgoing} in \cref*{ch-lemmas}), remembering that $\varepsilon(j) = n + j$.

Indeed, according to this Proposition $\omega \circ S\te$ is\[
\FP'(x_{-i}, x_{-n-j}) \wedge \omega(x_{-\varepsilon(j)}, x_{i})
\]\[
=\varphi(x_{-n-f(i)})\wedge E^{V}(x_{-n-j})\wedge \bigwedge_{i} x_{-i} = x_{-n-f(i)} \wedge \omega(x_{-n-j}, x_{i})
\]\[
=\varphi(x_{-n-f(i)})\wedge E^{V}(x_{-n-j})\wedge \bigwedge_{i} x_{-i} = x_{-n-f(i)} \wedge \paren{\Theta(x_{-j})\wedge \bigwedge_{i} x_{i} = x_{-f(i)}}\sbx{x_{-n-j}}{x
_{-j}}{x_{i}}{x_{i}}
\]\[
=\varphi(x_{-n-f(i)})\wedge E^{V}(x_{-n-j})\wedge \bigwedge_{i} x_{-i} = x_{-n-f(i)} \wedge \Theta(x_{-n-j})\wedge \bigwedge_{i} x_{i} = x_{-n-f(i)}
\]
The $\Theta(x_{-n-j})$ clause is redundant, being exactly the conjunction of the first two clauses, so this is just\[
\varphi(x_{-n-f(i)})\wedge E^{V}(x_{-n-j})\wedge \bigwedge_{i} x_{-i} = x_{-n-f(i)} \wedge \bigwedge_{i} x_{i} = x_{-n-f(i)}, \qquad\text{or}
\]
\[\FP'(x_{-i}, x_{-n-j}) \wedge \bigwedge_{i} x_{i} = x_{-n-f(i)}
\]
as a morphism. From this form we see it is indeed provably equivalent to $\pi_{\varphi}$, the only difference being which of the already equal $x_{-i}, x_{-n-f(i)}$ the last clause asserts $x_{i}$ to be.\qed\\

This concludes the proof of \cref*{dec-choosypully}.\qed\\

{\bf Proposition 1.5} (\cref*{dec-choosypully} for empty $V_{0}$).
Suppose $s$ is a sentence and $V$ is a nonempty finite set of variables. Then we have a pullback square\begin{center}
\begin{tikzcd}
s \wedge E^{V} \arrow[r]\arrow[d] & E^{V} \arrow[d]\\
s \arrow[r] & S
\end{tikzcd}
\end{center}

{\it Proof.} This is the construction of the fiber product given in \cref{seg-sentoversent} of \cref*{ch-fps}.\qed\\

\begin{prop}\label{dec-productiterator} Suppose $\varphi$, $\psi$ are nonsentences such that all free variables of $\psi$ are greater index than those of $\varphi$. Then

\begin{center}
\begin{tikzcd}
        & \varphi \wedge \psi \arrow[ld, "\pi_{\varphi}"'] \arrow[rd, "\pi_\psi"] &      \\
\varphi &                                                                     & \psi
\end{tikzcd}
\end{center}

is a product, where $\pi_{\varphi}, \pi_{\psi}$ are inclusion-special.\end{prop}
\pf Immediate from the verification of \cref{seg-oversentence} in \cref*{ch-fps}. \qed\\

\begin{cor}\label{dec-nfoldproduct} Let $V = \crul{X_{i}}$ be a nonempty finite set of variables. Then we have a diagram
\begin{center}\begin{tikzcd}
E^V \arrow[dd, "\pi_{1}", shift right=10mm] \arrow[dd, "\pi_{n}", shift left=10mm]\\
\dots \\
x_{0} = x_{0}                                                                            
\end{tikzcd}\end{center}
which exhibits $E^V$ as an $n$-fold product of $x_{0} = x_{0}$, where $\pi_{l}$ is defined by the formula $x_{-l} = x_{1}$.\end{cor}

\pf By iteration of \cref*{dec-productiterator}, we have that
\begin{center}
\begin{tikzcd}
        & E^V \arrow[ld, "S\ti_1"'] \arrow[rd, "S\ti_{n}"] &      \\
X_{1} = X_{1} &                                                 \dots                    & X_{n} = X_{n} 
\end{tikzcd}
\end{center} 
is a product. Composing with the obvious isomorphisms $X_{l} = X_{l}\to x_{0} = x_{0}$, we obtain an $n$-fold product of the latter formula with itself. By the postcomposition formula for special morphisms, the projection maps $\pi_{l}$ for this are given by\[
\pi_{l} = (E^V)(x_{-i})\wedge (x_{-1} = x_{1})(x_{-\iota_{l}(1)}, x_{1})
\]\[
= (E^V)(x_{-i})\wedge x_{-l} = x_{1}
\]
which as a morphism is the same as $x_{-l} = x_{1}$, as claimed, since the first clause is trivial to derive.\qed\\

\section{Model setup}\label{seg-modset}
Let us begin construction of the model $\calm$ for our consistent theory $T$ demanded by the GCT. By Deligne's Theorem, let $F:\syn(T)\to \set$ be a functor satisfying conditions 1, 2 and 3 of \cref{seg-recap}. For now let us assume that the set $M:=F(x_{0} = x_{0})$ is nonempty.

We need one final lemma to get started.
\subsection{Last lemma}
\begin{clm} \label{dec-lastlem}Let $\tf: V_{0} \to V_{1}$ be a map of nonempty finite sets of variables, and let $V_{0} = \crul{X_{i}}$, $V_{1} = \crul{Z_{j}}$ in increasing order. Then we have a commutative diagram
\begin{center}
\begin{tikzcd}
F(E^{V_{1}}) \arrow[rr, "b_{1}"] \arrow[d, "FS\tf"]& & M^{V_{1}}\arrow[d, "\Lambda \hspace{1mm} (\lambda\mapsto \lambda \circ \tf)"]\\
F(E^{V_{0}}) \arrow[rr, "b_{0}"] & & M^{V_{0}}
\end{tikzcd}
\end{center}
in \set, where $b_{0}$ and $b_{1}$ are bijections depending respectively only on $V_{0}, V_{1}$.\end{clm}
\pf The functor $F$ preserves products by condition 1, so applying it to the output of \cref*{dec-nfoldproduct} we have products

\begin{center}
\begin{tikzcd}
F(E^{V_{0}}) \arrow[d, "F\pi_{i}"]\\
M                                  
\end{tikzcd},
\begin{tikzcd}
F(E^{V_{1}}) \arrow[d, "F\varepsilon_{j}"]\\
M                    
\end{tikzcd}.
\end{center}

We also have products
\begin{tikzcd}
M^{V_{0}} \arrow[d, "p_{i}"]\\
M                                  
\end{tikzcd},
\begin{tikzcd}
M^{V_{1}} \arrow[d, "q_{j}"]\\
M                         
\end{tikzcd}
where $p_{i}$, $q_{j}$ are the evaluation maps (i.e. $p_{l}(\lambda) = \lambda(X_{l})$, and similarly for $q$'s). Therefore let $b_{0}$, $b_{1}$ respectively from $F(E^{V_{0}})\to M^{V_{0}}$, $F(E^{V_{1}})\to M^{V_{1}}$ be the unique bijections identifying these.

By the universal property of the product (a map into it is defined by its composites with all the projectors), to establish the square, it suffices to check that \[p_{l}\circ \Lambda\circ b_{1} = p_{l}\circ b_{0}\circ FS\tf\tag{1}\] for every $1 \leq l \leq n$. 

By definition chasing (recalling the relation of the index map $f$ to $\tf$) we see that\[
q_{f(l)}= p_{l} \circ \Lambda
\]
for every $l$. We claim also that $F\varepsilon_{f(l)} = F\pi_{l} \circ FS\tf$.

Indeed, it suffices to check $\varepsilon_{f(l)} = \pi_{l}\circ S\tf$ in $\syn(T)$. We have that $\varepsilon_{j}$ is defined by $x_{-j} = x_{1}$ by \cref*{dec-nfoldproduct}, and similarly $\pi_{l}$ is $x_{-l} = x_{1}$. So the composition on the right-hand side is $\pi_{l}\circ S\tf=$ \[
\exists T_{i}\rak{E^{V_{1}}(x_{-j})\wedge \bigwedge_i T_{i} = x_{-f(i)} \wedge T_{l} = x_{1}}
\]
which is provably equivalent to\[
\exists T_{i}\rak{ \bigwedge_i T_{i}=x_{-f(i)}\wedge T_{l} = x_{1}\wedge x_{-f(l)} = x_{1}}\]
p.e. to\[
\exists T_{i}\rak{ \bigwedge_i T_{i}=x_{-f(i)} \wedge x_{-f(l)} = x_{1}}
\]p.e. to \[
\exists T_{i}\rak{ \bigwedge_i T_{i}=x_{-f(i)}} \wedge x_{-f(l)} = x_{1}
\]
p.e. to $\varepsilon_{f(l)}$, as the first clause is vacuous. \qed

Now, we have  $p_{l}\circ \Lambda \circ b_{1} = q_{f(l)}\circ b_{1}$, which is equal by definition of $b_{1}$ to $F\varepsilon_{f(l)}$. On the other hand $p_{l}\circ b_{0} = F\pi_{l}$ by definition of $b_{0}$, so $p_{l}\circ b_{0} \circ FS\tf = F\pi_{l}\circ FS\tf = F\varepsilon_{f(l)}$ also. Thus we have checked (1).\qed\\

\subsection{Plan}
In terms of \cref{seg-semsetup}, $\calm$ will have underlying set $M$, and we will use $U = F(S)$ (which is a singleton as the validity $S$ is final in $\syn(T)$ and $F$ preserves final objects by condition 1). $\calm$ will also be such that the following holds for all $\varphi$:
\[\text{``$b_{\free(\varphi)}\circ F\iota_{\varphi}$ is a bijection $F(\varphi)\to \calm(\varphi)$''}  \tag{*\pr}
\]

where $b_{\free(\varphi)}$ denotes the bijection associated to $\free(\varphi)$ coming from \cref{dec-lastlem}, and $\iota_{\varphi}$ denotes the special mono $\varphi \to E^{\free(\varphi)}$. We consider $b_{\emptyset}$ to be the identity map of $U$ and $E^{\emptyset}$ to be $S$.

We note that $F$ preserves monos as a consequence of condition 1 (cf \cite{lu}, Lecture 4), so $F\iota_{\varphi}$ is always injective. Thus (*\pr) is equivalent to asserting that the image of $F(\varphi)$ under this composition is $\calm(\varphi)$.\\

Having obtained such a structure $\calm$, we will deduce that it is a model of $T$ as follows: if $\varphi$ is an element of $T$ (axiom), then $\varphi$ is a final object in $\syn(T)$, since it is isomorphic to the validity $S$. $F$ preserves final objects, so $F(\varphi)$ is a singleton. So by (*\pr), $\calm(\varphi)$ must be $U$. But according to \hyperlink{semanticdefs}{our definitions} this means $\calm \models \varphi$. \qed

\section{Definition}
If $|P| = n > 0$, let $\widetilde{P}$ be the formula $Px_{1}\dots x_{n}$. Let $s_{n}: M^{\{x_{1}, \dots x_{n}\}}\to M^{n}$ be the obvious bijection.\\

By \cref{dec-structurespecdby}, we define the structure $\calm$ by letting \hypertarget{defofmodel}{}\begin{itemize}
\item $P^{\calm} = (s_{n}\circ b_{\crul{x_{1}, \dots x_{n}}}\circ F\iota_{\widetilde{P}})(F(\widetilde{P}))$
\item for $P$ of arity 0, $P^{\calm} = F\iota_{P}(F(P))$.
\end{itemize}

\section{Induction}
We shall verify (*\pr) for all formulas $\varphi$ by induction. First we check it for atomics.
\subsection{Atomics}
\begin{itemize}
\item $x_{i} = x_{i}$: we have \begin{center}
\begin{tikzcd}
F(\varphi) \arrow[r, "F\iota_{\varphi}=\text{id}"] & F(E^{\crul{x_{i}}}) \arrow[r, "b_{\crul{x_{i}}}"] & M^{\crul{x_{i}}}\\
\end{tikzcd}
\end{center}
and both maps are bijective. By \hyperlink{suchthat-c}{\cref*{dec-missuchthat}(c)} (*\pr) therefore holds.\qstp

\item $x_{i} = x_{j}, i \neq j$: Let $\min = \min(i, j)$ and $\max = \max(i, j)$ and write
\begin{center}
\begin{tikzcd}
F(\varphi) \arrow[r, "F\iota_{\varphi}"] & F(E^{\crul{x_{i}, x_{j}}}) \arrow[r, "b_{\{x_{i}, x_{j}\}}"] \arrow[d, shift left, "F\pi_{\min}"] \arrow[d, shift right, "F\pi_{\max}"']& M^{\crul{x_{i}, x_{j}}} \arrow[ld, "p_{\min}"'] \arrow[ld, shift left=2, "p_{\max}"]\\
  & M
\end{tikzcd}
\end{center} where the inner and outer triangles (sharing a top leg) are commutative.
In $\syn(T)$, the composition $\pi_{\min}\circ \iota_{\varphi}$ is given by \cref{dec-soutgoing} as
\[
x_{-1}=x_{-2}\wedge \pi_{\min}
\] which is just\[
x_{-1}=x_{-2}\wedge x_{-1} = x_{1}. \tag{1}
\]
This is easily seen to be an isomorphism $\varphi \to x_{0} = x_{0}$, so $F(\pi_{\min}\circ \iota_{\varphi}) = F\pi_{\min}\circ F\iota_{\varphi} $ is bijective (an isomorphism in $\set$) merely as $F$ is a functor.

But by the inner triangle, we can write $F\pi_{\min}\circ F\iota_{\varphi} $ as $p_{\min}\circ b_{\crul{x_{i}, x_{j}}}\circ F\iota_{\varphi}$. The fact that this is a surjection to $M$ implies the image of $b_{\crul{x_{i}, x_{j}}}\circ F\iota_{\varphi}$ contains points with arbitrary $x_{\min}$-values. Therefore, if we can show that the compositions of this map with $p_{\min}$, $p_{\max}$ are equal, we will have that the image is exactly the diagonal in $M^{\crul{x_{i}, x_{j}}}$. This will be claim $(*\pr)$ for this case, by \hyperlink{suchthat-d}{\cref*{dec-missuchthat}(d)}.

But indeed, as noted $p_{\min}\circ b_{\crul{x_{i}, x_{j}}}\circ F\iota_{\varphi}$ is just $F\pi_{\min}\circ F\iota_{\varphi} $, and similarly for $\max$. Thus we need only show $\pi_{\min}\circ \iota_{\varphi} = \pi_{\max}\circ \iota_{\varphi}$ in $\syn(T)$.

But the formula (1) serves to define both of these; in the case of $\max$, we will obtain that the composition is defined by $x_{-1}=x_{-2}\wedge x_{-2} = x_{1}$, and this is provably equivalent. So the morphism is the same.\qstp

\item Predicate $P$ of arity 0: (*\pr) holds by \hyperlink{suchthat-b}{\cref*{dec-missuchthat}(b)} and \hyperlink{defofmodel}{the definition of $\calm$.}\qstp

\item $Px_{1}\dots x_{n}$: (*\pr) holds by \hyperlink{suchthat-a}{\cref*{dec-missuchthat}(a)} and the definition of $\calm$, taking preimages under $s_{n}$.\qstp
\end{itemize}

\bigskip
\bigskip

\begin{obs}\label{dec-theobs}Let $i_{1}, \dots i_{n}$ be integers. Then we have a surjection\[
\tf(x_{k}) = x_{i_{k}} \qquad f(k)= (\text{$i_{k}$'s rank among $i_{1}, \dots i_{n}$})
\]
for which the pullback square in \cref{dec-choosypully} is \begin{center}
\begin{tikzcd}
Px_{i_{1}}\dots x_{i_{n}}\wedge E^{V_{1}} \arrow[rr] \arrow[d, "\omega"]& & E^{V_{1}\arrow[d, "\chi"]}\\
Px_{1}\dots x_{n}\arrow[rr]& & E^{\crul{x_{1}, \dots x_{n}}}
\end{tikzcd}
\end{center}
where $\omega$, $\chi$ are the $\tf$-special morphisms. The first clause of $Px_{i_{1}}\dots x_{i_{n}}\wedge E^{V_{1}}$ already contains all the free variables ($V_{1}$ is just the set of $x_{i_{k}}$), so we can apply a special isomorphism to change the upper-left corner of this square to just $Px_{i_{1}}\dots x_{i_{n}}$, without changing anything else.

We can then apply $F$ to obtain a pullback square in $\set$, by condition 1. The right vertical of the result, is the left vertical of the commutative square from \cref{dec-lastlem}. Appending these and taking the outer rectangle yields a square in $\set$ as follows:\[
\begin{tikzcd}
F(Px_{i_{1}}\dots x_{i_{n}}) \arrow[rr, "b_{V_{1}}\circ F\iota_{1}"] \arrow[d, "F\omega"]& & M^{V_{1}\arrow[d, "\Lambda"]}\\
F(Px_{1}\dots x_{n})\arrow[rr, "b_{\crul{x_{1}, \dots x_{n}}}\circ F\iota_{0}"]& & M^{\crul{x_{1}, \dots x_{n}}}
\end{tikzcd}\tag{*}
\]
By the Pullback Lemma this is pullback, because the $\sett$ square we extended by from \cref{dec-lastlem} had isomorphisms on the top and bottom and was thus also pullback. To verify claim (*\pr) for general $Px_{i_{1}}\dots x_{i_{n}}$, we exactly need to show $\calm(Px_{i_{1}}\dots x_{i_{n}})$ is the image of the top row.\\ \end{obs}

\begin{itemize}
\item General $Px_{i_{1}}\dots x_{i_{n}}$: by \hyperlink{suchthat-a}{\cref*{dec-missuchthat}(a)} $\calm(Px_{i_{1}}\dots x_{i_{n}})$ is the subset of functions $h \in M^{\{x_{i_{1}, \dots x_{i_{n}}}\}}$ such that $(h(x_{i_{1}}), \dots h(x_{i_{n}}))\in P^{\calm}$; we claim this coincides with the image of the top row in (*).

Indeed, suppose $h \in M^{\{x_{i_{1}}, \dots x_{i_{n}}\}}$ is such. Then the unique $s_{n}$-lift of the point $(h(x_{i_{1}}), \dots h(x_{i_{n}}))\in P^{\calm}$ is the function $\widetilde{h}\in M^{\{x_{1}, \dots x_{n}\}}$ sending $x_{1}\mapsto h(x_{i_{1}})$, $\dots x_{n}\mapsto h(x_{i_{n}})$. So by the definition of $P^{\calm}$, this function $\widetilde{h}$ must be in $\calm (Px_{1}, \dots x_{n})$, which is the image of the bottom map in (*) by (*\pr) for $\widetilde{P}$. We note that as (*) is a pullback in $\set$, we have $\Lambda^{-1}(\text{im}(\text{bottom})) = \text{im}(\text{top})$ (this holds for the usual fiber product construction, and isomorphism over the diagram does not change either of these images). But by inspection $h$ is a lift of $\widetilde{h}$ under $\Lambda$, and thus $h$ must be in the image of the top row.

Conversely suppose $h \in \text{im}(\text{top})$. The function $\Lambda(h) = h\circ \tf$ returns $h(x_{i_{1}})$ on $x_{1}, \dots h(x_{i_{n}})$ on $x_{n}$, and this function must be in the image of the bottom map by commutativity. Therefore the pushdown of it by $s_{n}$ is in $P^{\calm}$ by definition, and this pushdown is just $(h(x_{i_{1}}), \dots h(x_{i_{n}}))$. But by \hyperlink{suchthat-a}{\cref*{dec-missuchthat}(a)} this means $h \in \calm(Px_{i_{1}}\dots x_{i_{n}})$.\qcl \qstp

\end{itemize}

\subsection{Inductive step}

\subsubsection{Case 1} First consider $\varphi = \alpha \vee \beta$ where $\alpha$ and $\beta$ are nonsentences. Let $V_{c} = \free(\varphi) = \free(\alpha) \cup \free(\beta)$. Clearly by \cref{seg-latticestruct}, we have that $\varphi \overset{\iota_{\varphi}}{\to} E^{V_{c}}$ is the join of the subobjects $\alpha \wedge E^{V_{c}}\overset{\iota_{1}}{\to} E^{V_{c}}$, $\beta \wedge E^{V_{c}}\overset{\iota_{2}}{\to} E^{V_{c}}$. As $F$ induces a map of Boolean algebras from $\sub(E^{V_{c}})$ to $ \sub(F(E^{V_{c}}))$, we conclude that 
\[
F(\varphi) \overset{F\iota_{\varphi}}{\rightarrow} F(E^{V_{c}})
\]
is the join of $F(\alpha \wedge E^{V_{c}})\overset{F\iota_{1}}{\to} F(E^{V_{c}})$ and $F(\beta \wedge E^{V_{c}})\overset{F\iota_{2}}{\to} F(E^{V_{c}})$.

But as we are in $\set$ now this is another way of stating that the image of $F\iota_{\varphi}$ is the union of the other two images.

We will show (using the inductive hypothesis and the combination of \cref{dec-choosypully} and \cref{dec-lastlem}) that under $b_{V_c}$ the first image goes to $\calm(\alpha)\times M^{\free(\beta)\setminus \free(\alpha)}$, and the second image to $\calm(\beta)\times M^{\free(\alpha)\setminus \free(\beta)}$. (Cf \hyperlink{suchthat-e}{\cref*{dec-missuchthat}(e)}.) This will show $\calm(\varphi)$ is equal to $b_{V_{c}}(F\iota_{\varphi}(F(\varphi)))$. This is the claim $(*\pr)$ for this case.

Indeed let us set $\tf: \free(\alpha) \to V_{c}$ to be inclusion, and apply \cref*{dec-choosypully}. The pullback square is\begin{center}
\begin{tikzcd}
\alpha \wedge E^{V_{c}}\arrow[rr, "\iota_{1}"] \arrow[d] & & E^{V_{c}}\arrow[d]\\
\alpha \arrow[rr, "\iota_{\alpha}"] & & E^{\free(\alpha)}
\end{tikzcd}
\end{center}

Applying $F$ and combining with \cref{dec-lastlem} as in \cref{dec-theobs}, we obtain \begin{center}
\begin{tikzcd}
F(\alpha\wedge E^{V_{c}}) \arrow[rr, "F\iota_{1} \circ b_{V_c}"]\arrow[d] & & M^{V_{c}}\arrow[d, "\Lambda"]\\
F(\alpha) \arrow[rr, "F\iota_{\alpha}\circ b_{\free(\alpha)}"]& & M^{\free(\alpha)}
\end{tikzcd}
\end{center}
As this is a pullback in $\set$, we have $\Lambda^{-1}(\text{im}(\text{bottom})) = \text{im}(\text{top})$. But by induction $\text{im}(\text{bottom})$ is just $\calm(\alpha)$. And $\Lambda^{-1}(\calm(\alpha))$ is just $\calm(\alpha)\times M^{\free(\beta) \setminus \free(\alpha)}$ by explicit inspection. So this is equal to the image of the top map, that is, to $b_{V_c}(F\iota_{1}(F(\alpha\wedge E^{V_{c}})))$, as desired. The verification for the second image is identical. $\qed_{\text{show}}$

\subsubsection{Case 2}
Suppose $\varphi = s \vee \beta$, $|s| = 0$, $|\beta| > 0$. We have by \hyperlink{suchthat-f}{\cref*{dec-missuchthat}(f)}  \[
\calm(\varphi) = \begin{cases}
\calm(\beta) & \calm(s) = \emptyset\\
M^{\free(\beta)} & \calm(s) = U.
\end{cases}
\]

Like in Case 1, $s \vee \beta \overset{\iota_{\varphi}}{\rightarrow} E^{\free (\beta)}$ is the join of 
 $s \wedge E^{\free (\beta)}\overset{\iota_{1}}{\rightarrow} E^{\free(\beta)}$ and $ \beta \overset{\iota_\beta}{\rightarrow} E^{\free(\beta)}$ as a subobject.
\renewcommand{\im}{\text{im}}
So we have \[F\iota_{\varphi}(F(s \vee \beta)) = F\iota_{1}(F(s \wedge E^{\free(\beta)}))\cup F\iota_\beta(F(\beta))\tag{1}\]
in \set.

First suppose $\calm(s) = \emptyset$. Then by induction $F(s) = \emptyset$. By Proposition 1.5, $\iota_{1}$ is a pullback of $s \to S$ in $\syn(T)$, and as $F$ preserves pullbacks this means $F\iota_{1} $ is a pullback of $F(s)\to U$ in \set, which is an empty map. Hence $F\iota_{1}$ is an empty map. So (1) becomes $F\iota_{\varphi}(F(s \vee \beta)) = F\iota_\beta(F(\beta))$, and the images of these under $b_{\free(\beta)}$ ($=b_{\free(\varphi)}$) are therefore the same. This shows the claim (*\pr), since $b_{\free(\beta)}(\RHS) =  \calm(\beta)$ by induction, and this is $\calm(\varphi)$.

On the other hand if $\calm (s) = U$ then by induction $F(s) \to U$ is a bijection. So $F\iota_{1}$ is a pullback of an {\it isomorphism} in \set, by the same reasoning as above. Therefore its image is all of $F(E^{\free(\beta)})$. (1) then gives that $F\iota_{\varphi}(F(s \vee \beta))$ is also equal to this. Applying $b_{\free(\beta)}$ takes this to all of $M^{\free(\beta)}$, so we have (*\pr).\qed\\

\subsubsection{Case 3}
Suppose $\varphi = s \vee s'$, $s$ and $s'$ both sentences. Then $s \vee s'$ is an evident subobject of $S$, and it is the join of $s$ and $s'$. So the image of $F(s \vee s')\to U$ is the union of images of $F(s) \to U$, $F(s')\to U$. By induction this is $\calm(s)\cup \calm(s')$. But $\calm(s \vee s')$ is indeed $\calm(s)\cup \calm(s')$ by \hyperlink{suchthat-g}{\cref*{dec-missuchthat}(g)}, so (*\pr) holds. \qed

\subsubsection{Case 4}
Supposing $\varphi = \lnot \alpha$, $\alpha$ a nonsentence: we have $\calm(\varphi) = \calm(\alpha)^{c}$, complement in $M^{\free(\alpha)}$. The objects $\neg \alpha $ and $ \neg \alpha \wedge E^{\free(\alpha)}$ represent the same subobject of $E^{\free(\alpha)}$, and the latter formula is the complement we have given in \cref{seg-latticestruct} for the subobject $\alpha$. As discussed in \cite{lu} Lecture 4, conditions 1 and 3 on $F$ together with the fact that $\sub(X)$ is a Boolean lattice and not just a join semilattice, force $F$ to preserve complements. Therefore $F\iota_{\varphi}(F(\varphi))$ is the complement of $F\iota_{\alpha}(F(\alpha))$ in $F(E^{\free(\alpha)})$, and the images under $b_{\free(\alpha)}$ are complements of each other in $M^{\free(\alpha)} = b_{\free(\alpha)}(F(E^{\free(\alpha)}))$. This is claim (*\pr) because by induction the image of $F\iota_{\alpha}(F(\alpha))$ is $\calm(\alpha)$.

Supposing $\varphi = \neg s$, $s$ a sentence: We have $\calm(\varphi) $ equal to the complement of $\calm(s)$ in $U$. By \cref{seg-subasentence}, $\varphi\to S$ is the complement of $s\to S$ in $\sub(S)$. Therefore 
the image of $F(\varphi)\to U$ is the complement of $\im(F(s)\to U)$ in $U$. By induction this image is $\calm(s)$, so we have (*\pr).\qed

\renewcommand{\im}{\text{Im }}
\newcommand{\inc}{\text{inclusion}}
\subsubsection{Case 5}
Suppose $\varphi = \exists X \psi$, and suppose that $X$ is in $\free(\psi) = \crul{X_{i}}$ and is not the only element of it. In this case we have an inclusion $\tf :\free(\varphi)\hookrightarrow \free(\psi)$ which gives a special morphism $S\tf: \psi \to \varphi$, defined by $\psi(x_{-i})\wedge\bigwedge_{j} x_{j} = x_{-f(j)}$. We claim $S\tf$ is an effective epimorphism.

We check this by checking (by \cref{dec-eesy} in \cref*{ch-stability}) that $\im S\tf$ is provably equivalent to $\varphi$. We have
\[
\im S\tf = \exists M_{i} \rak{\psi(M_{i})\wedge \bigwedge_{j} Z_{j} = M_{f(j)}} \tag{1}
\]
where $(M_{i})$ are unused and $Z_{j}$ are the free variables of $\varphi$. By \subs (1) is provably equivalent to\[
\exists M_{i} \rak{\psi(M_{i})\sbb{Z_{j}}{M_{f(j)}}\wedge \bigwedge_{j} Z_{j} = M_{f(j)}} 
\]\[
=\exists M_{i} \rak{\psi\sbb{M_{i}}{X_{i}}\sbb{Z_{j}}{M_{f(j)}}\wedge \bigwedge_{j} Z_{j} = M_{f(j)}} 
\]\[
=\exists M_{i} \rak{\psi\sbb{M_{I}}{X_{I}}\sbb{M_{f(j)}}{X_{f(j)}}\sbb{Z_{j}}{M_{f(j)}}\wedge \bigwedge_{j} Z_{j} = M_{f(j)}} \tag{1\pr}
\]
where $I$ is the unique index with $X_{I} = X$ (note if the input variables are distinct and unused, simultaneous substitution is just given by iterated substitution). This is equal to\[
\exists M_{i} \rak{\psi\sbb{M_{I}}{X}\sbb{Z_{j}}{X_{f(j)}}\wedge \bigwedge_{j} Z_{j} = M_{f(j)}} \tag{1\pr}
\]
and we may drop the substitution $\sbb{Z_{j}}{X_{f(j)}} = \sbb{Z_{j}}{\tf(Z_{j})}$ as it is trivial; $\tf$ is inclusion.

By \ex and \reord therefore (1\pr) is provably equivalent to \[
\exists M_{I}\rak{\psi\sbb{M_{I}}{X}\wedge \exists M_{f(j)}\rak{\bigwedge_{j} Z_{j} = M_{f(j)}}} \tag{2}
\]p.e. to \[
\exists M_{I}\rak{\psi\sbb{M_{I}}{X}}\wedge \exists M_{f(j)}\rak{\bigwedge_{j} Z_{j} = M_{f(j)}} \tag{2\pr}
\] p.e. to\[
\exists M_{I}\rak{\psi\sbb{M_{I}}{X}} \tag{2\pr\pr}
\]
as the last clause is vacuous. This is alpha-equivalent to $\varphi$.\qcl

Now by condition 2 on $F$, we have that $FS\tf: F(\psi) \to F(\varphi)$ is a surjection. From \cref{dec-lastlem} we have the (not pullback, just commutative) square\begin{center}
\begin{tikzcd}
F(\psi) \arrow[rr, "b_{\free(\psi)}"]\arrow[d, "FS\tf"] & & M^{\crul{X_{i}}}\arrow[d, "\lambda \mapsto \lambda \circ \inc \hspace{5mm} (\Lambda)"]\\
F(\varphi)\arrow[rr, "b_{\free(\varphi)}"] & & M^{\crul{X_{i}}\setminus X_{I}}
\end{tikzcd}
\end{center}
and by \hyperlink{suchthat-j}{\cref*{dec-missuchthat}(j)} $\calm(\varphi) = \Lambda(\calm(\psi))$. But by induction $\calm(\psi)$ is the image of the top map, and applying $\Lambda$ to this, we obtain exactly the image of the bottom map, by surjectivity of the left map. This shows (*\pr). \qed

\subsubsection{Case 6}
If $\varphi = \exists X \psi$ and $X$ {\it is} the only free variable of $\psi$, we also claim a type 2 effective epimorphism $m: \psi\to \varphi$. Indeed, existence of this map is obvious and its image as defined in \cref{seg-factsum} is just the formula $\varphi$, so it is e.e. \qcl

Therefore by condition 2, $Fm: F(\psi) \to F(\varphi)$ is a surjection. By \hyperlink{suchthat-k}{\cref*{dec-missuchthat}(k)} we have that $\calm(\varphi)$ is $U$ if $\calm(\psi)$ is nonempty, and $\emptyset$ otherwise.

In the first case, $F(\psi)$ is nonempty by induction. So $F(\varphi)$ must be nonempty in order to receive a map from it. Therefore the image in $U$ must be everything, and we have $(*\pr)$.

In the second case, $F(\psi)$ is empty by induction. So $F(\varphi)$ must be empty in order to receive a surjection from it. Therefore the image in $U$ must be empty, and we have $(*\pr)$.\qed

\subsubsection{Case 7}
Finally, suppose $\varphi = \exists X \psi$ and the quantified variable $X$ is not free in $\psi$. Then our proof theory easily gives that $\varphi$ is special-isomorphic to the formula $\widetilde{\varphi}:=\neg E \wedge \psi$.  Therefore, these are the same natural subobject of $E^{\free(\varphi)}$ ($=E^{\free(\widetilde{\varphi})}$). Using condition 2 we have $F\iota_{\varphi}(F(\varphi)) = F\iota_{\widetilde{\varphi}}(F(\widetilde{\varphi}))$.

The already checked induction steps suffice to show claim (*\pr) for the sentence $\neg E$; and hence by induction (using only the propositional steps) we have (*\pr) for $\widetilde{\varphi}$.

By \hyperlink{suchthat-l}{\cref*{dec-missuchthat}(l)}, $\calm(\varphi) = \calm(\psi)$; and by \cref*{dec-missuchthat}(f)-(i) (and some casework) we see that $\calm(\widetilde{\varphi})$ is also equal to $\calm(\psi)$. Therefore $\calm(\varphi) = \calm(\widetilde{\varphi})$.

But then we must have (*\pr) for $\varphi$ as well, since $b_{\free(\varphi)} = b_{\free(\widetilde{\varphi})}$ takes the above left-hand and right-hand images to the same thing. \qed

\section{Empty model case}
Now suppose the functor $F$ does in fact have $F(x_{0 }= x_{0}) = \emptyset$; thus we are back to the beginning of \cref{seg-modset}.  We claim in this case the theory $T\cup\crul{E}$ is consistent (in the weakened proof calculus we have used all along, of course).

Indeed, the claim is equivalent to saying $\neg E \not\cong S$, or $\neg E \cong \exists x (x = x)$ \underline{fails} to be final, in $\syn(T)$. If it were final, by preservation $F(\neg E)$ would be a singleton. So it suffices to show that this is not the case. We shall show that $F(\neg E)= \emptyset$.

Indeed, as noted in the Case 6 argument above, we have a type 2 e.e. $m: x_{0}=x_{0} \to \neg E$ in $\syn(T)$.

Therefore by condition 2 on $F$, $Fm$ is a surjection $F(x_{0} = x_{0})\to F(\neg E)$. But then $F(\neg E)$ must be $\emptyset$, as the domain is empty.\qcl

\section{Completeness}
This completes our derivation of the classical G\"{o}del Completeness Theorem. For if our theory $T$ is additionally consistent under the \underline{classical} (stronger) inference rules, then the theory $T' :=T\cup \crul{\neg E}$ is also classically consistent ($\neg E$ is a classical validity). Hence $T'$ is certainly a consistent theory in our weakened system. But $T' \cup \crul{E}$ is explicitly inconsistent, so by the above claim, taking a functor $F$ for $T\pr$ cannot yield the $M = \emptyset$ case. Therefore we construct a model $\calm$ for $T'$ which is fully classical (nonempty). It is also a model for $T$.\qed

\section{Ending remarks}
The attitude underlying this paper likely comes across at best as proof-relevant and at worst as aggressively formalist, constructivist, and anti-Platonic. My only defense is this: I am fascinated by the idea of finding structure in the apparently barren wasteland of syntax, through the algebraic {\it lens} of category theory. This is not what is achieved here. I can only hope that in observing the formal arguments, someone might be struck by genuine insight in this direction.

What we have done here, is shown that the constraints and possibilities within the psuedo-algebraic system which is first-order proof theory, admit a certain categorical delineation. By Deligne's theorem {\it this fact alone} implies a semantic conclusion.

Is it conceivable that in another situation, one might be able to define, on the basis of syntax, some categorical structure which did {\it not} come a posteriori from an ``intended interpretation''? Could one possibly thereby deduce a truly new semantic result? I find this question extremely interesting.

\printbibliography

\end{document}